\documentclass[10pt, letterpaper]{article}

\usepackage[english]{babel}
\usepackage{graphicx}
\usepackage{mathtools}
\usepackage{amsmath,nccmath}
\usepackage{amssymb}
\usepackage{breakcites}
\usepackage{arydshln}
\usepackage{cite}

\usepackage{comment}

\usepackage{amsthm}
\usepackage{hyperref}
\usepackage{subfigure}
\usepackage{tikz}
\usepackage{algorithm}

\usepackage{algorithmic}
\usepackage{pgfplots}
\usetikzlibrary{positioning,chains,fit,shapes,calc}
\usetikzlibrary{backgrounds}
\usetikzlibrary{arrows.meta}
\usetikzlibrary{shapes,snakes}
\pgfplotsset{compat=1.12}

\newtheorem{theorem}{Theorem}
\newtheorem{proposition}{Proposition}
\newtheorem{lemma}{Lemma}
\newtheorem{corollary}{Corollary}
\newtheorem{definition}{Definition}
\usepackage{lineno, blindtext, color}
\usepackage{cancel}

\newtheorem{remark}[theorem]{Remark}
\newcommand{\bigA}{\mbox{\normalfont\Large\bfseries $A$}}
\newcommand{\bigR}{\mbox{\normalfont\Large\bfseries $R$}}
\newcommand{\bigL}{\mbox{\normalfont\Large\bfseries $L$}}
\newcommand{\bigI}{\mbox{\normalfont\Large\bfseries $I$}}
\newcommand{\bigV}{\mbox{\normalfont\Large\bfseries $e$}}
\newcommand{\bigC}{\mbox{\normalfont\Large\bfseries $C$}}

\newcommand{\biglambdaI}{\mbox{\normalfont\Large\bfseries $\lambda I$}}
\newcommand{\bigzero}{\mbox{\normalfont\Large\bfseries $0$}}
\newcommand{\rvline}{\hspace*{-\arraycolsep}\vline\hspace*{-\arraycolsep}}

\definecolor{carmine}{rgb}{0.59, 0.0, 0.09}
\definecolor{asparagus}{rgb}{0.53, 0.66, 0.42}
\definecolor{ao}{rgb}{0.55, 0.71, 0.0}
\definecolor{bleudefrance}{rgb}{0.19, 0.55, 0.91}
\definecolor{dimgray}{rgb}{0.41, 0.41, 0.41}    
\definecolor{mediumorchid}{rgb}{0.73, 0.33, 0.83}
\definecolor{mediumtealblue}{rgb}{0.0, 0.33, 0.71}
\definecolor{harvestgold}{rgb}{0.85, 0.57, 0.0}
\definecolor{blue(pigment)}{rgb}{0.2, 0.2, 0.6}
\definecolor{forestgreen(traditional)}{rgb}{0.27, 0.35, 0.27}
\definecolor{cadmiumred}{rgb}{0.89, 0.0, 0.13}
\definecolor{orange(webcolor)}{rgb}{1.0, 0.5, 0.0}
\definecolor{c1}{rgb}{0.368417, 0.506779,0.709798}
\definecolor{c2}{rgb}{0.880722, 0.611041,0.142051}
\definecolor{c3}{rgb}{0.560181, 0.691569,0.194885}
\definecolor{c4}{rgb}{0.922526, 0.385626,0.209179}
\definecolor{c5}{rgb}{0.528488, 0.470624,0.701351}
\definecolor{c6}{rgb}{0.772079, 0.431554,0.102387}
\definecolor{c7}{rgb}{0.363898, 0.618501,0.782349}
\definecolor{c8}{rgb}{1, 0.75, 0}
\definecolor{c9}{rgb}{0.647624, 0.37816, 0.614037}
\definecolor{c10}{rgb}{0.571589, 0.586483, 0.}
\definecolor{c11}{rgb}{0.915, 0.3325, 0.2125}
\definecolor{c12}{rgb}{0.83, 0.46, 0.}
\definecolor{c13}{rgb}{0.9575, 0.545, 0.11475}
\definecolor{c14}{rgb}{1., 0.7575, 0.}
\definecolor{c15}{rgb}{0.6175, 0.715, 0.}
\definecolor{c16}{rgb}{0.15, 0.715, 0.595}
\definecolor{c17}{rgb}{0.3625, 0.545, 0.85}
\definecolor{c18}{rgb}{0.575, 0.4175, 0.85}
\definecolor{c19}{rgb}{0.677, 0.358, 0.595}
\definecolor{c20}{rgb}{0.7875, 0.358, 0.425}
\definecolor{c21}{rgb}{0.915, 0.3325, 0.2125}

\title{Distributed design of deterministic discrete-time privacy preserving average consensus for multi-agent systems through network augmentation\footnote{
{\scriptsize G. Ramos and A. Pedro Aguiar are with the Department of Electrical and Computer Engineering, Faculty of Engineering, University of Porto, Portugal. 
S. Kar is with the Department of Electrical and Computer Engineering, Carnegie Mellon University, USA. 
S. Pequito is a faculty member at the Delft University of Technology in the Delft Center for Systems and Control.
 This work was supported in part by projects: IMPROVE (POCI-01-0145-FEDER-031823), funded by FEDER funds through COMPETE2020 – POCI and by the Portuguese national FCT/MCTES (PIDDAC); RELIABLE (PTDC/EEI-AUT/3522/2020) funded by FCT/MCTES; and DynamiCITY (NORTE-01-0145-FEDER-000073), funded by NORTE2020/PORTUGAL2020, through the European Regional Development Fund.}}
}
\date{\today}

\author{
Guilherme Ramos, A. Pedro Aguiar, Soummya Kar, Sérgio Pequito}

\begin{document}
\maketitle

\begin{abstract}
Average consensus protocols emerge with a central role in distributed systems and decision-making such as distributed information fusion, distributed optimization, distributed estimation, and control. 
A key advantage of these protocols is that agents exchange and reveal their state information only to their neighbors. Yet, it can raise privacy concerns in situations where the agents' states contain sensitive information. 
%, which may be retrieved using estimators that consider both the access to an incoming stream of data from the neighbors and the knowledge of the consensus protocol performed by all the agents in the network. 
In this paper, we propose a novel (noiseless) privacy preserving distributed algorithms for multi-agent systems to reach an average consensus. 
The main idea of the algorithms is that each agent runs a (small) network with a crafted structure and dynamics to form a network of networks (i.e., the connection between the newly created networks and their interconnections respecting the initial network connections). Together with a re-weighting of the dynamic parameters dictating the inter-agent dynamics and the initial states, we show that it is possible to ensure that the value of each node converges to the consensus value of the original network. 
%each agent augments its state, which initial augmented state sums to the agent's initial state, and the augmented state dynamics is crafted to ensure that no agent in the network can retrieve all the augmented states of the remaining agents. 
Furthermore, we show that, under mild assumptions, it is possible to craft the dynamics such that the design can be achieved in a distributed fashion.
Finally, we illustrate the proposed algorithm with examples. 
\end{abstract}

% Also, regarding the augmentation of the state, why not something like each node becomes a network such that the network of networks (i.e., the connection between the newly created networks and their interconnections respecting the initial network connections), will be such that there exists a re-weighting of the dynamic parameters and the initial states such that the value of each node converges to the same values as the original network... or something like that...

% Note that keywords are not normally used for peerreview papers.
% \begin{keywords}
% privacy, observability, multi-agent systems, average consensus.
% \end{keywords}

% make the title area
\maketitle

% To allow for easy dual compilation without having to reenter the
% abstract/keywords data, the \IEEEtitleabstractindextext text will
% not be used in maketitle, but will appear (i.e., to be "transported")
% here as \IEEEdisplaynontitleabstractindextext when the compsoc 
% or transmag modes are not selected <OR> if conference mode is selected 
% - because all conference papers position the abstract like regular
% papers do.
%\IEEEdisplaynontitleabstractindextext
% \IEEEdisplaynontitleabstractindextext has no effect when using
% compsoc or transmag under a non-conference mode.

% For peer review papers, you can put extra information on the cover
% page as needed:
% \ifCLASSOPTIONpeerreview
% \begin{center} \bfseries EDICS Category: 3-BBND \end{center}
% \fi
%
% For peerreview papers, this IEEEtran command inserts a page break and
% creates the second title. It will be ignored for other modes.
%\IEEEpeerreviewmaketitle

% The very first letter is a 2 line initial drop letter followed
% by the rest of the first word in caps.
% 
% form to use if the first word consists of a single letter:
% \IEEEPARstart{A}{demo} file is ....
% 
% form to use if you need the single drop letter followed by
% normal text (unknown if ever used by the IEEE):
% \IEEEPARstart{A}{}demo file is ....
% 
% Some journals put the first two words in caps:
% \IEEEPARstart{T}{his demo} file is ....
% 
% Here we have the typical use of a "T" for an initial drop letter
% and "HIS" in caps to complete the first word.
\section{Introduction}\label{sec:introduction} 

%\IEEEPARstart{M}{ulti-agent} 
% Multi-agent systems (MAS) can address problems that may be challenging or inadequate for solving with a single agent or a monolithic system~\cite{dorri2018multi}. 
% MAS model a plethora of applications, such as consensus problems~\cite{9304107,ramos2020ijoc}, target surveillance~\cite{hu2019distributed}, network resistance~\cite{RAMOS2021104842}, online trading~\cite{luo2018distributed}, disaster response~\cite{nadi2017adaptive}, and wireless sensor networks (WSN)~\cite{akyildiz2002wireless}.

We are facing a growing interest in distributed computation, reflected by the increase of data being exchanged in large-scale systems that can be spatially distributed. 
Therefore, numerous systems have been developed to compute distributedly a function. The distributed function results from each agent computing a local function, using the neighbors' communicated information, and sharing the output with its neighbors, for instance, using consensus~\cite{bullo2009distributed,9304107,ramos2020ijoc,RAMOS2021104842} or state retrieval~\cite{sundaram2008distributed}. 
Under the described setting, for security reasons, the agents may want to ensure that their information is not disclosed. 
This privacy requirement leads to what is called \emph{secure multi-party computation}, the problem of a set of agents computing an agreed-upon function of their initial states in a private way~\cite{kolesnikov2009advances}. 

Computer science literature often focus on the problem of shielding sensitive information. 
Techniques using cryptographic schemes provide ways of ensuring privacy-preserving information, such as Yao’s garbled circuit~\cite{yao1982protocols}, Shamir’s secret sharing~\cite{shamir1979share}, secure multi-party computation~\cite{chaum1988multiparty}, and homomorphic encryption (HE)~\cite{smart2016cryptography}.  
In our work, we are interested in ensuring data privacy to reach an average consensus. 
There are mainly three research directions that aim tho reach this goal: 
\begin{itemize}
    \item[$(i)$] homomorphic encryption-based (HE-based);
    \item[$(ii)$] differential privacy-based (DP-based);
    \item[$(iii)$] observability-based (O-based). 
\end{itemize}

The authors of~\cite{lazzeretti2014secure,freris2016distributed,kishida2018encrypted,yin2019accurate,ruan2019secure,hadjicostis2020privacy} introduced HE-based average consensus methods. 
%In~\cite{shoukry2016privacy,zhang2018enabling}, the authors designed HE-based distributed optimization algorithms.  
Notwithstanding, a major drawback of the HE-based approaches is that they require expensive computation and communication, making them (potentially) unsustainable to use in real-world applications that have computation power and communication limitations~\cite{lagendijk2012encrypted,kogiso2015cyber}. %, e.g., controllability of drones or power grids. 

A different research path to ensure that data remains private is DP-based, which consists of leveraging differential privacy guarantees~\cite{cortes2016differential}.  
The main idea behind differential privacy is to add noise to the shared information. 
In the works of~\cite{huang2012differentially,nozari2017differentially,wang2018privacy,gao2018differentially,fiore2019resilient}, average consensus algorithms are devised to ensure differential privacy, aiming to converge to the average, while keeping the agents' initial values private. 
However, the introduced noise does not ensure the exact average consensus. 
To circumvent this problem, in~\cite{mo2016privacy,nozari2017differentially,he2019consensus,he2018preserving}, the authors propose to add correlated noise. 
However, from an implementation point of view, it is difficult to guarantee that the required assumptions hold, and therefore there will always be uncertainty in the final consensus value achieved. 
% In contrast this line of research, the algorithms that we propose do not require %the assumption that the %network is undirected neither it needs the 
% the generation of noise to be added to the agent's states. 
Furthermore, we observe that noise generation is typically achieved with a pseudo-random generation created, deterministically, from a seed. 
Therefore, the privacy guarantees depend on the seed that must be hidden from the other agents (i.e., personal or secret) or the use of an expensive random number generator device~\cite{yu2019survey}.

In contrast, it is possible to hide information without adding noise into the agents' states. 
The O-based approaches concept means there is a stream of observation/measurements of neighboring agents and the agent's own, generated by the underlying consensus protocol (assume to be known to all the agents). 
Thus, the notion of observability in dynamical systems describes necessary and sufficient conditions to obtain an algorithm (i.e., an estimator) to retrieve all agents' states~\cite{boutat2021observability}. 
We can categorize the methods of this line of research according to what it is proposed to be modified in the consensus dynamics, e.g.: 
\begin{itemize}
    \item changing the dynamics, while maintaining the structure (i.e., the network topology);
    \item editing the connections (i.e., adding/removing edges in the network and the associated  weights);
    \item augment the network by including more nodes/(i.e., agents) and connections between these (and, subsequently, the number of parameters specifying the dynamics);
    \item adjusting the initial agents' states;
    \item a combination of the previous.
\end{itemize}

In~\cite{gupta2016confidentiality}, the authors propose a distributed average information consensus algorithm that preserves each agent initial state confidentiality, without corrupting the state values of agents with noise. 
The authors ensure privacy with the help of \textit{concealing factors}, which are assigned to the agents in the network before the initiation of the proposed consensus algorithm by a central authority (CA). 
Further, the method requires a balancing constraint, implying that the sum of the edges' weights that end in a node equals the sum of the edges' weights that start in that node. 
In this work, we do not need a central authority, and we further do not require a balancing constraint.  

%zero-mean noise as practiced in differentially private average consensus algorithms and hence, evading the tradeoff between privacy and accuracy 
In~\cite{he2018preserving}, the authors  %present a theoretical framework to compute the optimal distributed estimation and the privacy analysis. With this framework, the authors 
obtained closed-form expressions of both the optimal distributed estimation and the privacy parameters. 
%For example, in~\cite{gupta2016confidentiality,he2018preserving}, the authors propose to scale the agents' information of each agent privately.  
The work in~\cite{wang2019privacy} proposes a privacy-preserving approach, based on state decomposition, for the network average consensus problem, where each node decomposes its state into two sub-states with random initial values with a mean equal to the original state initial value.   
In~\cite{ridgley2020private}, the authors proposed a dynamic average consensus algorithm that guarantees accuracy and preserves privacy of the initial values, under topological restrictions on the communication digraph. A main drawback is that the algorithm has to create a virtual network with $\mathcal O(n^2)$ nodes, where $n$ is the number of agents in the original consensus network. This contrasts with the methods that we propose here, which require creating a virtual network with $\mathcal O(n)$ nodes. 
%In~\cite{wang2019privacy,ridgley2020private}, each agent uses two states, one to be the real state and the other to be a dummy state to be the one shared with neighbor agents. 
%By following this strategy, the agents can achieve average consensus. 

The work in~\cite{pequito2014design} analyzes the interplay between a specific network topologies and the observability subspace generated by the underlying dynamics. The authors propose a communication protocol to reach weighted average consensus on the initial network agents' states, but where each agent is only able to retrieve a small subset of the initial states of the agents that are not its neighbors. 
Finally, in~\cite{alaeddini2017adaptive}, the authors use the concept of observability together with tools from graph theory and optimization to present an algorithm for network synthesis with privacy guarantees. The method selects the optimal weights for the communication graph to maximize the privacy of the network nodes. 
Notwithstanding, in the previous methods, the design is difficult to accomplish, and it is not guaranteed that \textit{all} agents manage to keep their state private. 

% Finally, the works in~\cite{pequito2014design,alaeddini2017adaptive} develop observability-based methods to secure privacy, designing the communication network to reduce the agents' capacity of recovering other agents' states. 

Our work aligns with an observability-based approach to achieve privacy in average consensus while ensuring the convergence to the \emph{exact} average of the initial agents' values, and that \emph{no agent is able to retrieve other agents' initial state}, which can be understood as full privacy. 
Moreover, contrarily to~\cite{ridgley2020private}, our nodes can share information with their neighbors in a broadcast setting.  
%Here, we allow that each agent in the network to observe the messages that other agents send except the local copies of each agent that further do not communicate with any other agent. 
In fact, in our scheme, we propose that an agent shares one of its states, and keep the remaining hidden from the network despite the local dynamics being publicly available (i.e., available to \textit{all} the other agents). Therefore, in a sense, the secret key to each agent is the way that it splits its initial condition across its augmented states, which can be done arbitrarily as long as some of its states are not set to zero.

\paragraph*{Main contributions}
We propose two network designs that consist in augmenting the agents' state and impose a given dynamics on it that allows a network of multi-agents to reach average consensus, while keeping the agents' initial states private. 
%By augmenting the state of an agent with a set of copy nodes that only interact locally, each agent can decide how to distribute its initial state to the copies while ensuring that the other agents cannot recover the initial information. 
Specifically, the first network design that we propose achieves the objective by increasing the agents' univariate initial state to four states (i.e., the original one plus three more), and re-defining the obtained network's weights, including the original weights of the consensus network, such that consensus is possible while ensuring privacy. 
Furthermore, it is possible to proceed to a rescaling of the initial condition through the left-eigenvector associated with the eigenvalue one of the matrix described by the consensus protocol such that \emph{average} consensus is possible. 

Notwithstanding, the computation of such eigenvectors can only be approximated to a certain precision-level with polynomial-time algorithms. 
%designing four local copies of each agent and the weights of the adjacency matrix. However, one of the steps required is to compute the left-eigenvector associated with the eigenvalue $1$ of the dynamics matrix, which in general may be only approximately computed with polynomial-time algorithms. 
Therefore, we propose a secondly (mildly) constrained design, which requires the network to be bidirectional and time-reversible. 
We emphasize that it does not need to be symmetric (in terms of the weights). 
In this scenario, we need to increase the agents' univariate initial state to five states (i.e., the original one plus four more). 
This new construction allows, withal, computing exactly the \mbox{left-eigenvector} associated with the eigenvalue one of the described by the consensus protocol matrix with linear-time complexity and in a distributed fashion. 

% Say that we propose a network with five states for the agents, and then say that under the mentioned assumptions we can compute the left eig... associated with the...  in linear-time complexity. Furthermore, this can be done in a distributed fashion, right?

%algorithm that requires the network to be bidirected. 
%The second algorithm also designs four local copies of each agent and the weights of the adjacency matrix, but it further ensures that the final network is time-reversible. This additional property allows computing the left-eigenvector associated with the eigenvalue 1 of the dynamics matrix explicitly with linear-time complexity. 

\paragraph*{Paper structure}

In Section~\ref{sec:notation}, we present the notation used in this paper. 
In Section~\ref{sec:prob_stat}, we formally state the two problem formulations addressed in this paper. 
In Section~\ref{sec:main_res}, we present solutions to the proposed problems. Lastly, in Section~\ref{sec:ill_exp}, we present several illustrative simulations that support the main results in the current paper. 
%We start by presenting an algorithm for the design of the dynamics matrix, that creates two local copies of each agent, which ensures privacy when the MAS has not the goal of reaching average consensus. 
%After, we show that the previous solution cannot be deployed when we had the constraint goal of reaching average consensus. 

\section{Notation}\label{sec:notation}
We use lower case letters do denote vectors (e.g., $x$) and upper case letters to denote matrices (e.g., $A$ and $C$). 
We use $x_k$ to denote the $k$th entry of vector $x\in\mathbb R^N$, with $k\in\{1,\ldots, N\}$, $A_i$ to denote the $i$th row of matrix $A\in\mathbb R^{N\times M}$ and $A_{ij}$ to denote the $j$th entry of the $i$th row of $A$, with $i\in\{1,\ldots, N\}$ and $j\in\{1,\ldots, M\}$. 

Given a matrix $A\in\mathbb R^{N\times N}$, we denote by $span(A)$ the \textit{linear span} of $A$, i.e., the smallest linear subspace that contains the set of column vectors of matrix $A$. 
Moreover, we denote by $e^N_i$ the $i$-th canonical $N$-dimensional column vector, i.e., a vector of size $N$ with all entries equal to zero, except the $i$-th entry that is one. We denote by $I_N$ the $N\times N$ identity matrix, and 
%if $\mathcal I\subset\{1,.\ldots,N\}$ then we denote 
by $I_N^{\mathcal I}$ the $N\times N$ diagonal matrix such that ${(I_N^{\mathcal I})}_{ii}=1$ if $i\in\mathcal I\subset\{1,.\ldots,N\}$ and ${(I_N^{\mathcal I})}_{ii}=0$ otherwise. 

We represent a (directed) \emph{network of agents} by the graph $\mathcal G=\langle \mathcal X,\mathcal E_{\mathcal X,\mathcal X}\rangle$, where $\mathcal X=\{1,\ldots,n\}$ are the \emph{nodes} that denote the set of $n>2$ agents, and $\mathcal E_{\mathcal X,\mathcal X}\subset \mathcal X\times \mathcal X$ are the \emph{edges} that correspond to pairs of agents (nodes).  
If $(i,j) \in \mathcal E_{\mathcal X,\mathcal X}$ then $i$ transmits to $j$, and we can associate a weight $w(i,j)=w_{ij}$ given by a weight function $w:\mathcal E_{\mathcal X,\mathcal X}\to\mathbb R^+$. %that communicate with each other, together with a nonzero real number that denotes the weight of the edge. 
In the case where the edges do not have associated a weight, we simply denote them as pairs of nodes. 
If a digraph $\mathcal G=\langle \mathcal X,\mathcal E_{\mathcal X,\mathcal X}\rangle$, where $\mathcal X=\{1,\ldots,n\}$ does not have weights associated with edges is such that $(i,j)\in\mathcal E_{\mathcal X,\mathcal X}$ \textit{if and only if} $(j,i)\in\mathcal E_{\mathcal X,\mathcal X}$, then we call it a \emph{graph}. 
Similarly, if a digraph $\mathcal G=\langle \mathcal X,\mathcal E_{\mathcal X,\mathcal X}\rangle$, where $\mathcal X=\{1,\ldots,n\}$ has weights associated with edges is such that $((i,j),w_{i,j})\in\mathcal E_{\mathcal X,\mathcal X}$ \textit{if and only if} $((j,i),w_{j,i})\in\mathcal E_{\mathcal X,\mathcal X}$ (with possibly $w_{i,j} \neq w_{j,i}$), then we call it a \emph{bidirected graph}. 

Additionally, a matrix $A\in\mathbb R^{N\times N}$, can induce a (weighted) network of agents through a digraph representation $\mathcal G(A)=\langle \mathcal X,\mathcal E_{\mathcal X,\mathcal X}\rangle$, where $\mathcal X=\{1,\ldots,N\}$ and $(i,j)\in\mathcal E_{\mathcal X,\mathcal X}$ (or $((i,j),w)\in\mathcal E_{\mathcal X,\mathcal X}$) if and only if $A_{j,i}\neq 0$ (or $A_{j,i}=w\neq 0$). 
Analogously, given a network of agents $\mathcal G=\langle \mathcal X,\mathcal E_{\mathcal X,\mathcal X}\rangle$, we define the matrix $A\equiv A(\mathcal G)\in\mathbb R^{n\times n}$ such that $A_{j,i}=w_{i,j}$ if $((i,j),w_{i,j})\in\mathcal E_{\mathcal X,\mathcal X}$ (or in the unweighted case $A_{j,i}=1$ if $(i,j)\in\mathcal E_{\mathcal X,\mathcal X}$) and $A_{j,i}=0$ otherwise. 
When only the structure of the network is defined, the edges, but not the weights, we highlight this fact using a bar on top of the matrix (e.g., $\bar A$).

For a network of agents (communication graph) $\mathcal G=\langle \mathcal X,\mathcal E_{\mathcal X,\mathcal X}\rangle$ and an agent $i\in\mathcal X$ in the network of agents $\mathcal G$, we denote the \emph{neighborhood} of agent $i$ by $\mathcal N_i$, where $\mathcal N_i=\{i\}\cup\{j\,:\,(i,j)\in\mathcal E_{\mathcal X,\mathcal X}\}$. 
Finally, we say that a network is \emph{strongly connected} if there is a path between each pair of nodes (i.e., if for each node $x\in\mathcal X$ there is a sequence of nodes $x,x_1,\ldots,x_k,y$ for all $y\in\mathcal X$ such that $(x,x_1),(x_k,y)$, $(x_i,x_{i+1})\in\mathcal E_{\mathcal X,\mathcal X}$ for all $i=1,\ldots,k-1$).

\section{Problem statement}\label{sec:prob_stat}

In this paper, we consider multi-agent systems whose goal is to reach consensus, through discrete-time linear time-invariant (DT-LTI) dynamics. 
In the context of multi-agents, the agents' interdependencies are described by a network, given by a digraph $\mathcal G=\langle\mathcal X,\mathcal E_{\mathcal X,\mathcal X}\rangle$, where each node corresponds to an agent state. 
First, notice that every agent $i$ is such that its dynamics are driven by input-output relations due to their state's exchange. Specifically, the agent's state $x_i[k]\in \mathbb{R}$ at a given time $k\in\mathbb N $ evolves as follows: 
\[
x_i[k+1]=w_{ii}x_i[k] +u_i[k]
\]
where $u_i\in \mathbb{R}$ is the (external) input to agent $i$, which, in the context of the present paper, is given as a linear combination of the agent's incoming neighbors, i.e.,
\begin{equation}\label{eq:output}
u_i[k]=\sum_{j\in\mathcal N_i\setminus\{i\}} w_jx_j[k].
\end{equation}
From the perspective of agent $i$, it measures the state of agents $j$, for $j\in\mathcal N_i$. 
Therefore, we can describe the measured output $y_i[k]\in\mathbb{R}^N$ as
\[
y_i[k]=I_N^{\mathcal N_i}x[k],
\]
where $I_N^{\mathcal N_i}$ denotes the $N\times N$ matrix with columns given by the canonical vector $e_N^j$, for $j\in\mathcal N_i$, and $x[k]$ is the vector gathering the agents' states. 
Henceforth, the external input to an agent $i$ is a closed-loop dynamical system %with state feedback 
described as
\[
x_i[k+1]=w_i y_i[k]=w_ix[k],
\]
 where $w_i$ is the combined set of weights. %(or gain in the context of state feedback). 
 %, and equivalent to equation 
 
Next, we explain how privacy can be understood as an observability-based (i.e., O-based) problem where the interdependencies between agents conditioning which states can be retrieved. %, i.e., are not private.
Therefore, we can now pose the question: 
``Can agent i
retrieve agent $j$’s state $x_j[0]$ for all $j = 1,\ldots,N$?'' 
If the answer is negative across all agents (i.e., for all $i$ and for all $j$), then we can say that all states are private. 
Specifically, we notice that the dynamics is a consequence of input/output relations between the agents. Thus, to assess how many states an agent has access to  over the entire dynamical network, we have 
\begin{equation}\label{eq:lti}
\begin{array}{c}
 \begin{array}{rcl}
    x[k+1] & = &  A(G) x[k],\\[0.2cm] 
    y_i[k] & = &  I_N^{\mathcal N_i} x[k],\\[0.2cm]
 \end{array}
 \\[0.3cm]
  x[0]\in\mathbb R^N,
  \end{array}
\end{equation}
where, for $k\in\mathbb Z_0^+$, $x[k]\in\mathbb{R}^{N}$ is the state, and $y_i[k]\in\mathbb R^{N\times N}$ is the output for agent $i$. 
%The system models the scenario that each agent $i$ can directly observe its neighbors (and itself), as indicated by $\mathbb I_n^{\mathcal N_i}$. 

% Therefore, we can now pose the question: 
% ``\textit{Can agent $i$ retrieve agent $j$'s state $x_0^j$ for all $j=1,...n$ ?}'' 
% If the answer to this question is negative across all agents (for all $i$) and for all $j$, then we can say that all states are private. 
% It turns out that this is related to the question of observability of linear time-invariant dynamical systems. 
% In fact, it has been shown that the answer to the question is possible, an algorithm can be placed (often referred to as observer) that uses the measured output and the model to retrieve the initial state. 
% Notice that, because the LTI is deterministic (i.e., there is no noise), knowing the initial state implies that all the states, at any time, are retrievable. 

In particular, we focus on the problem of consensus. It is the decision-making process that strives to achieve agreement among members of a group. 
The consensus process is restricted by the communication between agents, delimited by the network topology. 
Hence, the consensus value of any agent is the outcome of a distributed algorithm. 
This problem is of distinguished importance due to its wide range of applications. 
Examples of applications are such as distributed optimization~\cite{johansson2008subgradient,tsitsiklis1986distributed}; motion coordination tasks (e.g.,  flocking, leader following)~\cite{jadbabaie2003coordination,liu2015distributed}; rendezvous problems~\cite{cortes2006robust}; computer networks resource allocation~\cite{chiang2007layering}; and even in computing relative importance of webpages in the PageRank algorithm~\cite{silvestre2018pagerank}.

To determine the states that are retrievable when the system dynamics and measurements are available, we can invoke a system's property so-called observabilility. A system is \textit{observable} if, for any possible evolution of state and control vectors, the current system's state can be estimated using only the information from the outputs. 
In the context of dynamical systems, it is possible to retrieve the agents' states from the model using estimators/observers (e.g., in the case of deterministic systems, we could use~\cite{luenberger1971introduction}) by using the measured output (i.e., the neighbors' states) and the knowledge of the system's dynamics. Therefore, \textit{we can view privacy as an observability-based problem}. 

Several criteria exist to assess if a system is observable~\cite{antsaklis1997linear}. 
In the context of this paper, we relay on the Popov-Belevitch-Hautus eigenvector criteria (PBH test) defined next.

%Then I would say that this question has been solved through a system's property so-called observabilility (bring the formal definition of observability). Then, mention that several criteria exist, and in the context of this paper we relay on the PBH defined next:
\vspace{1mm} 
{\centering
\allowdisplaybreaks
\noindent \colorbox{c7!15}{\parbox{0.98\textwidth}{
\begin{proposition}[The Popov-Belevitch-Hautus (PBH test) eigenvector observability test~\cite{antsaklis1997linear}]\label{th:PBH}
 Consider the LTI system $x[k+1]= Ax[k]$ and $y[k]=Cx[k]$, where $C$ is the output matrix. 
 We can refer to the system by the matrices that describe them, i.e., the pair (A,C).  
 Given $(A,C)$, with $A\in\mathbb R^{N\times N}$ and $C\in\mathbb R^{M\times N}$, consider the matrix
 \[
     P_{O(A,C)}^\lambda \equiv
     \begin{bmatrix}
        C\\
        \lambda I_N - A
     \end{bmatrix},\,\forall \lambda\in\mathbb C.
 \]
 The pair $(A,C)$ is observable if and only if  $rank\left(P_{O(A,C)}^\lambda\right)=N$, $\forall \lambda\in\mathbb C$.\hfill$\circ$
\end{proposition}
}
}
}
\vspace{1mm}

% {\color{red}
% Def - privacy as an observability-based problem... write it as an unobservabiltiy problem... a state of an agent is private if it belongs to the unobservability subspace of all agents except itself.
% }
Next, we can define privacy through an observability criterion. 

%we define privacy as follows.
\vspace{1mm} 
{\centering
\allowdisplaybreaks
\noindent \colorbox{c2!18}{\parbox{0.98\textwidth}{
\begin{definition}\label{def:priv}
Given a communication digraph $\mathcal G=(\mathcal V,\mathcal E)$, $A(\mathcal G)\in\mathbb R^{N\times N}$ and an agent $i\in\mathcal V$, we say that the initial state of agent $j\in\mathcal V\setminus\{i\}$ is private to agent $i$ (or, equivalently, agent $i$ \textbf{cannot recover} the initial state of any agent $j$) if for $\lambda\in\mathbb C$
\[
e_j^N\notin span\left(P_{O\left(A,I_N^{\mathcal N_i}\right)}^\lambda\right).
\]\hfill$\circ$
%where $\mathcal N_j$ is the set of neig. 
\end{definition}
}
}
}
\vspace{1mm}

In other words, Definition~\ref{def:priv} reflects that if an agent, by receiving information from its neighbors (observing its neighbors' states),  can recover another agent's sate, then this other agent's state is not private to the first agent.  
Remark that, from this, it follows that the neighbors can always have access to the initial state of the neighbors, and, therefore, privacy in the strict sense above cannot be achieved. Thus, a new approach to the problem of ensuring privacy needs to be considered, as we detail hereafter.

Notwithstanding, it might be possible to augment the system dynamics to circumvent such implicit limitation, which is motivated by the following example.

% maybe make this assumption explicit as well.... also motivate the assumption from practical applications where it is not possible to change the interconnections naturally, due to short range communications between agents, or inherent infrastructure that is to costly to change.

%%%%%%%%%%%%%%%%%%%%%%%%%%%%%%%%%%%%%%%%%%%%%%%%%%%%%%%%%%%%%%%%
%%%%%%%%%%%%%%%%%%%%%%%%%%%%%%%%%%%%%%%%%%%%%%%%%%%%%%%%%%%%%%%%
%%%%%%%%%%%%%%%%%%%%%%%%%%%%%%%%%%%%%%%%%%%%%%%%%%%%%%%%%%%%%%%%
%%%%%%%%%%%%%%%%%%%%%%%%%%%%%%%%%%%%%%%%%%%%%%%%%%%%%%%%%%%%%%%%
%%%%% ILLUSTRATIVE EXAMPLE TO MOTIVATE PROBLEM STATEMENT %%%%%%%
%%%%%%%%%%%%%%%%%%%%%%%%%%%%%%%%%%%%%%%%%%%%%%%%%%%%%%%%%%%%%%%%
%%%%%%%%%%%%%%%%%%%%%%%%%%%%%%%%%%%%%%%%%%%%%%%%%%%%%%%%%%%%%%%%

\begin{figure}[H]
\centering
\subfigure[]{
\begin{tikzpicture}[scale=.56, transform shape,node distance=1.5cm]
\begin{scope}[every node/.style={circle,thick,draw},square/.style={regular polygon,regular polygon sides=4}]
\node[fill=green!20] (1) at (2.15387,2.50051) {\large \,};
\node[fill=green!20] (2) at (4.33095,3.72641) {\Large $x_1$};
\node[fill=green!20] (3) at (2.13179,0.) {\large \,};
\node[fill=green!20] (4) at (0.,3.77145) {\large \,};
\end{scope}
\begin{scope}[>={Stealth[black]},
              every edge/.style={draw=black, thick}]
\path [->] (1) edge[bend right=15] node {} (2);
\path [->] (2) edge[bend right=15] node {} (1);
\path [->] (1) edge[bend right=15] node {} (3);
\path [->] (3) edge[bend right=15] node {} (1);
\path [->] (1) edge[bend right=15] node {} (4);
\path [->] (4) edge[bend right=15] node {} (1);
\end{scope}
\end{tikzpicture}
}
\subfigure[]{
\begin{tikzpicture}[scale=.56, transform shape,node distance=1.5cm]
\begin{scope}[every node/.style={circle,thick,draw},square/.style={regular polygon,regular polygon sides=4}]
\node[fill=green!20] (1) at (3.88975,4.46742) {\large \,};
\node[fill=green!20] (2) at (6.19526,5.81822) {\Large $x_1$};
\node[fill=green!20] (3) at (1.56877,5.78742) {\large \,};
\node[fill=green!20] (4) at (3.90551,1.79731) {\large \,};
\node[draw=gray,fill=yellow!50] (5) at (7.74693,6.72753) {\Large $z_1^1$};
\node[draw=gray,fill=yellow!50] (6) at (5.3,7.4) {\Large $z_1^2$};
\node[draw=gray,fill=yellow!50] (7) at (6.4,4.1) {\Large $z_1^3$};
\node[draw=gray,fill=yellow!50] (8) at (1.4,4.2) {\large \,};
\node[draw=gray,fill=yellow!50] (9) at (2.5,7.0) {\large \,};
\node[draw=gray,fill=yellow!50] (10) at (0.0055904,6.67466) {\large \,};
\node[draw=gray,fill=yellow!50] (11) at (3.91459,0) {\large \,};
\node[draw=gray,fill=yellow!50] (12) at (5.40364,2.3) {\large \,};
\node[draw=gray,fill=yellow!50] (13) at (2.42016,2.3) {\large \,};
\node[draw=gray,fill=yellow!50] (14) at (3.88975,6.2) {\large \,};
\node[draw=gray,fill=yellow!50] (15) at (5.40364,3.9) {\large \,};
\node[draw=gray,fill=yellow!50] (16) at (2.42016,3.9) {\large \,};
\end{scope}
\begin{scope}[>={Stealth[black]},
              every edge/.style={draw=black, thick}]
\path [->] (1) edge[bend right=15] node {} (2);
\path [->] (2) edge[bend right=15] node {} (1);
\path [->] (1) edge[bend right=15] node {} (3);
\path [->] (3) edge[bend right=15] node {} (1);
\path [->] (1) edge[bend right=15] node {} (4);
\path [->] (4) edge[bend right=15] node {} (1);
\path [->] (2) edge[>={Stealth[gray!70]},gray!70,bend right=15] node {} (5);
\path [->] (5) edge[>={Stealth[gray!70]},gray!70,bend right=15] node {} (2);
\path [->] (2) edge[>={Stealth[gray!70]},gray!70,bend right=15] node {} (6);
\path [->] (6) edge[>={Stealth[gray!70]},gray!70,bend right=15] node {} (2);
\path [->] (2) edge[>={Stealth[gray!70]},gray!70,bend right=15] node {} (7);
\path [->] (7) edge[>={Stealth[gray!70]},gray!70,bend right=15] node {} (2);
\path [->] (3) edge[>={Stealth[gray!70]},gray!70,bend right=15] node {} (8);
\path [->] (8) edge[>={Stealth[gray!70]},gray!70,bend right=15] node {} (3);
\path [->] (3) edge[>={Stealth[gray!70]},gray!70,bend right=15] node {} (9);
\path [->] (9) edge[>={Stealth[gray!70]},gray!70,bend right=15] node {} (3);
\path [->] (3) edge[>={Stealth[gray!70]},gray!70,bend right=15] node {} (10);
\path [->] (10) edge[>={Stealth[gray!70]},gray!70,bend right=15] node {} (3);
\path [->] (4) edge[>={Stealth[gray!70]},gray!70,bend right=15] node {} (11);
\path [->] (11) edge[>={Stealth[gray!70]},gray!70,bend right=15] node {} (4);
\path [->] (4) edge[>={Stealth[gray!70]},gray!70,bend right=15] node {} (12);
\path [->] (12) edge[>={Stealth[gray!70]},gray!70,bend right=15] node {} (4);
\path [->] (4) edge[>={Stealth[gray!70]},gray!70,bend right=15] node {} (13);
\path [->] (13) edge[>={Stealth[gray!70]},gray!70,bend right=15] node {} (4);
\path [->] (1) edge[>={Stealth[gray!70]},gray!70,bend right=15] node {} (14);
\path [->] (14) edge[>={Stealth[gray!70]},gray!70,bend right=15] node {} (1);
\path [->] (1) edge[>={Stealth[gray!70]},gray!70,bend right=15] node {} (15);
\path [->] (15) edge[>={Stealth[gray!70]},gray!70,bend right=15] node {} (1);
\path [->] (1) edge[>={Stealth[gray!70]},gray!70,bend right=15] node {} (16);
\path [->] (16) edge[>={Stealth[gray!70]},gray!70,bend right=15] node {} (1);
\end{scope}
\end{tikzpicture}

}
\caption{Motivation idea of how to achieve privacy by creating a local network for each agent and distribute the original initial state among the local nodes (augmented states).}
\label{fig:motiv_exp}
\end{figure}
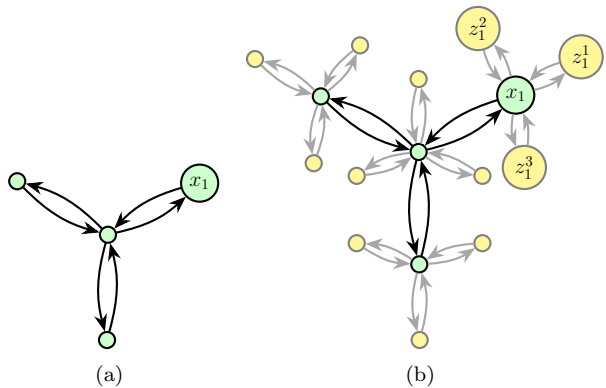

\noindent\textbf{Motivation example.} 
Consider the star network of agents depicted in Figure~\ref{fig:motiv_exp}~(a). 
For almost all parameters, it is possible to invoke the PBH test, to readily notice that the agent in the center of the star network can recover the initial state of every other agent. 
Therefore, we wonder if we can scramble the initial state of each agent in a network local to that agent such that the initial state is no longer recoverable by other agents, see Figure~\ref{fig:motiv_exp}~(b), where the local network of each agent is depicted by the yellow nodes and gray edges.  
In this example, the state of agent 1, $x_1[k]$ becomes a new vector state $\tilde x_1[k]=\left[x_1[k]\,\,\,z_1^1[k]\,\,\,z_1^2[k]\,\,\,z_1^3[k]\right]$. 
Observe that this may imply redesigning the original dynamics matrix not only to accommodate the local networks of each agent but also to re-scale the original weights, while keeping the original connections between agents. 

In fact, the previous motivation example (Figure~\ref{fig:motiv_exp}) is sufficient to satisfy Proposition~\ref{th:PBH} such that  Definition~\ref{def:priv} holds. 
However, it is not sufficient to ensure that the network reaches consensus, as we will, later on, further notice.

In what follows, we do the following assumptions. 

\vspace{1mm} 
{\centering
\allowdisplaybreaks
\noindent \colorbox{carmine!10}{\parbox{0.98\textwidth}{
\noindent$\mathbf{A_{1}}$ The structure of the network is available and cannot be changed (i.e., we cannot change which pairs of agents communicate), but we can choose its parameters.\hfill$\diamond$
}
}
}
\vspace{1mm} 

\vspace{1mm} 
{\centering
\allowdisplaybreaks
\noindent \colorbox{carmine!10}{\parbox{0.98\textwidth}{
\noindent$\mathbf{A_2}$ The network of agents $\mathcal G$ is strongly connected, and it has at least three agents.\hfill$\diamond$ % with at least two of them not working in coalition to recover the initial agents states. 
}
}
}
\vspace{1mm}

Observe that the first assumption ($\mathbf{A_{1}}$) portrays diverse practical applications where it is impossible to change the interconnections naturally due to short-range communications between agents or inherent infrastructure that is too costly to change. 
The second assumption ($\mathbf{A_2}$) reflects that fact that with only two agents, after reaching average consensus, an agent can trivially recover the initial state of the other agent. 
Intuitively, the same reasoning applies if the network has $N+1$ agents and $N$ of them work in coalition, being also able to recover the other agent initial state.

% assume that the structure of the network is available and cannot be changed (i.e., which pairs of agents communicate), and we can choose its parameters. 

Under these assumptions, we can formulate our formal problem statement as follows.

%%%%%%%%%%%%%%%%%%%%%%%%%%%%%%%%%%%%%%%%%%%%%%%%%%%%%%%%%%%%%%%%
%%%%%%%%%%%%%%%%%%%%%%%%%%%%%%%%%%%%%%%%%%%%%%%%%%%%%%%%%%%%%%%%
%%%%%%%%%%%%%%%%%%%%%%%%%%%%%%%%%%%%%%%%%%%%%%%%%%%%%%%%%%%%%%%%
%%%%%%%%%%%%%%%%%%%%%%%%%%%%%%%%%%%%%%%%%%%%%%%%%%%%%%%%%%%%%%%%
%%%%%%%%%%%%%%%%%%%%%%%%%%%%%%%%%%%%%%%%%%%%%%%%%%%%%%%%%%%%%%%%
%%%%%%%%%%%%%%%%%%%%%%%%%%%%%%%%%%%%%%%%%%%%%%%%%%%%%%%%%%%%%%%%

% Under the previous assumptions, we can formulate our formal problem statement as follows.

\vspace{1mm} 
{\centering
\allowdisplaybreaks
\noindent \colorbox{blue!10}{\parbox{0.98\textwidth}{
\noindent$\mathbf{P_{1}}$ Given $N$ agents with a communication digraph $\mathcal G=(\mathcal X,\mathcal E_{\mathcal X,\mathcal X})$, %and $A(\mathcal G)\in\mathbb R^{N\times N}$, 
for each $i=1,\ldots,N$
  \textit{is there} 
  (\emph{i}) an augmented state with \textit{dimension} $d_i\in\mathbb N$ denoted by $z_k^i\in\mathbb R^{d_i}$ which nodal dynamics' state is  $\tilde{x}_i[k]=[x_i[k]\,\,z_i[k]]$,
  (\emph{ii}) initial condition $\tilde{x}_i[0]\in\mathbb R^{d_i+1}$,
  and (\emph{iii}) \textit{augmented dynamics}
  $\tilde W_i\in\mathbb R^{(d_i+1)\times(d_i+1)}$ satisfying  
%$A=[a_{ij}]_{\substack{i=1\\j=1}}^N$ such that $\bar A=\bar A(\mathcal G)$ (i.e., the zero entries of $A(\mathcal G)$ cannot be changed), a set of dimensions $\{d_i\}_{i=1}^N$, and a set of functions $\{f_i,g_i\}_{i=1}^N$ such that
%\begin{equation}

\begin{center}
\textbf{-- Augmented  Dynamics --}
\end{center}
\begin{fleqn}[\parindent]
\begin{subequations}
\begin{align}
\label{eq:P1a}
\bullet\;\; \begin{array}{ll}\tilde x_i[k+1] &=\tilde W_i \tilde x_i[k] +\tilde u^i[k], \\[0.1cm]
\tilde x_i[k]&=\left[x_i[k]\;\,z_i[k]\right],\text{ with } z_i[k]\in\mathbb R^{d_i}\text{ and}\\[0.1cm]
\tilde u_i[k]&=\left[u_i[k]\;\,\mathbf{0}_{1\times d_i}\right],\,\,
% \end{array}
% \end{align}
% % \end{subequations}
% % \end{fleqn}
% % \begin{fleqn}[\parindent]
% % \begin{subequations}
% \begin{align}
% \label{eq:P1b}
u_i[k] =\displaystyle\sum_{j\in\mathcal N_i\setminus\{i\}} w_jx_j[k] 
%\\[0.15cm]
%& \text{(is the same input feedback as in~\eqref{eq:output})}
\end{array}
\end{align}
% \begin{align}
% \label{eq:P1c}
% \bullet\;\;\tilde x_{k+1} &= \tilde A(G) \tilde x_k+u_k
% \end{align}
% \end{subequations}
% \end{fleqn}
% \begin{fleqn}[\parindent]
% %x_0^i=f_i(x_0^i,\theta_0^i)-x_0^i\\[0.2cm]
% \begin{subequations}
such that the following holds  
\begin{center}
\textbf{-- Specifications --}
\end{center}
\textbf{Consensus}
\begin{align}
\label{eq:P1d}
\bullet\;\;\displaystyle\mathop{\lim}_{k\to\infty}\tilde{x}_i[k] &= \mathbf{1}_{d_i+1}^\intercal\left(\displaystyle\frac{1}{N}\sum_{j=1}^N x_j[0]\right)
\end{align}

\textbf{Privacy}
\begin{align}
\label{eq:P1c}
 \bullet\;\;\begin{array}{l}\text{agent }j \text{ cannot recover the initial state of any agent},\\ i\neq j,\text{ according to
%\mathbf{1}_{d_i+1}^\intercal \tilde{x}_i[0]
(Definition~\ref{def:priv})}
%\displaystyle\sum_{l=1}^{d_i}[z_0^i]^l,\, i\in\{1,\ldots,N\}\setminus\{j\}
\end{array}
\end{align}
% \end{subequations}
% \end{fleqn}
% \begin{fleqn}[\parindent]
% \begin{subequations}
\end{subequations}
\end{fleqn}
%\end{equation}
}
}
}
\vspace{1mm}

Notice that when $d_i = 0$, the agents do not have augmented states and, therefore, the third condition, \eqref{eq:P1c}, boils down to assess if the agent's state belongs to the unobservable subspace of the system. 
Furthermore, notice that we do not allow to change (add/remove) the connections between agents given in $\mathcal G$, but we allow each node to have an internal subnetwork (which may be seen as having memory). 
Also, each agent does not require extra communication capabilities, and it continues to send a scalar value in each iteration. 

That said, and briefly speaking, the proposed solution will be able to ensure privacy by crafting the augmented network dynamics and splitting the scalar initial state of an agent among its augmented states. 

Additionally, we would like, if possible, to find a solution for $\mathbf{P_1}$ that can be achieved distributively, which leads to the following. 

\vspace{1mm} 
{\centering
\allowdisplaybreaks
\noindent \colorbox{blue!10}{\parbox{0.98\textwidth}{
\noindent$\mathbf{P_{1}^D}$ consist of $\mathbf{P_1}$ enforcing the solution to be distributed.  
}
}
}
\vspace{1mm}

\subsection{Roadmap of the main results}

%{\color{red}\Large MUDAR! e colocar no Roadmap}
% Dizer quais os teoremas para correccao e complexidade e algoritmos

The rest of the paper focused on addressing the above mentioned problems. We achieve these solutions in a series of steps, focusing on various intermediate objectives and relaxed problem variants. In what follows we summarize the flow of the remainder of the paper by briefly reviewing these intermediate objectives and outlining their solution strategy by pointing to the corresponding key technical results that are developed in the subsequent sections.

\noindent\textbf{Objective 1:} Obtaining a solution to a relaxed version of $\mathbf{P_1}$, say $\mathbf{P_1^R}$, by eliminating the consensus requirement $(3b)$.

\begin{figure}[H]
\centering
\subfigure[]{
\begin{tikzpicture}[scale=.54, transform shape,node distance=1.5cm]
\begin{scope}[every node/.style={circle,thick,draw},square/.style={regular polygon,regular polygon sides=4}]
%\node[fill=green!20] (1) at (3.88975,4.46742) {\large \,};
\node[fill=green!20] (2) at (0,0) {\Large $x_1$};
% \node[fill=green!20] (3) at (1.56877,5.78742) {\large \,};
% \node[fill=green!20] (4) at (3.90551,1.79731) {\large \,};
\node[draw=gray,fill=yellow!50] (6) at (-2,0.7) {\Large $z_1^1$};
\node[draw=gray,fill=yellow!50] (7) at (2,0.7) {\Large $z_1^2$};
% \node[draw=gray,fill=yellow!50] (8) at (.4,4.8) {\large \,};
% \node[draw=gray,fill=yellow!50] (9) at (1.5,7.4) {\large \,};
% \node[draw=gray,fill=yellow!50] (12) at (5.40364,1.3) {\large \,};
% \node[draw=gray,fill=yellow!50] (13) at (2.42016,1.3) {\large \,};
% \node[draw=gray,fill=yellow!50] (15) at (5.40364,3.9) {\large \,};
% \node[draw=gray,fill=yellow!50] (16) at (2.42016,3.9) {\large \,};
\end{scope}
\begin{scope}[>={Stealth[black]},
              every edge/.style={draw=black, thick}]
% \path [->] (1) edge[bend right=15] node {} (2);
% \path [->] (2) edge[bend right=15] node {} (1);
% \path [->] (1) edge[bend right=15] node {} (3);
% \path [->] (3) edge[bend right=15] node {} (1);
% \path [->] (1) edge[bend right=15] node {} (4);
% \path [->] (4) edge[bend right=15] node {} (1);
\path [->] (2) edge[>={Stealth[gray!70]},gray!70,bend right=15] node {} (6);
\path [->] (6) edge[>={Stealth[gray!70]},gray!70,bend right=15] node {} (2);
\path [->] (2) edge[>={Stealth[gray!70]},gray!70,bend right=15] node {} (7);
\path [->] (7) edge[>={Stealth[gray!70]},gray!70,bend right=15] node {} (2);
% \path [->] (3) edge[>={Stealth[gray!70]},gray!70,bend right=15] node {} (8);
% \path [->] (8) edge[>={Stealth[gray!70]},gray!70,bend right=15] node {} (3);
% \path [->] (3) edge[>={Stealth[gray!70]},gray!70,bend right=15] node {} (9);
% \path [->] (9) edge[>={Stealth[gray!70]},gray!70,bend right=15] node {} (3);
% \path [->] (4) edge[>={Stealth[gray!70]},gray!70,bend right=15] node {} (12);
% \path [->] (12) edge[>={Stealth[gray!70]},gray!70,bend right=15] node {} (4);
% \path [->] (4) edge[>={Stealth[gray!70]},gray!70,bend right=15] node {} (13);
% \path [->] (13) edge[>={Stealth[gray!70]},gray!70,bend right=15] node {} (4);
% \path [->] (1) edge[>={Stealth[gray!70]},gray!70,bend right=15] node {} (15);
% \path [->] (15) edge[>={Stealth[gray!70]},gray!70,bend right=15] node {} (1);
% \path [->] (1) edge[>={Stealth[gray!70]},gray!70,bend right=15] node {} (16);
% \path [->] (16) edge[>={Stealth[gray!70]},gray!70,bend right=15] node {} (1);
\end{scope}
\end{tikzpicture}
}
\subfigure[]{
     $
     \begin{bmatrix}
        0 & a_{1,2} & a_{1,3}\\
        a_{2,1} & 0 & 0 \\
        a_{3,1} & 0 & 0
     \end{bmatrix}
     $
}
\caption{Solution to Challenge 1.}
\label{fig:road_2cp}
\end{figure}
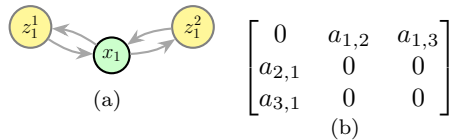

\noindent\textbf{Solution Outline (Theorem~\ref{th:main}):} Consider two additional memory states with the connections depicted in Figure~\ref{fig:road_2cp}~(a), with dynamics matrix in Figure~\ref{fig:road_2cp}~(b), following the same color scheme as in the motivation example (for each node of the original network a similar design is conducted).  
The solution yields for almost all possible numerical realization of the dynamics matrix parameters. 
Additionally, the initial state is distributed among the different states such that non-zero initial condition is given to any of the two augmented states

\noindent\textbf{Challenge 2:} Finding a solution to Problem $\mathbf{P_1}$.

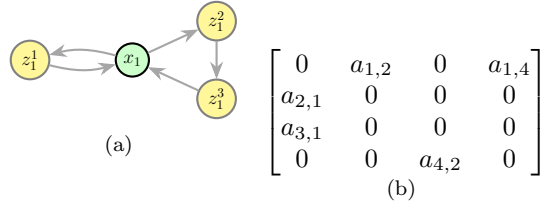
\begin{figure}[H]
\centering
\subfigure[]{
\begin{tikzpicture}[scale=.5, transform shape,node distance=1.5cm]
\begin{scope}[every node/.style={circle,thick,draw},square/.style={regular polygon,regular polygon sides=4}]
\node[fill=green!20] (1) at (2.72064,1.03433) {\Large $x_1$};
\node[draw=gray,fill=yellow!50] (2) at (0.,1.03407) {\Large $z_1^1$};
\node[draw=gray,fill=yellow!50] (3) at (4.96021,2.06811) {\Large $z_1^2$};
\node[draw=gray,fill=yellow!50] (4) at (4.96004,0.) {\Large $z_1^3$};
\end{scope}
\begin{scope}[>={Stealth[black]},
              every edge/.style={draw=black, thick}]
\path [->] (1) edge[>={Stealth[gray!70]},gray!70,bend right=15] node {} (2);
\path [->] (1) edge[>={Stealth[gray!70]},gray!70] node {} (3);
\path [->] (2) edge[>={Stealth[gray!70]},gray!70,bend right=15] node {} (1);
\path [->] (3) edge[>={Stealth[gray!70]},gray!70] node {} (4);
\path [->] (4) edge[>={Stealth[gray!70]},gray!70] node {} (1);
\end{scope}
\end{tikzpicture}
}
\subfigure[]{
     $
     \begin{bmatrix}
        0 & a_{1,2} & 0 & a_{1,4}\\
        a_{2,1} & 0 & 0 & 0 \\
        a_{3,1} & 0 & 0 & 0\\
        0 & 0 & a_{4,2} & 0
     \end{bmatrix}
     $
} 
\caption{Solution to Challenge 2.}
\label{fig:road_3cp}
\end{figure}

\noindent\textbf{Solution (Theorem~\ref{cor:3cpy}):} Consider three augmented states as depicted in Figure~\ref{fig:road_3cp}~(a) (with all the possibilities one can explore later on), where the dynamics is as in Figure~\ref{fig:road_3cp}~(b). 
The parameter selection is $a_{2,1}=2/Z_2$, $a_{1,2}=1/Z_1$, $a_{1,4}=1/Z_1$, $a_{3,1}=1/Z_3$, and $a_{4,2}=1/Z_4$, where $Z_i$ is the sum of the row~$i$ entries. 
In addition, because in this scenario the consensus value is $(v^0)^\intercal 1$, where $v^0$ is the left-eigenvector associated with the extended dynamics matrix, we need to distribute the initial value of the agent's state among the extended states, and scale it to account for the total number of nodes in the extended network and also to mitigate the  weight imposed by $v^0$. 
First, each agent distributes the initial condition across its augmented states (setting as $0$ for the original state), the distribution can be tailored by each agent, and scaling it to make the new average of initial states equal to the original one, i.e., 
\[\frac{4 x_i[0]}{z_i^1[0] + z_i^2[0] + z_i^3[0]} \begin{bmatrix}0 & z_i^1[0] & z_i^2[0] & z_i^3[0]\end{bmatrix}.
\] 
Second, each agent and augmented states re-scale the new initial condition, $x_j[0]$, to be 
\[
\frac{x_j[0]}{4N v_j^0}\sum_{i=1}^{4N}v_i^0.
\]

%  4 # x[[i]]/Total[#] &@RandomReal[{0, 1}, 3]
%  ((Total[#]/(Length[x0] #)) &@left) x0

Because the computation of $v^0$ is centralized, and it may lead to numerical challenges, we also propose the following.

\noindent\textbf{Challenge 3:} Finding a solution to $\mathbf{P_1^D}$.  

\begin{figure}[H]
\centering
\subfigure[]{
\begin{tikzpicture}[scale=.5, transform shape,node distance=1.5cm]
\begin{scope}[every node/.style={circle,thick,draw},square/.style={regular polygon,regular polygon sides=4}]
\node[fill=green!20] (1) at (2.15841,0.) {\Large $x_1$};
\node[draw=gray,fill=yellow!50] (2) at (2.156,3.01787) {\Large $z_1^1$};
\node[draw=gray,fill=yellow!50] (3) at (2.8,1.50681) {\Large $z_1^2$};
\node[draw=gray,fill=yellow!50] (4) at (4.67655,1.51032) {\Large $z_1^3$};
\node[draw=gray,fill=yellow!50] (5) at (0.,1.50851) {\Large $z_1^4$};
\end{scope}
\begin{scope}[>={Stealth[black]},
              every edge/.style={draw=gray, thick}]
\path [-] (1) edge node {} (2);
\path [-] (1) edge node {} (3);
\path [-] (1) edge node {} (4);
\path [-] (1) edge node {} (5);
\path [-] (2) edge node {} (3);
\path [-] (2) edge node {} (4);
\path [-] (2) edge node {} (5);
\path [-] (3) edge node {} (4);
\path [-] (3) edge node {} (5);
\end{scope}
\end{tikzpicture}
}
\subfigure[]{
     $
     \begin{bmatrix}
        0 & a_{1,2} & 0 & a_{1,4}\\
        a_{2,1} & 0 & 0 & 0 \\
        a_{3,1} & 0 & 0 & 0\\
        0 & 0 & a_{4,2} & 0
     \end{bmatrix}
     $
}
\caption{Solution to Challenge 3.}
\label{fig:road_4cp}
\end{figure}
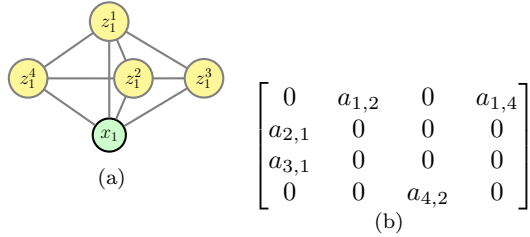

\noindent\textbf{Solution (Theorem~\ref{th:main_P2}):} 
Under the assumption that the graph is bidirectional, we consider the augmented state as in Figure~\ref{fig:road_4cp}~(a) with dynamics as in Figure~\ref{fig:road_4cp}~(b). %, where the weights are as follows:...
From this, we can compute the left-eigenvector using only local information, see details in Algorithm~\ref{alg:buildAP_rev}.  %$[v_0]_i=....$ thus enabling a distributed design. 	

\section{Deterministic privacy preserving average consensus for multi-agent systems}\label{sec:main_res}

We first propose to address a restricted version $\mathbf{P^R_1}$, where we do waive the consensus condition (Section~\ref{sub:main_res_lti}). 
In particular, to build an intuition we show that simple solutions can be proposed to address $\mathbf{P_1^R}$, but these are insufficient to address $\mathbf{P_1}$. 
Subsequently, by identifying the conditions that limit the achievement of consensus, we identify a class of augmented network structures and particular realizations of the dynamics that solve $\mathbf{P_1}$ (Section~\ref{sub:priv1}). 
Lastly, we show that under some additional assumptions, it is possible to deploy a solution to $\mathbf{P_1^D}$, i.e., a solution to $\mathbf{P_1}$ that can be achieved in a distributed fashion (Section~\ref{sub:rev_mat}). 

\subsection{Achieving privacy in discrete-time LTI multi-agent systems}\label{sub:main_res_lti}

The first question that emerges is if it is possible to extend the network such that neighbors' states could be perceived as neighbors of neighbors, to guarantee that some of the new states are private -- see Remark~\ref{rmk2}. 

\begin{remark}\label{rmk2}
The neighbors can always have access to the initial state of the neighbors, and therefore privacy in the strict sense above cannot be achieved.
\end{remark}

\begin{figure}[ht!]
\begin{center}
\subfigure[]{
$
\begin{array}{c}
A=\begin{bmatrix}
 0 & \frac{1}{2} & \frac{1}{2}\\[1mm]
 \frac{1}{2} & 0 & \frac{1}{2}\\[1mm]
 \frac{1}{2} & \frac{1}{2} & 0
\end{bmatrix}
 \\
 \,
\end{array}
$
}
\subfigure[]{
    \begin{tikzpicture}[scale=.6, transform shape,node distance=1.5cm]
\begin{scope}[every node/.style={circle,thick,draw},square/.style={regular polygon,regular polygon sides=4}]
\node[fill=green!20] (1) at (4.95793,5.94893) {\large $1$};
\node[fill=green!20] (2) at (6.70249,2.40956) {\large $2$};
\node[fill=green!20] (3) at (2.76639,2.67303) {\large $3$};
\end{scope}
\begin{scope}[>={Stealth[black]},
              every edge/.style={draw=black, thick}]
\path [->] (1) edge[bend right=15] node {} (2);
\path [->] (1) edge[bend right=15] node {} (3);
\path [->] (2) edge[bend right=15] node {} (1);
\path [->] (2) edge[bend right=15] node {} (3);
\path [->] (3) edge[bend right=15] node {} (1);
\path [->] (3) edge[bend right=15] node {} (2);
\end{scope}
\end{tikzpicture}
}
\caption{Dynamics matrix A in~(a) and its digraph representation (of a bidirected cycle network with 3 agents)~(b).}
\label{fig:cycle3_trivial}
\end{center}
\end{figure}
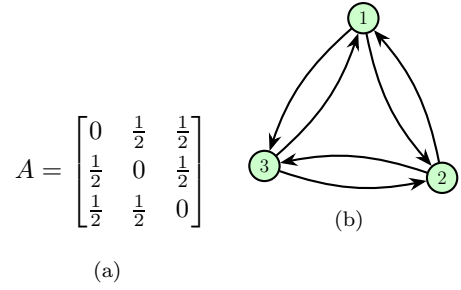

% The problems that we address in this paper are the following.
%First, because for an arbitrary dynamics it may happen that the states are observable from every agent. 
For instance, consider A as in Figure~\ref{fig:cycle3_trivial}~(a), with digraph representation depicted in Figure~\ref{fig:cycle3_trivial}~(b). 
In this network, it is trivial to see that any agent can recover the initial state of all the agents in the network. 

%.. and maybe introduce a simple example where this occurs, maybe a path with three agents like the simulations we did... then this motivates the need to look into P2

%that said, theorem should come earlier, i.e., before P2 and likely before the motivating example

Using the Popov criteria, Theorem~\ref{th:PBH}, we have that, for any agent,
 \[
     P_{O(A,C)}^\lambda =
     \begin{bmatrix}
        C\\
        \lambda I_N - A
     \end{bmatrix},
 \]
 where $C=I_{3}$. And we get 
  \[
     P_{O(A,C)}^\lambda =
     \begin{bmatrix}
        1 & 0 & 0\\
        0 & 1 & 0\\
        0 & 0 & 1\\
        \lambda & -\frac{1}{2} & -\frac{1}{2}\\
        -\frac{1}{2} & \lambda & -\frac{1}{2}\\
        -\frac{1}{2} & -\frac{1}{2} & \lambda     
    \end{bmatrix},
 \]
 and $rank\left(P_{O(A,C)}^\lambda\right)$ is full for any $\lambda\in\mathbb C$. Hence, the system is observable. 

% maybe we should use the Remark framework to emphasize the limitation. Then we may want to do the following:(i) give an example where the state of an agent is not retrievable by other agent, say a star network where one of the corners cannot retrieve the other corners. Even show the popov criteria on this to make it crystal clear. 

Intuitively, due to the structure, we may have agents in one place that guarantee that the relative position of the other agent is such that it cannot retrieve other agent's state. 
This property raises the question: 
Is it possible to extend the network by equipping the agents with memory such that such relative positions hold for the newly added states (at least two)? 
This way, we could distribute the initial state of the agent among the augmented states (memory) such that the sum is not retrieved. This intuition leads us to the problem below, a restricted version of problem $\mathbf{P_1}$, where we remove condition~\eqref{eq:P1c}.

\vspace{1mm} 
{\centering
\allowdisplaybreaks
\noindent \colorbox{blue!10}{\parbox{0.98\textwidth}{
\noindent$\mathbf{P_1^R}$ consists of $\mathbf{P_1}$ only with the privacy specification, i.e., satisfying \eqref{eq:P1a} and \eqref{eq:P1c}.

}
}
}
\vspace{1mm}

Next, we propose a design following the procedure described in Algorithm~\ref{alg:buildAP} that attains the feasibility question posed in $\mathbf{P_1^R}$, i.e., that constructs a matrix $A^\text{\textsf{P}}$ such that 
$$\tilde x^\text{\textsf{P}}[k+1]=A^\text{\textsf{P}}\tilde x^\text{\textsf{P}}[k]$$ 
and the conditions of $\mathbf{P_1^R}$ hold. 
$\tilde x^\text{\textsf{P}}$ is the $\tilde x$ of Problem~$\mathbf{P_1}$ with agents in a permuted order, and the entries of $A^\text{\textsf{P}}\tilde x^\text{\textsf{P}}$ are the entries of the $\tilde{W}_i$ following the same permuted order of agents. 
In the problem statement, to better understand the problem, we relabel the agents and augmented states such that they sequentially follow the order of an agent followed by its augmented states, and so on. 
Notwithstanding, for design purposes, we leave the original agents with the original label and label, sequentially, the augmented states of each original agent.

\begin{algorithm}[H]
		{%\small
		\caption{Construction of a privacy dynamics matrix $A^\text{\textsf{P}}$ that solves $\mathbf{P_1^R}$}
		\label{alg:buildAP}
		\begin{algorithmic}[1]
			\STATE{\textbf{input}: network of agents $\mathcal G=\langle \mathcal X,\mathcal E_{\mathcal X,\mathcal X}\rangle$ and initial state $x[0]\in\mathbb R^{N}$}
			\STATE{\textbf{output}: dynamics matrix of the augmented network $A^\text{\textsf{P}}\in\mathbb R^{3N\times 3N}$ and augmented initial state $\tilde x^{\text{\textsf{P}}}[0]\in\mathbb R^{3N}$}
		\STATE{\textbf{set} $A\in\mathbb R^{N\times N}$ as the adjacency matrix of $\mathcal G$}
		\STATE{\textbf{fill} the entries of $A^{\text{\textsf{P}}}$ with zeros}
		\STATE{\textbf{set}
		{\color{gray}\hfill $\rhd$ 
		
		The digraph representation is augmented placing two additional nodes per each original node that only connect bidirectly to the original node, %whcih already 
		setting some of the dynamics matrix entries to zero, and then the non-zero entries will be as follows. 
		\underline{Copy the matrix $A$ to the first $n$ rows and $n$ columns}}
		$$A^{\text{\textsf{P}}}_{ij}=A_{ij}, \forall i,j\in\{1,\ldots,N\}\text{ and }i\neq j$$}
		\STATE{\textbf{for} $i=1,\ldots,N$\hfill{\color{gray}$\rhd$ \underline{Set additional entries values}}
		$$
		\begin{array}{ll}
		A^{\text{\textsf{P}}}_{i,N+2i-1} = a_{i,N+2i-1}, & A^{\text{\textsf{P}}}_{i,N+2i}=a_{i,N+2i},  \\
		A^{\text{\textsf{P}}}_{N+2i-1,1} = a_{N+2i-1,1},  & 
		A^{\text{\textsf{P}}}_{N+2i,1} = a_{N+2i,1}
		\end{array}
		$$
		}
		\STATE{\textbf{distribute} the initial value of each original agent across its augmented states to build $\tilde x^{\text{\textsf{P}}}[0]$}
		\end{algorithmic} 
		}
\end{algorithm} 

% The role of assigning different values in step~5 and step~6 of Algorithm~\ref{alg:buildAP} will be crucial to ensure privacy, as we subsequently show. 

Now, we show that given a matrix $A\in\mathbb R^{N\times N}$ and initial state $x[0]\in\mathbb R^{N}$ the new parametric matrix $A^\text{\textsf{P}}\in\mathbb R^{3N\times 3N}$ and augmented initial state $\tilde x^{\text{\textsf{P}}}[0]\in\mathbb R^{3N}$, constructed with Algorithm~\ref{alg:buildAP} from $A$, ensures that the initial states of the original network are kept private.

\vspace{1mm} 
{\centering
\allowdisplaybreaks
\noindent \colorbox{c3!15}{\parbox{0.98\textwidth}{
\begin{theorem}\label{th:main}
    The parametric matrix $A^\text{\textsf{P}}\in\mathbb R^{3N\times 3N}$ resulting from Algorithm~\ref{alg:buildAP} with input matrix $A$, ensures that there is not a node $i\in\{1,\ldots,N\}$ capable of recovering every initial state in $\{x_j[0],x_{N+2j-1}[0],x_{N+2j}[0]\}$ by observing every direct neighbor node $j\in\mathcal N_i$. 
 Formally, the canonical vectors $$e_{3N}^{N+2j-1},e_{3N}^{N+2j}\notin span\left( P_{O(A^\text{\textsf{P}},C)}^\lambda\right)\; \forall \lambda\in\mathbb R\,\forall j\in\mathcal N_i,$$ 
 with  $C=\left[\begin{smallmatrix}
 e^i_N & e^{n_1}_N & \cdots e^{n_k}_N& 
\end{smallmatrix}\right]^\intercal$, where $\mathcal N_i=\{n_1,\cdots,n_k\}$. 
Further, if we have the additional goal that the original initial state (distributed across the augmented states) cannot be recovered, i.e., if
$$e_{3N}^{N+2j-1}+e_{3N}^{N+2j}\notin span\left( P_{O(A^\text{\textsf{P}},C)}^\lambda\right)\; \forall \lambda\in\mathbb R\,\forall j\in\mathcal N_i,$$
then $0\neq a_{j,N+2j-1}\neq a_{j,N+2j-1}\neq 0$.
\hfill$\circ$   
\end{theorem}
}
}
}
\vspace{1mm}

\begin{proof}
    Let $A^\text{\textsf{P}}\in\mathbb R^{3N\times 3N}$ be the parametric matrix resulting from Algorithm~\ref{alg:buildAP} with input matrix $A$. 
    We have that 
    $$
\bigA^\text{\textsf{P}}=
\begin{bmatrix}
\bigA & \rvline  & \bigR
%   \begin{matrix}
%   a_{1,N+1} & a_{1,N+2} & 0 & \cdots & 0 \\
%     & \ddots & & & \vdots \\
%   0 & \cdots & 0 & a_{N,3N-1} & a_{N,3N} \\
%   \end{matrix}
  \\
\hline
\bigL
%   \begin{matrix}
%   a_{N+1,1} & 0 & \cdots & 0 \\
%   a_{N+2,1} & 0 & \cdots & 0 \\
%   & \ddots & & \vdots \\
%   0 & \cdots & 0 & a_{3N-1,N} \\
%   0 & \cdots & 0 & a_{3N,N} \\
%   \end{matrix}
  & \rvline & \bigzero
\end{bmatrix},
    $$
where $\mathbf{0}$ is the $2N\times 2N$ matrix of zeros, the matrix 
\begin{equation}\label{eq:matR}
\bigR_{N\times 2N} = 
  \begin{bmatrix}
  a_{1,N+1} & a_{1,N+2} & 0 & \cdots & 0 \\
    & \ddots & & & \vdots \\
  0 & \cdots & 0 & a_{N,3N-1} & a_{N,3N} \\
  \end{bmatrix},
\end{equation}
and the matrix 
$$
\bigL_{2N\times N} =
  \begin{bmatrix}
  a_{N+1,1} & 0 & \cdots & 0 \\
  a_{N+2,1} & 0 & \cdots & 0 \\
  & \ddots & & \vdots \\
  0 & \cdots & 0 & a_{3N-1,N} \\
  0 & \cdots & 0 & a_{3N,N} \\
  \end{bmatrix}.
$$
Next, we compute for all $\lambda\in\mathbb R$ the span of the $P_{O(A^\text{\textsf{P}},C)}^\lambda$, with the output matrix $C\in\mathbb R^{N\times 3N}$ that consists in observing each state variable that is a neighbor node of $i$ (including node $i$), 
$$
\bigC=
\begin{bmatrix}
\bigI_N^{\mathcal N_i} & \rvline  & \bigzero_{N\times 2N} \\
\bigzero_{2\times N} & \rvline & \bigV_{c_i^1,c_i^2} \\
\end{bmatrix}.$$
Notwithstanding, we will do a proof that shows a stronger property, by placing an output in each state variable of the original $A$ plus the two augmented states of node $i$. 
In other words, we do the proof using the an output matrix with more output capabilities than the required ones. 
Without loss of generality, suppose that $i=N$. If this is not the case, we may apply a permutation in the node labels and obtain this case. 
The output matrix is
$$
\bigC=
\begin{bmatrix}
\bigI_N & \rvline  & \bigzero_{N\times 2N} \\
\bigzero_{2\times N} & \rvline & \bigV_{c_i^1,c_i^2} \\
\end{bmatrix}, $$
with
$$
\bigV_{c_i^1,c_i^2} = 
\begin{bmatrix}
\bigzero_{2\times 2(N-1)} & \bigI_{2}
\end{bmatrix}.
    $$
% Observe that this is stronger than the stated observability capabilities. We would only require to use the rows of $C$ that correspond to node $i$ and its neighbors. 
Therefore, we have that 

$$
\begin{array}{rl}
span\left(P_{O(A^\text{\textsf{P}},C)}^\lambda\right)   & =  
span\left(
\begin{bmatrix}
\bigI_N & \rvline & \bigzero_{N\times 2N} \\ 
\bigzero_{2\times N} & \rvline & \bigV_{c_i^1,c_i^2} \\
\hline
\biglambdaI_{N}-\bigA & \rvline  & \bigR_{N\times 2N} \\ %\hline
\bigL_{2N\times N}  & \rvline & \biglambdaI_{2N}
\end{bmatrix}
\right)\\
 & = span\left(
\begin{bmatrix}
\bigI_N & \rvline & \bigzero_{N\times 2N} \\ 
\bigzero_{2\times N} & \rvline & \bigV_{c_i^1,c_i^2} \\
\hline
\bigzero_{N\times N} & \rvline  & \bigR_{N\times 2N} \\ %\hline
\bigzero_{2N\times N} & \rvline & \biglambdaI_{2N}
\end{bmatrix}
\right),
\end{array}
$$
for all $\lambda\in\mathbb R$. 
Now, we identify two possible scenarios for $\lambda$, which are $\lambda\neq 0$ or $\lambda = 0$. 
If $\lambda\neq 0$, we get 
$$
span\left(P_{O(A^\text{\textsf{P}},C)}^\lambda\right)   = 
span\left(
\begin{bmatrix}
\bigI_{N}    & \rvline & \bigzero_{N\times 2N} \\ \hline
\bigzero_{N\times N} & \rvline  & \bigzero_{N\times 2N} \\ %\hline
\bigzero_{2N\times N} & \rvline & \bigI_{2N}
\end{bmatrix}
\right)
 = 
\mathbb R^{3N}.
$$
However, when $\lambda=0$, we obtain  
$$
\begin{array}{rl}
span\left(P_{O(A^\text{\textsf{P}},C)}^0\right)   & = \\[0.4cm]
span\left(
\begin{bmatrix}
\bigI_N    & \rvline & \multicolumn{2}{c}{\bigzero_{N\times 2N}}\\
\bigzero_{2\times N} & \rvline & \bigzero_{2\times 2(N-1)} & \bigI_{2}
 \\
 \hline
\bigzero_{N\times N} & \rvline  & \bigR_{N\times 2(N-1)} & \bigzero_{N\times 2} \\ %\hline
\bigzero_{2N\times N} & \rvline &\multicolumn{2}{c}{\bigzero_{N\times 2N}}
\end{bmatrix}
\right) & =
  \\[0.4cm] 
 
span\left(
 \begin{bmatrix}
\bigI_N    & \rvline & \multicolumn{2}{c}{\bigzero_{N\times 2N}}\\
\bigzero_{2\times N} & \rvline & \bigzero_{2\times 2(N-1)} & \bigI_{2}
 \\ \hline
\bigzero_{N\times N} & \rvline  & \bigR_{N\times 2(N-1)} & \bigzero_{N\times 2} \\ %\hline
%\multicolumn{2}{c}{\bigzero_{N\times 2N}} \\ %\hline
\end{bmatrix}
\right).
\end{array}
$$
Hence, for all $ j\in\{1,\ldots,N\}\setminus\{i\}$, we have that $\alpha e_{3N}^{N+2j-1}+\beta e_{3N}^{N+2j}\in span\left( P_{O(A^\text{\textsf{P}},C)}^\lambda\right)$, and from~\eqref{eq:matR}, with $\alpha=a_{j,N+2j-1}\neq 0$ and $\beta=a_{j,N+2j-1}\neq 0$, it readily follows that $e_{3N}^{N+2j-1},e_{3N}^{N+2j}\notin span\left( P_{O(A^\text{\textsf{P}},C)}^\lambda\right)$. 
If we have the additional goal that $e_{3N}^{N+2j-1}+e_{3N}^{N+2j}\notin span\left( P_{O(A^\text{\textsf{P}},C)}^\lambda\right)$ then we further require that $\alpha\neq\beta$. 
\end{proof}

If two or more agents form a coalition, will they be able to retrieve other agents' states? The answer is no, as noted in the following observation. 

\begin{remark}\label{rmk:1}
From Theorem~\ref{th:main}, we have that coalition between a set of agents does not compromise the privacy of the remaining agents.  
\end{remark}

The next step is to address $\mathbf{P_1}$ by including the average consensus restriction~\eqref{eq:P1d}.

\subsection{Solution to problem $\mathbf{P_1}$}\label{sub:priv1}

Next, our goal is to apply the previous strategy to achieve privacy in average consensus, using an LTI system. 
Given the initial states of $N$ agents, $x_i[0]\in\mathbb R$, for $i=1,\ldots,N$ and a matrix $A$, the state update of the agents at time $k+1$ with $k\in\mathbb N$ is described as 
\begin{equation}\label{eq:LTI_consensus}
    x_i[k+1] = \sum_{j\in\mathcal N\cup\{i\}}x_j[k].
\end{equation}
The following result from~\cite{xiao2004fast} states the required conditions so that Eq.~\eqref{eq:LTI_consensus} reaches asymptotic consensus.

{\centering
\allowdisplaybreaks
 \noindent \colorbox{c7!10}{\parbox{.935\columnwidth}{
\begin{proposition}[\hspace{-0.1mm}\cite{xiao2004fast}]\label{th:consensus_reqs}
    The LTI system described in~\eqref{eq:LTI_consensus} reaches  distributed consensus under the technical requirement that $v_L$ is normalized so that $v_L^\intercal\mathbf{1}=1$ if and only if the matrix $A$ satisfies:
\begin{enumerate}
    \item $A$ has a simple eigenvalue $1$, and all other eigenvalues have magnitude less than 1;
    \item the left and right eigenvectors of $A$ corresponding to the eigenvalue $1$ are $v_L^\intercal$ and $\mathbf{1}$, respectively.\hfill$\circ$ 
\end{enumerate}
\end{proposition}
}
}
}

Hence, if we start with the matrix $A\in\mathbb R^{N\times N}$ and then apply Algorithm~\ref{alg:buildAP} to produce $A^\text{\textsf{P}}\in\mathbb R^{3N\times 3N}$, we need to ensure that the $A^\text{\textsf{P}}$ entries' values ensure:
\begin{enumerate}
    \item[(i)] the conditions of Proposition~\ref{th:consensus_reqs};
    \item[(ii)] the additional condition in Theorem~\ref{th:main} that $\forall j=1,\ldots,N$ 
    $a_{j,N+2j-1}\neq a_{j,N+2j-1}$. 
\end{enumerate}

In other words, the matrix $A^\text{\textsf{P}}$ needs to be row-stochastic to ensure $(1)$, and we can normalize the initial state of each agent $x_i[0]$ to be $\tilde x_i^{\text{\textsf{P}}}[0]= \frac{1}{3N{[v_L]}_i}x_i[0]$ to ensure that Eq.~\eqref{eq:LTI_consensus} converges to $
v_L^\intercal \tilde x^{\text{\textsf{P}}}[0]$, with $\tilde x^{\text{\textsf{P}}}=\begin{bmatrix}
   \tilde x_1^{\text{\textsf{P}}}[0] & \cdots & \tilde x_{3N}^{\text{\textsf{P}}}[0]
\end{bmatrix}^\intercal$. 
Notice that $
v_L^\intercal \tilde x^{\text{\textsf{P}}}[0] = \displaystyle\frac{1}{3N}\sum_{i=1}^{3N} {[v_L]}_i\frac{1}{{[v_L]}_i}\tilde{x}_i^{\text{\textsf{P}}}[0]=\frac{1}{3N}\tilde{x}^{\text{\textsf{P}}}[0]$, as envisioned. 

The first result is a negative one, stating that the previous construction does not fulfill the necessary conditions to reach average consensus stated in Proposition~\ref{th:consensus_reqs}. 
%This result builds up on a result of the Perron–Frobenius theorem, which is the following.

% \vspace{1mm} 
% {\centering
% \allowdisplaybreaks
% \noindent \colorbox{c3!10}{\parbox{0.98\textwidth}{
% \begin{theorem}[Perron–Frobenius]\label{th:period}
%     Let $A$ be an irreducible row stochastic matrix. 
%     Then, $1$ is a simple eigenvalue of $A$. 
%     For any other eigenvalue $\lambda$ of $A$, we have that $|\lambda|\leq 1$.
%     If $A$ is aperiodic, then $|\lambda|<1$. 
%     If $A$ is periodic with period $\delta\geq 2$, then there are $\delta$ eigenvalues with absolute value equal to $1$, all distinct, and they are 
%     \[
%     \lambda_j=e^{2\pi i j/\delta},\,j=0,\ldots,\delta-1,
%     \]
%     where $i=\sqrt{-1}$.\hfill$\circ$
% \end{theorem}
% }
% }
% }
% \vspace{1mm}

In the next result, we show that any matrix with the same structure as the one of a matrix obtained from Algorithm~\ref{alg:buildAP} has period $\delta\geq 2$.

%{\color{blue}
Let $A$ be non-negative. Fix an index $i$ and define the period of index $i$ to be the greatest common divisor of all natural numbers $m$ such that $(A_{i,i})^m > 0$. When $A$ is irreducible, the period of every index is the same and is called the period of $A$. 
In fact, when $A$ is irreducible, the period can be defined as the greatest common divisor of the lengths of the closed directed paths in $\mathcal G(A)$.

\vspace{1mm} 
{\centering
\allowdisplaybreaks
\noindent \colorbox{c2!18}{\parbox{0.98\textwidth}{
\begin{lemma}\label{prop:2copies}
    The construction of $A^\text{\textsf{P}}$, obtained from a matrix $A$ satisfying Proposition~\ref{th:consensus_reqs} with Algorithm~\ref{alg:buildAP} does not suffice to ensure consensus. 
    %Given any matrix $A$ satisfying Theorem~\ref{th:consensus_reqs}, if  $A^\text{\textsf{P}}$ is the matrix obtained with Algorithm~\ref{alg:buildAP}, then the matrix $A^\text{\textsf{P}}$ can be periodic. \hfill$\circ$
\end{lemma}
 }
 }
 }
 \vspace{1mm}
 
\begin{proof}
     Consider the matrix 
     \[
     A=
     \begin{bmatrix}
        0 & a_{1,2} \\
        a_{1,2} & 0
     \end{bmatrix}.
     \]
     The matrix obtained with Algorithm~\ref{alg:buildAP} is
     \[
     A^\text{\textsf{P}}=
     \begin{bmatrix}
        0 & a_{1,2} & \rvline & a_{1,3}  & a_{1,4} & 0       & 0\\
        a_{1,2} & 0 & \rvline & 0       & 0       & a_{2,5} & a_{2,6}\\ \hline
        a_{3,1} & 0 & \rvline& 0       & 0       & 0       & 0\\
        a_{4,1} & 0 & \rvline & 0       & 0       & 0       & 0\\
        0       & a_{5,2}  & \rvline   & 0       & 0       & 0       & 0\\
        0       & a_{6,2}  & \rvline   & 0       & 0       & 0       & 0\\
     \end{bmatrix},
     \]
     with digraph representation $\mathcal G(A^\text{\textsf{P}})$ depicted in Figure~\ref{fig:graph_1}.

\begin{figure}[ht!]     
     \begin{center}
     \begin{tikzpicture}[scale=.6, transform shape,node distance=1.5cm]
\begin{scope}[every node/.style={circle,thick,draw},square/.style={regular polygon,regular polygon sides=4}]
\node[fill=green!20] (2) at (5.08148,2.12398) {\large $2$};
\node[fill=green!20] (1) at (1.89803,2.12458) {\large $1$};
\node[draw=gray,fill=yellow!50] (7) at (6.95,0.869942) {\large $6$};
\node[draw=gray,fill=yellow!50] (6) at (6.95,3.3791) {\large $5$};
\node[draw=gray,fill=yellow!50] (4) at (0.,3.3791) {\large $3$};
\node[draw=gray,fill=yellow!50] (5) at (0.,0.869942) {\large $4$};
\end{scope}
\begin{scope}[>={Stealth[black]},
              every edge/.style={draw=black, thick}]
\path [->] (1) edge[bend right=15] node {} (2);
\path [->] (2) edge[bend right=15] node {} (1);
\path [->] (4) edge[>={Stealth[gray!70]},gray!70,bend right=15] node {} (1);
\path [->] (1) edge[>={Stealth[gray!70]},gray!70,bend right=15] node {} (4);
\path [->] (5) edge[>={Stealth[gray!70]},gray!70,bend right=15] node {} (1);
\path [->] (1) edge[>={Stealth[gray!70]},gray!70,bend right=15] node {} (5);
\path [->] (6) edge[>={Stealth[gray!70]},gray!70,bend right=15] node {} (2);
\path [->] (2) edge[>={Stealth[gray!70]},gray!70,bend right=15] node {} (6);
\path [->] (7) edge[>={Stealth[gray!70]},gray!70,bend right=15] node {} (2);
\path [->] (2) edge[>={Stealth[gray!70]},gray!70,bend right=15] node {} (7);
\end{scope}
\end{tikzpicture}
\caption{Digraph representation of $A^\text{\textsf{P}}$, where the yellow nodes are local augmented states of the agents and the gray edges are the extra edges between agents and their augmented states.}
\label{fig:graph_1}
\end{center}
\end{figure}
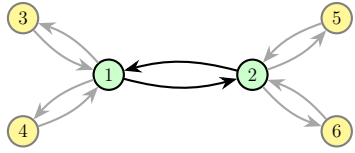
It is easy to observe that the closed directed paths have sizes multiple of 2. Hence, the matrix is periodic, with period $\delta = 2$. 
Hence, by the Perron-Frobenius theorem, there is more than one eigenvalue of $A^\text{\textsf{P}}$ with magnitude 1, and the system with this matrix does not reach consensus. 
\end{proof}

We further notice that, in the previous example, even if the agents were able to reach a consensus, an agent can infer the initial value of the other agent because it knows the values of its augmented states and the average consensus value. Hence, an agent can recover the sum of the other agent augmented states initial values, which is the initial value of the other agent.  

% \vspace{1mm} 
% {\centering
% \allowdisplaybreaks
% \noindent \colorbox{c4!15}{\parbox{0.98\textwidth}{
% \begin{corollary}\label{cor:2copies_not_work}
%     Given a matrix $A$ satisfying Theorem~\ref{th:consensus_reqs}, if  $A^\text{\textsf{P}}$ is the normalized to be row stochastic matrix resulting from the matrix obtained with Algorithm~\ref{alg:buildAP}, then the LTI system described in Eq.~\eqref{eq:LTI_consensus} using the matrix $A^\text{\textsf{P}}$ may not reach asymptotic consensus.\hfill$\circ$ 
% \end{corollary}
% }
% }
% }
% \vspace{1mm}

% \begin{proof}
%     From Proposition~\ref{prop:2copies}, we know that $A^\text{\textsf{P}}$ may be periodic. 
%     Therefore, by the Perron-Frobenius theorem, there is more than one eigenvalue of $A^\text{\textsf{P}}$ with absolute value equal to $1$. 
%     Finally, the condition $(1)$ of Theorem~\ref{th:consensus_reqs} is not satisfied and the LTI system does not reach consensus. 
% \end{proof}

% \vspace{1mm} 
% {\centering
% \allowdisplaybreaks
% \noindent \colorbox{gray!15}{\parbox{0.98\textwidth}{
% \noindent$\mathbf{Q_1}$ The question that follows is the following: \emph{``Since two copies of the nodes do not suffice to ensure privacy in average consensus, what if we use three copies of the nodes?''}\hfill$\triangleleft$
% }
% }
% }
% \vspace{1mm}

%\subsubsection{Each agent with three copies}

\subsubsection{Solution to Problem $\mathbf{P_1}$}

%To address $\mathbf{Q_1}$, 
The first step is to overcome the limitations previously encountered. 
Therefore, we need to verify what are the augmented state networks that achieve the desired goal. 
We summarize this result in the following.

\vspace{1mm} 
{\centering
\allowdisplaybreaks
\noindent \colorbox{c2!18}{\parbox{0.98\textwidth}{
\begin{lemma}\label{prop:enumerate_3cpy_}
The only networks (up to isomorphism) with three augmented states that satisfy the necessary conditions to address $\mathbf{P_1}$ are the ones depicted in Figure~\ref{fig:3copies_that_work}. 
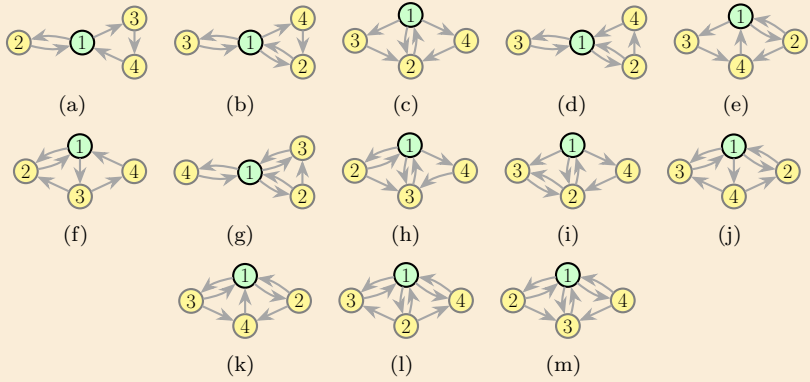
\begin{figure}[H]
\centering
\subfigure[]{
\begin{tikzpicture}[scale=.31, transform shape,node distance=1.5cm]
\begin{scope}[every node/.style={circle,thick,draw},square/.style={regular polygon,regular polygon sides=4}]
\node[fill=green!20] (1) at (2.72064,1.03433) {\Huge $1$};
\node[draw=gray,fill=yellow!50] (2) at (0.,1.03407) {\Huge $2$};
\node[draw=gray,fill=yellow!50] (3) at (4.96021,2.06811) {\Huge $3$};
\node[draw=gray,fill=yellow!50] (4) at (4.96004,0.) {\Huge $4$};
\end{scope}
\begin{scope}[>={Stealth[black]},
              every edge/.style={draw=black, thick}]
\path [->] (1) edge[>={Stealth[gray!70]},gray!70,bend right=15] node {} (2);
\path [->] (1) edge[>={Stealth[gray!70]},gray!70] node {} (3);
\path [->] (2) edge[>={Stealth[gray!70]},gray!70,bend right=15] node {} (1);
\path [->] (3) edge[>={Stealth[gray!70]},gray!70] node {} (4);
\path [->] (4) edge[>={Stealth[gray!70]},gray!70] node {} (1);
\end{scope}
\end{tikzpicture}
}
\subfigure[]{
\begin{tikzpicture}[scale=.31, transform shape,node distance=1.5cm]
\begin{scope}[every node/.style={circle,thick,draw},square/.style={regular polygon,regular polygon sides=4}]
\node[fill=green!20] (1) at (2.72645,1.03698) {\Huge $1$};
\node[draw=gray,fill=yellow!50] (2) at (4.97051,0.) {\Huge $2$};
\node[draw=gray,fill=yellow!50] (3) at (0.,1.03674) {\Huge $3$};
\node[draw=gray,fill=yellow!50] (4) at (4.97146,2.07299) {\Huge $4$};
\end{scope}
\begin{scope}[>={Stealth[black]},
              every edge/.style={draw=black, thick}]
\path [->] (1) edge[>={Stealth[gray!70]},gray!70,bend right=15] node {} (2);
\path [->] (1) edge[>={Stealth[gray!70]},gray!70,bend right=15] node {} (3);
\path [->] (1) edge[>={Stealth[gray!70]},gray!70] node {} (4);
\path [->] (2) edge[>={Stealth[gray!70]},gray!70,bend right=15] node {} (1);
\path [->] (3) edge[>={Stealth[gray!70]},gray!70,bend right=15] node {} (1);
\path [->] (4) edge[>={Stealth[gray!70]},gray!70] node {} (2);
\end{scope}
\end{tikzpicture}
}
\subfigure[]{
\begin{tikzpicture}[scale=.31, transform shape,node distance=1.5cm]
\begin{scope}[every node/.style={circle,thick,draw},square/.style={regular polygon,regular polygon sides=4}]
\node[fill=green!20] (1) at (2.3835,2.2193) {\Huge $1$};
\node[draw=gray,fill=yellow!50] (2) at (2.38428,0.) {\Huge $2$};
\node[draw=gray,fill=yellow!50] (3) at (0.,1.11038) {\Huge $3$};
\node[draw=gray,fill=yellow!50] (4) at (4.76698,1.11074) {\Huge $4$};
\end{scope}
\begin{scope}[>={Stealth[black]},
              every edge/.style={draw=black, thick}]
\path [->] (1) edge[>={Stealth[gray!70]},gray!70,bend right=15] node {} (2);
\path [->] (1) edge[>={Stealth[gray!70]},gray!70] node {} (3);
\path [->] (1) edge[>={Stealth[gray!70]},gray!70] node {} (4);
\path [->] (2) edge[>={Stealth[gray!70]},gray!70,bend right=15] node {} (1);
\path [->] (3) edge[>={Stealth[gray!70]},gray!70] node {} (2);
\path [->] (4) edge[>={Stealth[gray!70]},gray!70] node {} (2);
\end{scope}
\end{tikzpicture}
}
\subfigure[]{
\begin{tikzpicture}[scale=.31, transform shape,node distance=1.5cm]
\begin{scope}[every node/.style={circle,thick,draw},square/.style={regular polygon,regular polygon sides=4}]
\node[fill=green!20] (1) at (2.72645,1.03698) {\Huge $1$};
\node[draw=gray,fill=yellow!50] (2) at (4.97051,0.) {\Huge $2$};
\node[draw=gray,fill=yellow!50] (3) at (0.,1.03674) {\Huge $3$};
\node[draw=gray,fill=yellow!50] (4) at (4.97146,2.07299) {\Huge $4$};
\end{scope}
\begin{scope}[>={Stealth[black]},
              every edge/.style={draw=black, thick}]
\path [->] (1) edge[>={Stealth[gray!70]},gray!70,bend right=15] node {} (2);
\path [->] (1) edge[>={Stealth[gray!70]},gray!70,bend right=15] node {} (3);
\path [->] (2) edge[>={Stealth[gray!70]},gray!70,bend right=15] node {} (1);
\path [->] (2) edge[>={Stealth[gray!70]},gray!70] node {} (4);
\path [->] (3) edge[>={Stealth[gray!70]},gray!70,bend right=15] node {} (1);
\path [->] (4) edge[>={Stealth[gray!70]},gray!70] node {} (1);
\end{scope}
\end{tikzpicture}
}
\subfigure[]{
\begin{tikzpicture}[scale=.31, transform shape,node distance=1.5cm]
\begin{scope}[every node/.style={circle,thick,draw},square/.style={regular polygon,regular polygon sides=4}]
\node[fill=green!20] (1) at (2.32149,2.16099) {\Huge $1$};
\node[draw=gray,fill=yellow!50] (2) at (4.64645,1.08153) {\Huge $2$};
\node[draw=gray,fill=yellow!50] (3) at (0.,1.07938) {\Huge $3$};
\node[draw=gray,fill=yellow!50] (4) at (2.32612,0.) {\Huge $4$};
\end{scope}
\begin{scope}[>={Stealth[black]},
              every edge/.style={draw=black, thick}]
\path [->] (1) edge[>={Stealth[gray!70]},gray!70,bend right=15] node {} (2);
\path [->] (1) edge[>={Stealth[gray!70]},gray!70] node {} (3);
\path [->] (2) edge[>={Stealth[gray!70]},gray!70,bend right=15] node {} (1);
\path [->] (2) edge[>={Stealth[gray!70]},gray!70] node {} (4);
\path [->] (3) edge[>={Stealth[gray!70]},gray!70] node {} (4);
\path [->] (4) edge[>={Stealth[gray!70]},gray!70] node {} (1);
\end{scope}
\end{tikzpicture}
}
\subfigure[]{
\begin{tikzpicture}[scale=.31, transform shape,node distance=1.5cm]
\begin{scope}[every node/.style={circle,thick,draw},square/.style={regular polygon,regular polygon sides=4}]
\node[fill=green!20] (1) at (2.32091,2.16192) {\Huge $1$};
\node[draw=gray,fill=yellow!50] (2) at (0.,1.08168) {\Huge $2$};
\node[draw=gray,fill=yellow!50] (3) at (2.32351,0.) {\Huge $3$};
\node[draw=gray,fill=yellow!50] (4) at (4.64744,1.0829) {\Huge $4$};
\end{scope}
\begin{scope}[>={Stealth[black]},
              every edge/.style={draw=black, thick}]
\path [->] (1) edge[>={Stealth[gray!70]},gray!70,bend right=15] node {} (2);
\path [->] (1) edge[>={Stealth[gray!70]},gray!70] node {} (3);
\path [->] (2) edge[>={Stealth[gray!70]},gray!70,bend right=15] node {} (1);
\path [->] (3) edge[>={Stealth[gray!70]},gray!70] node {} (2);
\path [->] (3) edge[>={Stealth[gray!70]},gray!70] node {} (4);
\path [->] (4) edge[>={Stealth[gray!70]},gray!70] node {} (1);
\end{scope}
\end{tikzpicture}
}
\subfigure[]{
\begin{tikzpicture}[scale=.31, transform shape,node distance=1.5cm]
\begin{scope}[every node/.style={circle,thick,draw},square/.style={regular polygon,regular polygon sides=4}]
\node[fill=green!20] (1) at (2.73068,1.03786) {\Huge $1$};
\node[draw=gray,fill=yellow!50] (2) at (4.97888,0.) {\Huge $2$};
\node[draw=gray,fill=yellow!50] (3) at (4.97869,2.07617) {\Huge $3$};
\node[draw=gray,fill=yellow!50] (4) at (0.,1.03805) {\Huge $4$};
\end{scope}
\begin{scope}[>={Stealth[black]},
              every edge/.style={draw=black, thick}]
\path [->] (1) edge[>={Stealth[gray!70]},gray!70,bend right=15] node {} (2);
\path [->] (1) edge[>={Stealth[gray!70]},gray!70,bend right=15] node {} (3);
\path [->] (1) edge[>={Stealth[gray!70]},gray!70,bend right=15] node {} (4);
\path [->] (2) edge[>={Stealth[gray!70]},gray!70,bend right=15] node {} (1);
\path [->] (2) edge[>={Stealth[gray!70]},gray!70] node {} (3);
\path [->] (3) edge[>={Stealth[gray!70]},gray!70,bend right=15] node {} (1);
\path [->] (4) edge[>={Stealth[gray!70]},gray!70,bend right=15] node {} (1);
\end{scope}
\end{tikzpicture}
}
\subfigure[]{
\begin{tikzpicture}[scale=.31, transform shape,node distance=1.5cm]
\begin{scope}[every node/.style={circle,thick,draw},square/.style={regular polygon,regular polygon sides=4}]
\node[fill=green!20] (1) at (2.36388,2.20194) {\Huge $1$};
\node[draw=gray,fill=yellow!50] (2) at (0.,1.10171) {\Huge $2$};
\node[draw=gray,fill=yellow!50] (3) at (2.36653,0.) {\Huge $3$};
\node[draw=gray,fill=yellow!50] (4) at (4.73348,1.10295) {\Huge $4$};
\end{scope}
\begin{scope}[>={Stealth[black]},
              every edge/.style={draw=black, thick}]
\path [->] (1) edge[>={Stealth[gray!70]},gray!70,bend right=15] node {} (2);
\path [->] (1) edge[>={Stealth[gray!70]},gray!70,bend right=15] node {} (3);
\path [->] (1) edge[>={Stealth[gray!70]},gray!70] node {} (4);
\path [->] (2) edge[>={Stealth[gray!70]},gray!70,bend right=15] node {} (1);
\path [->] (2) edge[>={Stealth[gray!70]},gray!70] node {} (3);
\path [->] (3) edge[>={Stealth[gray!70]},gray!70,bend right=15] node {} (1);
\path [->] (4) edge[>={Stealth[gray!70]},gray!70,bend right=15] node {} (3);
\end{scope}
\end{tikzpicture}
}
\subfigure[]{
\begin{tikzpicture}[scale=.31, transform shape,node distance=1.5cm]
\begin{scope}[every node/.style={circle,thick,draw},square/.style={regular polygon,regular polygon sides=4}]
\node[fill=green!20] (1) at (2.36551,2.20256) {\Huge $1$};
\node[draw=gray,fill=yellow!50] (2) at (2.36629,0.) {\Huge $2$};
\node[draw=gray,fill=yellow!50] (3) at (0.,1.102) {\Huge $3$};
\node[draw=gray,fill=yellow!50] (4) at (4.73101,1.10236) {\Huge $4$};
\end{scope}
\begin{scope}[>={Stealth[black]},
              every edge/.style={draw=black, thick}]
\path [->] (1) edge[>={Stealth[gray!70]},gray!70,bend right=15] node {} (2);
\path [->] (1) edge[>={Stealth[gray!70]},gray!70] node {} (3);
\path [->] (1) edge[>={Stealth[gray!70]},gray!70] node {} (4);
\path [->] (2) edge[>={Stealth[gray!70]},gray!70,bend right=15] node {} (1);
\path [->] (2) edge[>={Stealth[gray!70]},gray!70,bend right=15] node {} (3);
\path [->] (3) edge[>={Stealth[gray!70]},gray!70,bend right=15] node {} (2);
\path [->] (4) edge[>={Stealth[gray!70]},gray!70] node {} (2);
\end{scope}
\end{tikzpicture}
}
\subfigure[]{
\begin{tikzpicture}[scale=.31, transform shape,node distance=1.5cm]
\begin{scope}[every node/.style={circle,thick,draw},square/.style={regular polygon,regular polygon sides=4}]
\node[fill=green!20] (1) at (2.31333,2.1534) {\Huge $1$};
\node[draw=gray,fill=yellow!50] (2) at (4.63012,1.07772) {\Huge $2$};
\node[draw=gray,fill=yellow!50] (3) at (0.,1.07558) {\Huge $3$};
\node[draw=gray,fill=yellow!50] (4) at (2.31794,0.) {\Huge $4$};
\end{scope}
\begin{scope}[>={Stealth[black]},
              every edge/.style={draw=black, thick}]
\path [->] (1) edge[>={Stealth[gray!70]},gray!70,bend right=15] node {} (2);
\path [->] (1) edge[>={Stealth[gray!70]},gray!70,bend right=15] node {} (3);
\path [->] (1) edge[>={Stealth[gray!70]},gray!70] node {} (4);
\path [->] (2) edge[>={Stealth[gray!70]},gray!70,bend right=15] node {} (1);
\path [->] (3) edge[>={Stealth[gray!70]},gray!70,bend right=15] node {} (1);
\path [->] (4) edge[>={Stealth[gray!70]},gray!70] node {} (2);
\path [->] (4) edge[>={Stealth[gray!70]},gray!70] node {} (3);
\end{scope}
\end{tikzpicture}
}
\subfigure[]{
\begin{tikzpicture}[scale=.31, transform shape,node distance=1.5cm]
\begin{scope}[every node/.style={circle,thick,draw},square/.style={regular polygon,regular polygon sides=4}]
\node[fill=green!20] (1) at (2.31333,2.1534) {\Huge $1$};
\node[draw=gray,fill=yellow!50] (2) at (4.63012,1.07772) {\Huge $2$};
\node[draw=gray,fill=yellow!50] (3) at (0.,1.07558) {\Huge $3$};
\node[draw=gray,fill=yellow!50] (4) at (2.31794,0.) {\Huge $4$};
\end{scope}
\begin{scope}[>={Stealth[black]},
              every edge/.style={draw=black, thick}]
\path [->] (1) edge[>={Stealth[gray!70]},gray!70,bend right=15] node {} (2);
\path [->] (1) edge[>={Stealth[gray!70]},gray!70,bend right=15] node {} (3);
\path [->] (2) edge[>={Stealth[gray!70]},gray!70,bend right=15] node {} (1);
\path [->] (2) edge[>={Stealth[gray!70]},gray!70] node {} (4);
\path [->] (3) edge[>={Stealth[gray!70]},gray!70,bend right=15] node {} (1);
\path [->] (3) edge[>={Stealth[gray!70]},gray!70] node {} (4);
\path [->] (4) edge[>={Stealth[gray!70]},gray!70] node {} (1);
\end{scope}
\end{tikzpicture}
}
\subfigure[]{
\begin{tikzpicture}[scale=.31, transform shape,node distance=1.5cm]
\begin{scope}[every node/.style={circle,thick,draw},square/.style={regular polygon,regular polygon sides=4}]
\node[fill=green!20] (1) at (2.35263,2.19056) {\Huge $1$};
\node[draw=gray,fill=yellow!50] (2) at (2.3534,0.) {\Huge $2$};
\node[draw=gray,fill=yellow!50] (3) at (0.,1.096) {\Huge $3$};
\node[draw=gray,fill=yellow!50] (4) at (4.70525,1.09636) {\Huge $4$};
\end{scope}
\begin{scope}[>={Stealth[black]},
              every edge/.style={draw=black, thick}]
\path [->] (1) edge[>={Stealth[gray!70]},gray!70,bend right=15] node {} (2);
\path [->] (1) edge[>={Stealth[gray!70]},gray!70,bend right=15] node {} (3);
\path [->] (1) edge[>={Stealth[gray!70]},gray!70,bend right=15] node {} (4);
\path [->] (2) edge[>={Stealth[gray!70]},gray!70,bend right=15] node {} (1);
\path [->] (2) edge[>={Stealth[gray!70]},gray!70] node {} (3);
\path [->] (2) edge[>={Stealth[gray!70]},gray!70] node {} (4);
\path [->] (3) edge[>={Stealth[gray!70]},gray!70,bend right=15] node {} (1);
\path [->] (4) edge[>={Stealth[gray!70]},gray!70,bend right=15] node {} (1);
\end{scope}
\end{tikzpicture}
}
\subfigure[]{
\begin{tikzpicture}[scale=.31, transform shape,node distance=1.5cm]
\begin{scope}[every node/.style={circle,thick,draw},square/.style={regular polygon,regular polygon sides=4}]
\node[fill=green!20] (1) at (2.35045,2.18944) {\Huge $1$};
\node[draw=gray,fill=yellow!50] (2) at (0.,1.09545) {\Huge $2$};
\node[draw=gray,fill=yellow!50] (3) at (2.35309,0.) {\Huge $3$};
\node[draw=gray,fill=yellow!50] (4) at (4.7066,1.09668) {\Huge $4$};
\end{scope}
\begin{scope}[>={Stealth[black]},
              every edge/.style={draw=black, thick}]
\path [->] (1) edge[>={Stealth[gray!70]},gray!70,bend right=15] node {} (2);
\path [->] (1) edge[>={Stealth[gray!70]},gray!70,bend right=15] node {} (3);
\path [->] (1) edge[>={Stealth[gray!70]},gray!70,bend right=15] node {} (4);
\path [->] (2) edge[>={Stealth[gray!70]},gray!70,bend right=15] node {} (1);
\path [->] (2) edge[>={Stealth[gray!70]},gray!70] node {} (3);
\path [->] (3) edge[>={Stealth[gray!70]},gray!70,bend right=15] node {} (1);
\path [->] (4) edge[>={Stealth[gray!70]},gray!70,bend right=15] node {} (1);
\path [->] (4) edge[>={Stealth[gray!70]},gray!70] node {} (3);
\end{scope}
\end{tikzpicture}
}
\caption{Networks with one agent (node 1) and three augmented states that can be used to address problem $\mathbf{P_1}$.\hfill$\circ$}
\label{fig:3copies_that_work}
\end{figure}%\hfill$\circ$
\end{lemma}
 }
 }
 }
 \vspace{1mm}

\begin{proof}
We start by enumerating all the strongly connected networks with 4 agents (3 augmented states).
Second, we filtered the obtained networks that satisfy the following:
\begin{itemize}
    \item do not correspond to periodic matrices;
    \item are not observable by sensing only node 1 (the original node); 
    \item allow a parameterization (numerical instance of the nonzeros) of the adjacency matrix which ensures that an agent cannot recover the initial state of at least one of the other agents augmented states. 
\end{itemize}
Last, we filtered the later list to have only one graph for each class of isomorphic graphs, obtaining the networks depicted in Figure~\ref{fig:3copies_that_work}. 
\end{proof}

Subsequently, without loss of generality, we focus on the first network of augmented states, Figure~\ref{fig:3copies_that_work}~(a), that allow to address problem $\mathbf{P_1}$, and we provide an algorithm that selects a parameterization of the nonzero values of the adjacency matrix to illustrate how we can ensure privacy in average consensus.

% (i -> (n + 3*(i - 1) + 1)),
% (i -> (n + 3*(i - 1) + 2)),
% ((n + 3*(i - 1) + 1) -> i),
% ((n + 3*(i - 1) + 2) -> (n + 3*(i - 1) + 3)),
% ((n + 3*(i - 1) + 3) -> i)

\begin{algorithm}[!ht]
		{%\small
		\caption{Construction of the privacy dynamics matrix $A^\text{\textsf{P}}$ for consensus}
		\label{alg:buildAP2}
		\begin{algorithmic}[1]
			\STATE{\textbf{input}: network of agents $\mathcal G=\langle \mathcal X,\mathcal E_{\mathcal X,\mathcal X}\rangle$ and initial agents state $x[0]\in\mathbb R^N$}
			\STATE{\textbf{output}: dynamics matrix $A^\text{\textsf{P}}\in\mathbb R^{4N\times 4N}$, and initial agents state $\tilde x^{\text{\textsf{P}}}[0]\in\mathbb R^{4N}$}
		\STATE{\textbf{set} $A\in\mathbb R^{N\times N}$ as the adjacency matrix of $\mathcal G$}
		\STATE{\textbf{fill} the entries of $A^{\text{\textsf{P}}}$ with zeros}
		\STATE{\textbf{set}
		\hfill{\color{gray}$\rhd$ \underline{Copy the matrix $A$ to the first $n$ rows and $n$ columns}}
		$$A^{\text{\textsf{P}}}_{ij}=A_{ij}, \forall i,j\in\{1,\ldots,N\}\text{ and }i\neq j$$}
		\STATE{\textbf{for} $i=1,\ldots,N$\hfill{\color{gray}$\rhd$ \underline{Set additional entries values}}
		$$
		\begin{array}{ll}A^{\text{\textsf{P}}}_{i, N + 3i - 2} = 1, & A^{\text{\textsf{P}}}_{N + 3i - 2, 1} = 2,
		\\
		A^{\text{\textsf{P}}}_{i, n + 3i - 1} = 1, &
		A^{\text{\textsf{P}}}_{n + 3i - 1, n + 3i - 0} = 1, 
		\\
		A^{\text{\textsf{P}}}_{n + 3i - 0, i} = 1 & 
		\end{array}
		$$
		}
		\STATE{\textbf{normalize} the rows of $A^{\text{\textsf{P}}}$ dividing by their sum}
		\STATE{\textbf{compute} the left-eigenvector $v^0$ of $A^{\text{\textsf{P}}}$ associated with the eigenvalue 1}
		\STATE{\textbf{distribute} the initial condition of each agent across its augmented states (setting as $0$ for the original state), the distribution can be tailored by each agent, and scaling it to make the new average of initial states equal to the original one, i.e., select $\alpha_i, \beta_i, \gamma_i\neq 0$
\[
\begin{split}
\begin{bmatrix}\tilde x_i^{\text{\textsf{P}}}[0] & \tilde x_{N+3i-2}^{\text{\textsf{P}}}[0] & \tilde x_{N+3i-1}^{\text{\textsf{P}}}[0] & \tilde x_{N+3i-0}^{\text{\textsf{P}}}[0]\end{bmatrix} =\\
\frac{4 x_i^{\text{\textsf{P}}}[0]}{\alpha_i+ \beta_i+ \gamma_i} \begin{bmatrix}0 & \alpha_i & \beta_i & \gamma_i\end{bmatrix}.
\end{split}
\] 
}
\STATE{\textbf{re-scale}, for each agent the new initial condition, $\tilde x_j^{\text{\textsf{P}}}[0]$, to be 
\[
\tilde{x}_j^{\text{\textsf{P}}}[0]=\frac{x_j^{\text{\textsf{P}}}[0]}{4N v^0_j}\sum_{i=1}^{4N}v^0_i.
\]
}
		\end{algorithmic} 
		}
\end{algorithm} 

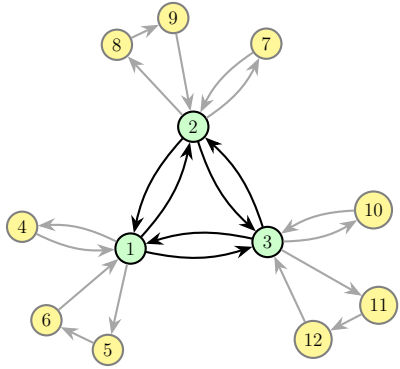
\begin{figure}[ht!]
\begin{center}
\begin{tikzpicture}[scale=.6, transform shape,node distance=1.5cm]
\begin{scope}[every node/.style={circle,thick,draw},square/.style={regular polygon,regular polygon sides=4}]
\node[fill=green!20] (1) at (2.40241,2.24721) {\large $1$};
\node[fill=green!20] (2) at (3.79239,4.94879) {\large $2$};
\node[fill=green!20] (3) at (5.43773,2.39296) {\large $3$};
\node[draw=gray,fill=yellow!50] (4) at (0.,2.73031) {\large $4$};
\node[draw=gray,fill=yellow!50] (5) at (1.9,0.) {\large $5$};
\node[draw=gray,fill=yellow!50] (6) at (0.538332,0.657477) {\large $6$};
\node[draw=gray,fill=yellow!50] (7) at (5.41037,6.79029) {\large $7$};
\node[draw=gray,fill=yellow!50] (8) at (2.1,6.76612) {\large $8$};
\node[draw=gray,fill=yellow!50] (9) at (3.34898,7.35913) {\large $9$};
\node[draw=gray,fill=yellow!50] (10) at (7.78551,3.10302) {\large $10$};
\node[draw=gray,fill=yellow!50] (11) at (7.9,0.98889) {\large $11$};
\node[draw=gray,fill=yellow!50] (12) at (6.44872,0.232587) {\large $12$};
\end{scope}
\begin{scope}[>={Stealth[black]},
              every edge/.style={draw=black, thick}]
\path [->] (1) edge[bend right=15] node {} (2);
\path [->] (2) edge[bend right=15] node {} (1);
\path [->] (2) edge[bend right=15] node {} (3);
\path [->] (3) edge[bend right=15] node {} (2);
\path [->] (1) edge[bend right=15] node {} (3);
\path [->] (3) edge[bend right=15] node {} (1);
\path [->] (1) edge[>={Stealth[gray!70]},gray!70,bend right=15] node {} (4);
\path [->] (1) edge[>={Stealth[gray!70]},gray!70] node {} (5);
\path [->] (4) edge[>={Stealth[gray!70]},gray!70,bend right=15] node {} (1);
\path [->] (5) edge[>={Stealth[gray!70]},gray!70] node {} (6);
\path [->] (6) edge[>={Stealth[gray!70]},gray!70] node {} (1);
\path [->] (2) edge[>={Stealth[gray!70]},gray!70,bend right=15] node {} (7);
\path [->] (2) edge[>={Stealth[gray!70]},gray!70] node {} (8);
\path [->] (7) edge[>={Stealth[gray!70]},gray!70,bend right=15] node {} (2);
\path [->] (8) edge[>={Stealth[gray!70]},gray!70] node {} (9);
\path [->] (9) edge[>={Stealth[gray!70]},gray!70] node {} (2);
\path [->] (3) edge[>={Stealth[gray!70]},gray!70,bend right=15] node {} (10);
\path [->] (3) edge[>={Stealth[gray!70]},gray!70] node {} (11);
\path [->] (10) edge[>={Stealth[gray!70]},gray!70,bend right=15] node {} (3);
\path [->] (11) edge[>={Stealth[gray!70]},gray!70] node {} (12);
\path [->] (12) edge[>={Stealth[gray!70]},gray!70] node {} (3);
\end{scope}
\end{tikzpicture}
\caption{Digraph representation of $A^{\text{\textsf{P}}}$ resulting from Algorithm~\ref{alg:buildAP2} for a bidirected cycle network with 3 agents.}
\label{fig:cycle3_3copies}
\end{center}
\end{figure}

Now, we use the running example of a bidirected cycle with 3 agents with Algorithm~\ref{alg:buildAP2}. 
We obtain the network of agents depicted in Figure~\ref{fig:cycle3_3copies} and the matrix $A^{\text{\textsf{P}}}$ is the following:

\[
\begin{array}{l}
A^{\textsf{P}}  =  
{\small
\begin{bmatrix}
 0 & \frac{1}{5} & \frac{1}{5} & \rvline & \frac{2}{5} & 0 & \frac{1}{5} & 0 & 0 & 0 & 0 & 0 & 0 \\
 \frac{1}{5} & 0 & \frac{1}{5} &\rvline & 0 & 0 & 0 & \frac{2}{5} & 0 & \frac{1}{5} & 0 & 0 & 0 \\
 \frac{1}{5} & \frac{1}{5} & 0 &\rvline & 0 & 0 & 0 & 0 & 0 & 0 & \frac{2}{5} & 0 & \frac{1}{5} \\ \hline 
 1 & 0 & 0 &\rvline & 0 & 0 & 0 & 0 & 0 & 0 & 0 & 0 & 0 \\
 1 & 0 & 0 &\rvline & 0 & 0 & 0 & 0 & 0 & 0 & 0 & 0 & 0 \\
 0 & 0 & 0 &\rvline & 0 & 1 & 0 & 0 & 0 & 0 & 0 & 0 & 0 \\
 0 & 1 & 0 &\rvline & 0 & 0 & 0 & 0 & 0 & 0 & 0 & 0 & 0 \\
 0 & 1 & 0 &\rvline & 0 & 0 & 0 & 0 & 0 & 0 & 0 & 0 & 0 \\
 0 & 0 & 0 &\rvline & 0 & 0 & 0 & 0 & 1 & 0 & 0 & 0 & 0 \\
 0 & 0 & 1 &\rvline & 0 & 0 & 0 & 0 & 0 & 0 & 0 & 0 & 0 \\
 0 & 0 & 1 &\rvline & 0 & 0 & 0 & 0 & 0 & 0 & 0 & 0 & 0 \\
 0 & 0 & 0 &\rvline & 0 & 0 & 0 & 0 & 0 & 0 & 0 & 1 & 0 \\
\end{bmatrix}
}
\end{array}.
\]

The network of agents depicted in Figure~\ref{fig:cycle3_3copies} is a repeated pattern for each agent. 
Hence, we can check what agent 1 can recover by observing its direct neighbors, i.e., with 
\[
C =
\begin{bmatrix}
\bigI_6 & \rvline  & \bigzero_{6\times 6} \\
\end{bmatrix}.
\]

Using Theorem~\ref{th:PBH}, we have that 
\[
span\left(P_{O(A^\text{\textsf{P}},C)}^0\right)=
span\left(\begin{bmatrix}
\bigI_6 & \rvline  & \bigzero_{3\times 6} \\ \hline
\bigzero_{6\times 6} & \rvline & \bigR
\end{bmatrix}
\right),
\]
with 

\[
R=
\begin{bmatrix}
 1 & 0 & 0 & 0 & 0 & 0 & 0 & 0 & 0 \\
 0 & 1 & 0 & 0 & 0 & 0 & 0 & 0 & 0 \\
 0 & 0 & 1 & 0 & 0 & 0 & 0 & 0 & 0 \\
 0 & 0 & 0 & 1 & 0 & \frac{1}{2} & 0 & 0 & 0 \\
 0 & 0 & 0 & 0 & 1 & 0 & 0 & 0 & 0 \\
 0 & 0 & 0 & 0 & 0 & 0 & 1 & 0 & \frac{1}{2} \\
 0 & 0 & 0 & 0 & 0 & 0 & 0 & 1 & 0 \\
\end{bmatrix}.
\]

Therefore, we have that $e_{12}^{10},e_{12}^{12}\notin span\left(P_{O(A^\text{\textsf{P}},C)}^0\right)$ and $e_{12}^{7}+e_{12}^{8}+e_{12}^{9}\notin span\left(P_{O(A^\text{\textsf{P}},C)}^0\right)$ and $e_{12}^{10}+e_{12}^{11}+e_{12}^{12}\notin span\left(P_{O(A^\text{\textsf{P}},C)}^0\right)$

In other words, agent 1 cannot recover the original initial states of agents 2 and 3. 

However, to reach average consensus, the initial state of the agents needs to be re-scaled according to the left-eigenvector associated with the eigenvalue 1, which is:
\[
v_L = \left[\begin{smallmatrix}
 \frac{5}{27} & \frac{5}{27} & \frac{5}{27} & \frac{2}{27} & \frac{1}{27} & \frac{1}{27} & \frac{2}{27} & \frac{1}{27} & \frac{1}{27} & \frac{2}{27} & \frac{1}{27} & \frac{1}{27} \\
\end{smallmatrix}\right]^\intercal.
\]

Hence, another property that we have to ensure, is that agent 1 is not able to recover a combination of other agents augmented states that is given by the corresponding values in~$v_L$. 

Regarding agent 2, the agent 1 can recover $e_{12}^{7}+e_{12}^{8}+\frac{1}{2}e_{12}^{9}$ with coefficients that are not multiple of the respective entries in $v_L$. 
Analogously, regarding agent 3, the agent 1 can recover $e_{12}^{10}+e_{12}^{11}+\frac{1}{2}e_{12}^{12}$ with coefficients that are not multiple of the respective entries in $v_L$. 

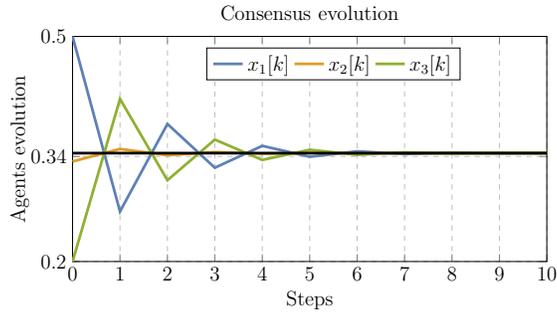
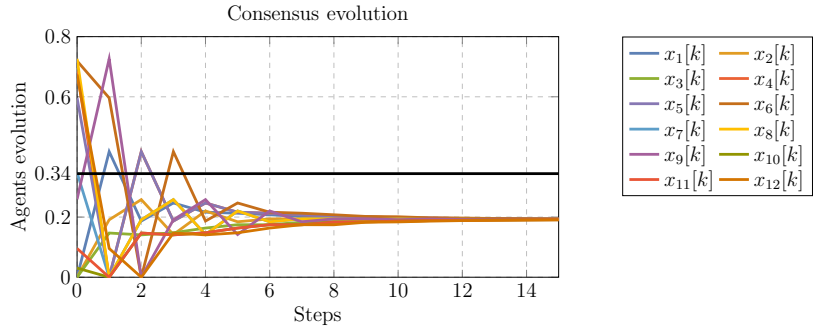
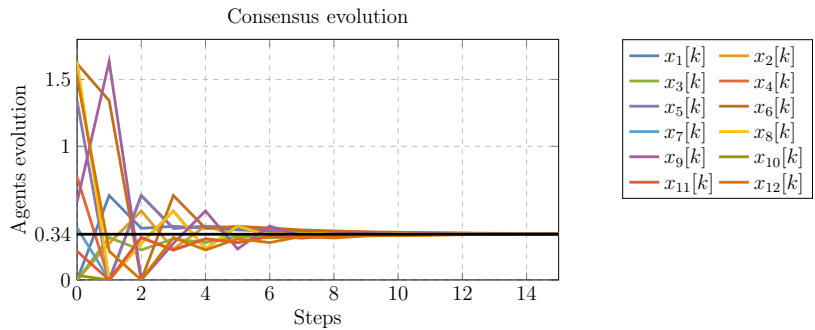
\begin{figure}[H]
\centering
\subfigure[
{State evolution of agents from network $\mathcal G(A)$ with weight matrix $A$ and initial states $x[0]$.}
]{
\tikzset{
  font={\fontsize{12pt}{10}\selectfont}}

\begin{tikzpicture}[scale=.65]
\begin{axis}[
    title={Consensus evolution},
    xlabel={Steps},
    ylabel={Agents evolution},
    ylabel shift=-6pt,
    xlabel shift=-1pt,
    xmin=0, xmax=10,
    ymin=0.2, ymax=0.5, %-3.15, 5.4
    ytick={0.2,0.34,0.5},
    yticklabel style={/pgf/number format/.cd,fixed,precision=2},
    legend pos=north west,
    ymajorgrids=true,
    xmajorgrids=true,
    grid style=dashed,
    % only scale the axis, not the axis including the ticks and labels
    scale only axis=true,
    % set `width' and `height' to the desired values
    width=0.8\columnwidth,
    height=0.38\columnwidth,
    legend style={at={(0.55,0.95)},anchor=north,legend cell align=left},
    legend columns=3,
    %transpose legend
]

\addplot[color=c1,ultra thick,]
coordinates {(0,0.5)(1,0.266667)(2,0.383333)(3,0.325)(4,0.354167)(5,0.339583)(6,0.346875)(7,0.343229)(8,0.345052)(9,0.344141)(10,0.344596)};
\addplot[color=c2,ultra thick,]
coordinates {(0,0.333333)(1,0.35)(2,0.341667)(3,0.345833)(4,0.34375)(5,0.344792)(6,0.344271)(7,0.344531)(8,0.344401)(9,0.344466)(10,0.344434)};
\addplot[color=c3,ultra thick,]
coordinates {(0,0.2)(1,0.416667)(2,0.308333)(3,0.3625)(4,0.335417)(5,0.348958)(6,0.342188)(7,0.345573)(8,0.34388)(9,0.344727)(10,0.344303)};

\addplot[color=black,ultra thick,]
coordinates {(0,0.344444)(10,0.344444)};

\legend{{ $x_{1}[k]$}, { $x_{2}[k]$}, { $x_{3}[k]$}}

\end{axis}
\end{tikzpicture}
}

\subfigure[{State evolution of agents from network $\mathcal G(A^\text{\textsf{P}})$ with weight matrix $A^\text{\textsf{P}}$ and initial states $x'[0]$.}]{
\tikzset{
  font={\fontsize{12pt}{10}\selectfont}}

\begin{tikzpicture}[scale=.66]
\begin{axis}[
    title={Consensus evolution},
    xlabel={Steps},
    ylabel={Agents evolution},
    ylabel shift=-4pt,
    xlabel shift=-1pt,
    xmin=0, xmax=15,
    ymin=0., ymax=0.8, %-3.15, 5.4
    ytick={0,0.2,0.344444,0.6,0.8},
    legend pos=north west,
    ymajorgrids=true,
    xmajorgrids=true,
    grid style=dashed,
    % only scale the axis, not the axis including the ticks and labels
    scale only axis=true,
    % set `width' and `height' to the desired values
    width=0.8\columnwidth,
    height=0.4\columnwidth,
    legend style={at={(1.33,1)},anchor=north,legend cell align=left},
    legend columns=2,
    %transpose legend
]

\addplot[color=c1,ultra thick,]
coordinates {(0,0)(1,0.417681)(2,0.186838)(3,0.246878)(4,0.216218)(5,0.2129)(6,0.207656)(7,0.202096)(8,0.201067)(9,0.197458)(10,0.196858)(11,0.195316)(12,0.194629)(13,0.194052)(14,0.193578)(15,0.193304)};
\addplot[color=c2,ultra thick,]
coordinates {(0,0)(1,0.191013)(2,0.258165)(3,0.141945)(4,0.220402)(5,0.184355)(6,0.194052)(7,0.194236)(8,0.191489)(9,0.193509)(10,0.192244)(11,0.192743)(12,0.192509)(13,0.192524)(14,0.192521)(15,0.19248)};
\addplot[color=c3,ultra thick,]
coordinates {(0,0)(1,0.146975)(2,0.140861)(3,0.147791)(4,0.163504)(5,0.174613)(6,0.174411)(7,0.182888)(8,0.183953)(9,0.186548)(10,0.188352)(11,0.18923)(12,0.190262)(13,0.19079)(14,0.191266)(15,0.191588)};
\addplot[color=c4,ultra thick,]
coordinates {(0,0.684605)(1,0.)(2,0.417681)(3,0.186838)(4,0.246878)(5,0.216218)(6,0.2129)(7,0.207656)(8,0.202096)(9,0.201067)(10,0.197458)(11,0.196858)(12,0.195316)(13,0.194629)(14,0.194052)(15,0.193578)};
\addplot[color=c5,ultra thick,]
coordinates {(0,0.596201)(1,0.)(2,0.417681)(3,0.186838)(4,0.246878)(5,0.216218)(6,0.2129)(7,0.207656)(8,0.202096)(9,0.201067)(10,0.197458)(11,0.196858)(12,0.195316)(13,0.194629)(14,0.194052)(15,0.193578)};
\addplot[color=c6,ultra thick,]
coordinates {(0,0.719194)(1,0.596201)(2,0.)(3,0.417681)(4,0.186838)(5,0.246878)(6,0.216218)(7,0.2129)(8,0.207656)(9,0.202096)(10,0.201067)(11,0.197458)(12,0.196858)(13,0.195316)(14,0.194629)(15,0.194052)};
\addplot[color=c7,ultra thick,]
coordinates {(0,0.347897)(1,0.)(2,0.191013)(3,0.258165)(4,0.141945)(5,0.220402)(6,0.184355)(7,0.194052)(8,0.194236)(9,0.191489)(10,0.193509)(11,0.192244)(12,0.192743)(13,0.192509)(14,0.192524)(15,0.192521)};
\addplot[color=c8,ultra thick,]
coordinates {(0,0.726167)(1,0.)(2,0.191013)(3,0.258165)(4,0.141945)(5,0.220402)(6,0.184355)(7,0.194052)(8,0.194236)(9,0.191489)(10,0.193509)(11,0.192244)(12,0.192743)(13,0.192509)(14,0.192524)(15,0.192521)};
\addplot[color=c9,ultra thick,]
coordinates {(0,0.25927)(1,0.726167)(2,0.)(3,0.191013)(4,0.258165)(5,0.141945)(6,0.220402)(7,0.184355)(8,0.194052)(9,0.194236)(10,0.191489)(11,0.193509)(12,0.192244)(13,0.192743)(14,0.192509)(15,0.192524)};
\addplot[color=c10,ultra thick,]
coordinates {(0,0.0304891)(1,0.)(2,0.146975)(3,0.140861)(4,0.147791)(5,0.163504)(6,0.174613)(7,0.174411)(8,0.182888)(9,0.183953)(10,0.186548)(11,0.188352)(12,0.18923)(13,0.190262)(14,0.19079)(15,0.191266)};
\addplot[color=c11,ultra thick,]
coordinates {(0,0.0956126)(1,0.)(2,0.146975)(3,0.140861)(4,0.147791)(5,0.163504)(6,0.174613)(7,0.174411)(8,0.182888)(9,0.183953)(10,0.186548)(11,0.188352)(12,0.18923)(13,0.190262)(14,0.19079)(15,0.191266)};
\addplot[color=c12,ultra thick,]
coordinates {(0,0.673898)(1,0.0956126)(2,0.)(3,0.146975)(4,0.140861)(5,0.147791)(6,0.163504)(7,0.174613)(8,0.174411)(9,0.182888)(10,0.183953)(11,0.186548)(12,0.188352)(13,0.18923)(14,0.190262)(15,0.19079)};

\addplot[color=black,ultra thick,]
coordinates {(0,0.344444)(15,0.344444)};

\legend{{ $x_{1}[k]$}, { $x_{2}[k]$}, { $x_{3}[k]$}, { $x_{4}[k]$}, { $x_{5}[k]$}, { $x_{6}[k]$}, { $x_{7}[k]$}, { $x_{8}[k]$}, { $x_{9}[k]$}, { $x_{10}[k]$}, { $x_{11}[k]$}, { $x_{12}[k]$}}

\end{axis}
\end{tikzpicture}
}
\qquad
\subfigure[{State evolution of agents from network $\mathcal G(A^\text{\textsf{P}})$ with weight matrix $A^\text{\textsf{P}}$ and re-scaled initial states $\tilde{x}^{\text{\textsf{P}}}[0]$.}]{
\tikzset{
  font={\fontsize{12pt}{10}\selectfont}}

\begin{tikzpicture}[scale=.66]
\begin{axis}[
    title={Consensus evolution},
    xlabel={Steps},
    ylabel={Agents evolution},
    ylabel shift=-4pt,
    xlabel shift=-1pt,
    xmin=0, xmax=15,
    ymin=0., ymax=1.8, %-3.15, 5.4
    ytick={0,0.344444,1,1.5},
    legend pos=north west,
    ymajorgrids=true,
    xmajorgrids=true,
    grid style=dashed,
    % only scale the axis, not the axis including the ticks and labels
    scale only axis=true,
    % set `width' and `height' to the desired values
    width=0.8\columnwidth,
    height=0.4\columnwidth,
    legend style={at={(1.33,1)},anchor=north,legend cell align=left},
    legend columns=2,
    %transpose legend
]

\addplot[color=c1,ultra thick,]
coordinates {(0,0)(1,0.63171)(2,0.38633)(3,0.400789)(4,0.388617)(5,0.374043)(6,0.366705)(7,0.359618)(8,0.357015)(9,0.352193)(10,0.350957)(11,0.348808)(12,0.347695)(13,0.34685)(14,0.346137)(15,0.345716)};
\addplot[color=c2,ultra thick,]
coordinates {(0,0)(1,0.273225)(2,0.516512)(3,0.231358)(4,0.402879)(5,0.329455)(6,0.347441)(7,0.348487)(8,0.342618)(9,0.346768)(10,0.344057)(11,0.345138)(12,0.344584)(13,0.344625)(14,0.344604)(15,0.34451)};
\addplot[color=c3,ultra thick,]
coordinates {(0,0)(1,0.316974)(2,0.224013)(3,0.307358)(4,0.279429)(5,0.326045)(6,0.313943)(7,0.329133)(8,0.332407)(9,0.334368)(10,0.338582)(11,0.339231)(12,0.341095)(13,0.341865)(14,0.342579)(15,0.343113)};
\addplot[color=c4,ultra thick,]
coordinates {(0,0.770181)(1,0.)(2,0.63171)(3,0.38633)(4,0.400789)(5,0.388617)(6,0.374043)(7,0.366705)(8,0.359618)(9,0.357015)(10,0.352193)(11,0.350957)(12,0.348808)(13,0.347695)(14,0.34685)(15,0.346137)};
\addplot[color=c5,ultra thick,]
coordinates {(0,1.34145)(1,0.)(2,0.63171)(3,0.38633)(4,0.400789)(5,0.388617)(6,0.374043)(7,0.366705)(8,0.359618)(9,0.357015)(10,0.352193)(11,0.350957)(12,0.348808)(13,0.347695)(14,0.34685)(15,0.346137)};
\addplot[color=c6,ultra thick,]
coordinates {(0,1.61819)(1,1.34145)(2,0.)(3,0.63171)(4,0.38633)(5,0.400789)(6,0.388617)(7,0.374043)(8,0.366705)(9,0.359618)(10,0.357015)(11,0.352193)(12,0.350957)(13,0.348808)(14,0.347695)(15,0.34685)};
\addplot[color=c7,ultra thick,]
coordinates {(0,0.391384)(1,0.)(2,0.273225)(3,0.516512)(4,0.231358)(5,0.402879)(6,0.329455)(7,0.347441)(8,0.348487)(9,0.342618)(10,0.346768)(11,0.344057)(12,0.345138)(13,0.344584)(14,0.344625)(15,0.344604)};
\addplot[color=c8,ultra thick,]
coordinates {(0,1.63388)(1,0.)(2,0.273225)(3,0.516512)(4,0.231358)(5,0.402879)(6,0.329455)(7,0.347441)(8,0.348487)(9,0.342618)(10,0.346768)(11,0.344057)(12,0.345138)(13,0.344584)(14,0.344625)(15,0.344604)};
\addplot[color=c9,ultra thick,]
coordinates {(0,0.583356)(1,1.63388)(2,0.)(3,0.273225)(4,0.516512)(5,0.231358)(6,0.402879)(7,0.329455)(8,0.347441)(9,0.348487)(10,0.342618)(11,0.346768)(12,0.344057)(13,0.345138)(14,0.344584)(15,0.344625)};
\addplot[color=c10,ultra thick,]
coordinates {(0,0.0343002)(1,0.)(2,0.316974)(3,0.224013)(4,0.307358)(5,0.279429)(6,0.326045)(7,0.313943)(8,0.329133)(9,0.332407)(10,0.334368)(11,0.338582)(12,0.339231)(13,0.341095)(14,0.341865)(15,0.342579)};
\addplot[color=c11,ultra thick,]
coordinates {(0,0.215128)(1,0.)(2,0.316974)(3,0.224013)(4,0.307358)(5,0.279429)(6,0.326045)(7,0.313943)(8,0.329133)(9,0.332407)(10,0.334368)(11,0.338582)(12,0.339231)(13,0.341095)(14,0.341865)(15,0.342579)};
\addplot[color=c12,ultra thick,]
coordinates {(0,1.51627)(1,0.215128)(2,0.)(3,0.316974)(4,0.224013)(5,0.307358)(6,0.279429)(7,0.326045)(8,0.313943)(9,0.329133)(10,0.332407)(11,0.334368)(12,0.338582)(13,0.339231)(14,0.341095)(15,0.341865)};

\addplot[color=black,ultra thick,]
coordinates {(0,0.341865)(15,0.341865)};

\legend{{ $x_{1}[k]$}, { $x_{2}[k]$}, { $x_{3}[k]$}, { $x_{4}[k]$}, { $x_{5}[k]$}, { $x_{6}[k]$}, { $x_{7}[k]$}, { $x_{8}[k]$}, { $x_{9}[k]$}, { $x_{10}[k]$}, { $x_{11}[k]$}, { $x_{12}[k]$}}

\end{axis}
\end{tikzpicture}
}
\caption{Privacy for bidirected cycle of size 3.}
\label{fig:cycle3_consensus}
\end{figure}

In Figure~\ref{fig:cycle3_consensus}~(a), we depict the agents state evolution under network $\mathcal G(A)$, weight matrix $A$, and initial agents' states $x[0]=\left[\begin{smallmatrix}
   \frac{1}{2} & \frac{1}{3} & \frac{1}{5}
\end{smallmatrix}\right]^\intercal$.

In Figure~\ref{fig:cycle3_consensus}~(b), we depict the agents state evolution under network $\mathcal G(A^{\text{\textsf{P}}})$, weight matrix $A^{\text{\textsf{P}}}$, and initial agents' states \[\begin{split}
x'[0]=\left[\begin{smallmatrix}
  0 & 0 & 0 & 0.684605 & 0.596201 & 0.719194 & 0.347897 & 0.726167 & 0.25927  \end{smallmatrix}\right.
  \\
  \left.\begin{smallmatrix}0.0304891 & 0.0956126 & 0.673898 
\end{smallmatrix}\right]^\intercal.\end{split}
\]

In Figure~\ref{fig:cycle3_consensus}~(c), we depict the agents state evolution under network $\mathcal G(A^{\text{\textsf{P}}})$, weight matrix $A^{\text{\textsf{P}}}$, and re-scaled initial agents' states 
\[\begin{split}\tilde x^{\text{\textsf{P}}}[0]=\left[\begin{smallmatrix}
 0 & 0 & 0 & 0.770181 & 1.34145 & 1.61819 & 0.391384 & 1.63388 & 0.583356\end{smallmatrix}\right.\\
 \left.\begin{smallmatrix} 0.0343002 & 0.215128 & 1.51627 
\end{smallmatrix}\right]^\intercal.\end{split}
\]

The previous example illustrates the following result. 

\vspace{1mm} 
{\centering
\allowdisplaybreaks
\noindent \colorbox{c2!18}{\parbox{0.98\textwidth}{
\begin{lemma}\label{prop:3cpy_obs}
     Let $A\in\mathbb R^{N\times N}$ be a matrix such that $\mathcal G(\bar A)$ is a strongly connected graph. 
     Then, the system $\tilde{x}^{\text{\textsf{P}}}[k+1]=A^\text{\textsf{P}} \tilde{x}^{\text{\textsf{P}}}[k]$ and $y[k]=C\tilde{x}^{\text{\textsf{P}}}[k]$ , where the matrix $A^\text{\textsf{P}}\in\mathbb R^{4N\times 4N}$ results from Algorithm~\ref{alg:buildAP2} and 
\[
C =
\begin{bmatrix}
\bigI_{N+3} & \rvline  & \bigzero_{(N+3) \times (3N-3)} \\
\end{bmatrix}
\]
has the following properties: 
     \begin{enumerate}
         \item it is not observable;
         \item $\left( e_{5N}^j+\displaystyle\sum_{k=N+4j-3}^{N+4j} e_{5N}^{k}\right)\notin span\left( P_{O(A^\text{\textsf{P}},C)}^\lambda\right)$, for all $j\neq i$.\hfill$\circ$  
        %  \item $\displaystyle\left(\frac{1}{{(v_L)}_j} e_{5N}^j+\displaystyle\sum_{k=N+4j-3}^{N+4j} \frac{1}{{(v_L)}_k} e_{5N}^{k}\right)\notin span\left( P_{O(A^\text{\textsf{P}},C)}^\lambda\right)$, for all $j\neq i$,
     \end{enumerate} 
    % where $v_L$ is the left-eigenvector of $A^\text{\textsf{P}}$ associated with the eigenvalue $1$. 
 \end{lemma}  
 }
 }
 }
 \vspace{1mm}
 
\begin{proof}
    We construct the observabillity matrix $P_{O(A^\text{\textsf{P}},C)}^0$ and compute its span. 
    \[
    \begin{array}{l}
        span\left(P_{O(A^\text{\textsf{P}},C)}^0\right)= \\
        span\left(\begin{bmatrix}
        \bigI_{N+3} & \rvline  & \bigzero_{(N+3) \times (3N-3)} \\ \hline
        \bigzero_{(3N-3)\times (N+3)} & \rvline & \bigR
        \end{bmatrix}
        \right),
        \end{array}
    \]
where $R\in\mathbb R^{(3N-3)\times (3N-3)}$ is $R=diag(D,\ldots,D)$ with
\[
D=
\begin{bmatrix}
 1 & 0 & \frac{1}{2} \\
 0 & 1 & 0 \\
 0 & 0 & 0
\end{bmatrix}
\]
Hence, $P_{O(A^\text{\textsf{P}},C)}^0$ is not full rank and $(1)$ holds.  Furthermore, $e_{4N}^{N+3j-2}+e_{4N}^{N+4j-1}+e_{4N}^{N+3j-0}+e_{4N}^{j}\notin span\left( P_{O(A^\text{\textsf{P}},C)}^\lambda\right)$ for all $j\neq i$, which implies that $(2)$ holds. 
\end{proof}

\vspace{1mm} 
{\centering
\allowdisplaybreaks
\noindent \colorbox{c3!15}{\parbox{0.98\textwidth}{
\begin{theorem}\label{cor:3cpy}
    Under the conditions of Lemma~\ref{prop:3cpy_obs}, if $1)$ and $2)$ hold and 
     $\displaystyle\left(\frac{1}{{(v_L)}_j} e_{4N}^j+\displaystyle\sum_{k=N+3j-2}^{N+3j} \frac{1}{{(v_L)}_k} e_{4N}^{k}\right)\notin span\left( P_{O(A^\text{\textsf{P}},C)}^\lambda\right)$, for all $j\neq i$,  
    where $v_L$ is the left-eigenvector of $A^\text{\textsf{P}}$ obtained with Algorithm~\ref{alg:buildAP2}, associated with the eigenvalue $1$, then the initial state of agent $j$ that is distributed among its augmented states nodes is private.~\hfill$\circ$
\end{theorem}
}
}
}
\vspace{1mm}

Observe that a disadvantage of the solution is that we need to compute the left-eigenvector associated with the eigenvalue 1. 
In general, we can only compute an approximation of this eigenvector in polynomial-time. 
Further, notice that a similar path can be followed using any other of the augmented states networks depicted in Figure~\ref{fig:3copies_that_work}. 

In the quest of eliminating the aforementioned drawback, we explore the situation where the considered network is bidirected but with possible distinct weights, modeled by a time-reversible transition matrix. 
This scenario will allow us to have an expedient and efficient way of computing the left-eigenvector of the dynamics matrix.  

\subsection{Solution to Problem $\mathbf{P_1^D}$}\label{sub:rev_mat} 

Let $A\in\mathbb R^{N\times N}$ be the system dynamics matrix, the matrix $A$ is said to be \textit{irreversible} is there is $s\in\mathbb R$ such that $\sum_{i=1}^N s_i=1$ and 
\[
s_i A_{ij} = s_j A_{ji}
\]
for all $i,j=1,\ldots,N$. 

Notice that, in this scenario, the network of agents has to be bidirected. 
The first result is that we cannot use solely three augmented states of each agent to both address $\mathbf{P_1}$ and ensure a time-reversible transition matrix.

\vspace{1mm} 
{\centering
\allowdisplaybreaks
\noindent \colorbox{c2!18}{\parbox{0.98\textwidth}{
\begin{lemma}\label{prop:rev_3_cpy_no}
     Given a bidirected network of agents, we cannot use three augmented states of each agent connected by a bidirected network to address $\mathbf{P_1}$.\hfill$\circ$ 
\end{lemma}
}
}
}
\vspace{1mm}

\begin{proof} 
The previous result trivially follows from the fact that to ensure reversibility the network should be bidirected and there is not a bidirected network in Figure~\ref{fig:3copies_that_work}, i.e., the networks with three augmented states of an agent that allows to address $\mathbf{P_1}$. 
\end{proof}

Therefore, we need to consider four augmented states of each agent, as detailed in the following.

\vspace{1mm} 
{\centering
\allowdisplaybreaks
\noindent \colorbox{c2!18}{\parbox{0.98\textwidth}{
\begin{lemma}\label{prop:enumerate_3cpy}
The only networks (up to isomorphism) with four augmented states that satisfy the necessary conditions to address $\mathbf{P_1^D}$ are the ones depicted in Figure~\ref{fig:4copies_that_work}. 
\begin{figure}[H]
\centering
\subfigure[]{
\begin{tikzpicture}[scale=.27, transform shape,node distance=1.5cm]
\begin{scope}[every node/.style={circle,thick,draw},square/.style={regular polygon,regular polygon sides=4}]
\node[fill=green!20] (1) at (2.15841,0.) {\Huge $1$};
\node[draw=gray,fill=yellow!50] (2) at (2.156,3.01787) {\Huge $2$};
\node[draw=gray,fill=yellow!50] (3) at (2.60552,1.50681) {\Huge $3$};
\node[draw=gray,fill=yellow!50] (4) at (4.67655,1.51032) {\Huge $4$};
\node[draw=gray,fill=yellow!50] (5) at (0.,1.50851) {\Huge $5$};
\end{scope}
\begin{scope}[>={Stealth[black]},
              every edge/.style={draw=gray, thick}]
\path [-] (1) edge node {} (2);
\path [-] (1) edge node {} (3);
\path [-] (1) edge node {} (4);
\path [-] (1) edge node {} (5);
\path [-] (2) edge node {} (3);
\path [-] (2) edge node {} (4);
\path [-] (2) edge node {} (5);
\path [-] (3) edge node {} (4);
\path [-] (3) edge node {} (5);
\end{scope}
\end{tikzpicture}
 }
% \subfigure[]{
% \begin{tikzpicture}[scale=.27, transform shape,node distance=1.5cm]
% \begin{scope}[every node/.style={circle,thick,draw},square/.style={regular polygon,regular polygon sides=4}]
% \node[fill=green!20] (1) at (2.04124,4.08248) {\Huge $1$};
% \node[draw=gray,fill=yellow!50] (2) at (0.,2.04124) {\Huge $2$};
% \node[draw=gray,fill=yellow!50] (3) at (2.04124,2.04124) {\Huge $3$};
% \node[draw=gray,fill=yellow!50] (4) at (4.08248,2.04124) {\Huge $4$};
% \node[draw=gray,fill=yellow!50] (5) at (0.,0.) {\Huge $5$};
% \end{scope}
% \begin{scope}[>={Stealth[black]},
%               every edge/.style={draw=gray, thick}]
% \path [-] (1) edge node {} (2);
% \path [-] (1) edge node {} (3);
% \path [-] (1) edge node {} (4);
% \path [-] (2) edge node {} (5);
% \end{scope}
% \end{tikzpicture}
% }
% \subfigure[]{
% \begin{tikzpicture}[scale=.27, transform shape,node distance=1.5cm]
% \begin{scope}[every node/.style={circle,thick,draw},square/.style={regular polygon,regular polygon sides=4}]
% \node[fill=green!20] (1) at (5.,0.) {\Huge $1$};
% \node[draw=gray,fill=yellow!50] (2) at (2.5,0.) {\Huge $2$};
% \node[draw=gray,fill=yellow!50] (3) at (7.5,0.) {\Huge $3$};
% \node[draw=gray,fill=yellow!50] (4) at (0.,0.) {\Huge $4$};
% \node[draw=gray,fill=yellow!50] (5) at (10.,0.) {\Huge $5$};
% \end{scope}
% \begin{scope}[>={Stealth[black]},
%               every edge/.style={draw=gray, thick}]
% \path [-] (1) edge node {} (2);
% \path [-] (1) edge node {} (3);
% \path [-] (2) edge node {} (4);
% \path [-] (3) edge node {} (5);
% \end{scope}
% \end{tikzpicture}
% }
\subfigure[]{
\begin{tikzpicture}[scale=.27, transform shape,node distance=1.5cm]
\begin{scope}[every node/.style={circle,thick,draw},square/.style={regular polygon,regular polygon sides=4}]
\node[fill=green!20] (1) at (2.35885,2.00831) {\Huge $1$};
\node[draw=gray,fill=yellow!50] (2) at (0.00036023,2.94345) {\Huge $2$};
\node[draw=gray,fill=yellow!50] (3) at (0.,1.07348) {\Huge $3$};
\node[draw=gray,fill=yellow!50] (4) at (4.20429,0.) {\Huge $4$};
\node[draw=gray,fill=yellow!50] (5) at (4.20593,4.01517) {\Huge $5$};
\end{scope}
\begin{scope}[>={Stealth[black]},
              every edge/.style={draw=gray, thick}]
\path [-] (1) edge node {} (2);
\path [-] (1) edge node {} (3);
\path [-] (1) edge node {} (4);
\path [-] (1) edge node {} (5);
\path [-] (2) edge node {} (3);
\end{scope}
\end{tikzpicture}
}
\subfigure[]{
\begin{tikzpicture}[scale=.27, transform shape,node distance=1.5cm]
\begin{scope}[every node/.style={circle,thick,draw},square/.style={regular polygon,regular polygon sides=4}]
\node[fill=green!20] (1) at (5.30846,0.713606) {\Huge $1$};
\node[draw=gray,fill=yellow!50] (2) at (2.5603,0.711594) {\Huge $2$};
\node[draw=gray,fill=yellow!50] (3) at (3.93079,2.41316) {\Huge $3$};
\node[draw=gray,fill=yellow!50] (4) at (7.87142,0.00131917) {\Huge $4$};
\node[draw=gray,fill=yellow!50] (5) at (0.,0.) {\Huge $5$};
\end{scope}
\begin{scope}[>={Stealth[black]},
              every edge/.style={draw=gray, thick}]
\path [-] (1) edge node {} (2);
\path [-] (1) edge node {} (3);
\path [-] (1) edge node {} (4);
\path [-] (2) edge node {} (3);
\path [-] (2) edge node {} (5);
\end{scope}
\end{tikzpicture}
}
\subfigure[]{
\begin{tikzpicture}[scale=.27, transform shape,node distance=1.5cm]
\begin{scope}[every node/.style={circle,thick,draw},square/.style={regular polygon,regular polygon sides=4}]
\node[fill=green!20] (1) at (2.25767,0.947948) {\Huge $1$};
\node[draw=gray,fill=yellow!50] (2) at (0.,1.89521) {\Huge $2$};
\node[draw=gray,fill=yellow!50] (3) at (0.000373252,0.) {\Huge $3$};
\node[draw=gray,fill=yellow!50] (4) at (5.27042,0.947848) {\Huge $4$};
\node[draw=gray,fill=yellow!50] (5) at (7.83813,0.947628) {\Huge $5$};
\end{scope}
\begin{scope}[>={Stealth[black]},
              every edge/.style={draw=gray, thick}]
\path [-] (1) edge node {} (2);
\path [-] (1) edge node {} (3);
\path [-] (1) edge node {} (4);
\path [-] (2) edge node {} (3);
\path [-] (4) edge node {} (5);
\end{scope}
\end{tikzpicture}
}
\subfigure[]{
\begin{tikzpicture}[scale=.27, transform shape,node distance=1.5cm]
\begin{scope}[every node/.style={circle,thick,draw},square/.style={regular polygon,regular polygon sides=4}]
\node[fill=green!20] (1) at (4.00477,1.41476) {\Huge $1$};
\node[draw=gray,fill=yellow!50] (2) at (1.83636,2.82885) {\Huge $2$};
\node[draw=gray,fill=yellow!50] (3) at (1.83716,0.) {\Huge $3$};
\node[draw=gray,fill=yellow!50] (4) at (6.6693,1.41423) {\Huge $4$};
\node[draw=gray,fill=yellow!50] (5) at (0.,1.4139) {\Huge $5$};
\end{scope}
\begin{scope}[>={Stealth[black]},
              every edge/.style={draw=gray, thick}]
\path [-] (1) edge node {} (2);
\path [-] (1) edge node {} (3);
\path [-] (1) edge node {} (4);
\path [-] (2) edge node {} (5);
\path [-] (3) edge node {} (5);
\end{scope}
\end{tikzpicture}
}
% \subfigure[]{
% \begin{tikzpicture}[scale=.27, transform shape,node distance=1.5cm]
% \begin{scope}[every node/.style={circle,thick,draw},square/.style={regular polygon,regular polygon sides=4}]
% \node[fill=green!20] (1) at (-2.37764,0.772542) {\Huge $1$};
% \node[draw=gray,fill=yellow!50] (2) at (0.,2.5) {\Huge $2$};
% \node[draw=gray,fill=yellow!50] (3) at (-1.46946,-2.02254) {\Huge $3$};
% \node[draw=gray,fill=yellow!50] (4) at (2.37764,0.772542) {\Huge $4$};
% \node[draw=gray,fill=yellow!50] (5) at (1.46946,-2.02254) {\Huge $5$};
% \end{scope}
% \begin{scope}[>={Stealth[black]},
%               every edge/.style={draw=gray, thick}]
% \path [-] (1) edge node {} (2);
% \path [-] (1) edge node {} (3);
% \path [-] (2) edge node {} (4);
% \path [-] (3) edge node {} (5);
% \path [-] (4) edge node {} (5);
% \end{scope}
% \end{tikzpicture}
% }
\subfigure[]{
\begin{tikzpicture}[scale=.27, transform shape,node distance=1.5cm]
\begin{scope}[every node/.style={circle,thick,draw},square/.style={regular polygon,regular polygon sides=4}]
\node[fill=green!20] (1) at (2.34089,2.07973) {\Huge $1$};
\node[draw=gray,fill=yellow!50] (2) at (0.,2.0795) {\Huge $2$};
\node[draw=gray,fill=yellow!50] (3) at (0.753276,4.16075) {\Huge $3$};
\node[draw=gray,fill=yellow!50] (4) at (0.754475,0.) {\Huge $4$};
\node[draw=gray,fill=yellow!50] (5) at (5.26227,2.08073) {\Huge $5$};
\end{scope}
\begin{scope}[>={Stealth[black]},
              every edge/.style={draw=gray, thick}]
\path [-] (1) edge node {} (2);
\path [-] (1) edge node {} (3);
\path [-] (1) edge node {} (4);
\path [-] (1) edge node {} (5);
\path [-] (2) edge node {} (3);
\path [-] (2) edge node {} (4);
\end{scope}
\end{tikzpicture}
}
\subfigure[]{
\begin{tikzpicture}[scale=.27, transform shape,node distance=1.5cm]
\begin{scope}[every node/.style={circle,thick,draw},square/.style={regular polygon,regular polygon sides=4}]
\node[fill=green!20] (1) at (2.5214,1.0049) {\Huge $1$};
\node[draw=gray,fill=yellow!50] (2) at (5.04137,0.) {\Huge $2$};
\node[draw=gray,fill=yellow!50] (3) at (5.0435,2.00777) {\Huge $3$};
\node[draw=gray,fill=yellow!50] (4) at (0.,2.00701) {\Huge $4$};
\node[draw=gray,fill=yellow!50] (5) at (0.00201513,0.00081012) {\Huge $5$};
\end{scope}
\begin{scope}[>={Stealth[black]},
              every edge/.style={draw=gray, thick}]
\path [-] (1) edge node {} (2);
\path [-] (1) edge node {} (3);
\path [-] (1) edge node {} (4);
\path [-] (1) edge node {} (5);
\path [-] (2) edge node {} (3);
\path [-] (4) edge node {} (5);
\end{scope}
\end{tikzpicture}
}
\subfigure[]{
\begin{tikzpicture}[scale=.27, transform shape,node distance=1.5cm]
\begin{scope}[every node/.style={circle,thick,draw},square/.style={regular polygon,regular polygon sides=4}]
\node[fill=green!20] (1) at (5.3556,1.86716) {\Huge $1$};
\node[draw=gray,fill=yellow!50] (2) at (5.35509,0.) {\Huge $2$};
\node[draw=gray,fill=yellow!50] (3) at (2.84573,0.933673) {\Huge $3$};
\node[draw=gray,fill=yellow!50] (4) at (7.57265,0.933233) {\Huge $4$};
\node[draw=gray,fill=yellow!50] (5) at (0.,0.933388) {\Huge $5$};
\end{scope}
\begin{scope}[>={Stealth[black]},
              every edge/.style={draw=gray, thick}]
\path [-] (1) edge node {} (2);
\path [-] (1) edge node {} (3);
\path [-] (1) edge node {} (4);
\path [-] (2) edge node {} (3);
\path [-] (2) edge node {} (4);
\path [-] (3) edge node {} (5);
\end{scope}
\end{tikzpicture}
}
\subfigure[]{
\begin{tikzpicture}[scale=.27, transform shape,node distance=1.5cm]
\begin{scope}[every node/.style={circle,thick,draw},square/.style={regular polygon,regular polygon sides=4}]
\node[fill=green!20] (1) at (2.16129,2.34506) {\Huge $1$};
\node[draw=gray,fill=yellow!50] (2) at (2.16136,0.0904008) {\Huge $2$};
\node[draw=gray,fill=yellow!50] (3) at (0.,1.21772) {\Huge $3$};
\node[draw=gray,fill=yellow!50] (4) at (4.86183,2.4354) {\Huge $4$};
\node[draw=gray,fill=yellow!50] (5) at (4.86168,0.) {\Huge $5$};
\end{scope}
\begin{scope}[>={Stealth[black]},
              every edge/.style={draw=gray, thick}]
\path [-] (1) edge node {} (2);
\path [-] (1) edge node {} (3);
\path [-] (1) edge node {} (4);
\path [-] (2) edge node {} (3);
\path [-] (2) edge node {} (5);
\path [-] (4) edge node {} (5);
\end{scope}
\end{tikzpicture}
}
\subfigure[]{
\begin{tikzpicture}[scale=.27, transform shape,node distance=1.5cm]
\begin{scope}[every node/.style={circle,thick,draw},square/.style={regular polygon,regular polygon sides=4}]
\node[fill=green!20] (1) at (2.51535,1.12158) {\Huge $1$};
\node[draw=gray,fill=yellow!50] (2) at (2.5193,2.55272) {\Huge $2$};
\node[draw=gray,fill=yellow!50] (3) at (4.28015,0.) {\Huge $3$};
\node[draw=gray,fill=yellow!50] (4) at (4.28379,3.67421) {\Huge $4$};
\node[draw=gray,fill=yellow!50] (5) at (0.,1.8406) {\Huge $5$};
\end{scope}
\begin{scope}[>={Stealth[black]},
              every edge/.style={draw=gray, thick}]
\path [-] (1) edge node {} (2);
\path [-] (1) edge node {} (3);
\path [-] (1) edge node {} (4);
\path [-] (1) edge node {} (5);
\path [-] (2) edge node {} (3);
\path [-] (2) edge node {} (4);
\path [-] (2) edge node {} (5);
\end{scope}
\end{tikzpicture}
}
\subfigure[]{
\begin{tikzpicture}[scale=.27, transform shape,node distance=1.5cm]
\begin{scope}[every node/.style={circle,thick,draw},square/.style={regular polygon,regular polygon sides=4}]
\node[fill=green!20] (1) at (3.14545,1.32878) {\Huge $1$};
\node[draw=gray,fill=yellow!50] (2) at (5.00361,2.65536) {\Huge $2$};
\node[draw=gray,fill=yellow!50] (3) at (6.25937,1.32931) {\Huge $3$};
\node[draw=gray,fill=yellow!50] (4) at (5.00228,0.) {\Huge $4$};
\node[draw=gray,fill=yellow!50] (5) at (0.,1.32893) {\Huge $5$};
\end{scope}
\begin{scope}[>={Stealth[black]},
              every edge/.style={draw=gray, thick}]
\path [-] (1) edge node {} (2);
\path [-] (1) edge node {} (3);
\path [-] (1) edge node {} (4);
\path [-] (1) edge node {} (5);
\path [-] (2) edge node {} (3);
\path [-] (2) edge node {} (4);
\path [-] (3) edge node {} (4);
\end{scope}
\end{tikzpicture}
}
\subfigure[]{
\begin{tikzpicture}[scale=.27, transform shape,node distance=1.5cm]
\begin{scope}[every node/.style={circle,thick,draw},square/.style={regular polygon,regular polygon sides=4}]
\node[fill=green!20] (1) at (2.74682,1.8229) {\Huge $1$};
\node[draw=gray,fill=yellow!50] (2) at (1.42275,0.00185561) {\Huge $2$};
\node[draw=gray,fill=yellow!50] (3) at (4.06809,0.) {\Huge $3$};
\node[draw=gray,fill=yellow!50] (4) at (0.,1.92765) {\Huge $4$};
\node[draw=gray,fill=yellow!50] (5) at (5.49673,1.9265) {\Huge $5$};
\end{scope}
\begin{scope}[>={Stealth[black]},
              every edge/.style={draw=gray, thick}]
\path [-] (1) edge node {} (2);
\path [-] (1) edge node {} (3);
\path [-] (1) edge node {} (4);
\path [-] (1) edge node {} (5);
\path [-] (2) edge node {} (3);
\path [-] (2) edge node {} (4);
\path [-] (3) edge node {} (5);
\end{scope}
\end{tikzpicture}
}
\subfigure[]{
\begin{tikzpicture}[scale=.27, transform shape,node distance=1.5cm]
\begin{scope}[every node/.style={circle,thick,draw},square/.style={regular polygon,regular polygon sides=4}]
\node[fill=green!20] (1) at (0.00188258,2.11569) {\Huge $1$};
\node[draw=gray,fill=yellow!50] (2) at (0.,0.332117) {\Huge $2$};
\node[draw=gray,fill=yellow!50] (3) at (2.07868,0.) {\Huge $3$};
\node[draw=gray,fill=yellow!50] (4) at (2.07772,2.44686) {\Huge $4$};
\node[draw=gray,fill=yellow!50] (5) at (4.43902,1.22385) {\Huge $5$};
\end{scope}
\begin{scope}[>={Stealth[black]},
              every edge/.style={draw=gray, thick}]
\path [-] (1) edge node {} (2);
\path [-] (1) edge node {} (3);
\path [-] (1) edge node {} (4);
\path [-] (2) edge node {} (3);
\path [-] (2) edge node {} (4);
\path [-] (3) edge node {} (5);
\path [-] (4) edge node {} (5);
\end{scope}
\end{tikzpicture}
}
\subfigure[]{
\begin{tikzpicture}[scale=.27, transform shape,node distance=1.5cm]
\begin{scope}[every node/.style={circle,thick,draw},square/.style={regular polygon,regular polygon sides=4}]
\node[fill=green!20] (1) at (2.24531,1.92738) {\Huge $1$};
\node[draw=gray,fill=yellow!50] (2) at (2.24653,0.103459) {\Huge $2$};
\node[draw=gray,fill=yellow!50] (3) at (0.00199298,2.03206) {\Huge $3$};
\node[draw=gray,fill=yellow!50] (4) at (0.,0.) {\Huge $4$};
\node[draw=gray,fill=yellow!50] (5) at (4.82138,1.01589) {\Huge $5$};
\end{scope}
\begin{scope}[>={Stealth[black]},
              every edge/.style={draw=gray, thick}]
\path [-] (1) edge node {} (2);
\path [-] (1) edge node {} (3);
\path [-] (1) edge node {} (4);
\path [-] (1) edge node {} (5);
\path [-] (2) edge node {} (3);
\path [-] (2) edge node {} (4);
\path [-] (2) edge node {} (5);
\path [-] (3) edge node {} (4);
\end{scope}
\end{tikzpicture}
}
\subfigure[]{
\begin{tikzpicture}[scale=.27, transform shape,node distance=1.5cm]
\begin{scope}[every node/.style={circle,thick,draw},square/.style={regular polygon,regular polygon sides=4}]
\node[fill=green!20] (1) at (1.48609,1.48877) {\Huge $1$};
\node[draw=gray,fill=yellow!50] (2) at (0.,2.88431) {\Huge $2$};
\node[draw=gray,fill=yellow!50] (3) at (2.88441,2.97336) {\Huge $3$};
\node[draw=gray,fill=yellow!50] (4) at (0.0922107,0.) {\Huge $4$};
\node[draw=gray,fill=yellow!50] (5) at (2.97861,0.0871823) {\Huge $5$};
\end{scope}
\begin{scope}[>={Stealth[black]},
              every edge/.style={draw=gray, thick}]
\path [-] (1) edge node {} (2);
\path [-] (1) edge node {} (3);
\path [-] (1) edge node {} (4);
\path [-] (1) edge node {} (5);
\path [-] (2) edge node {} (3);
\path [-] (2) edge node {} (4);
\path [-] (3) edge node {} (5);
\path [-] (4) edge node {} (5);
\end{scope}
\end{tikzpicture}
}
\subfigure[]{
\begin{tikzpicture}[scale=.27, transform shape,node distance=1.5cm]
\begin{scope}[every node/.style={circle,thick,draw},square/.style={regular polygon,regular polygon sides=4}]
\node[fill=green!20] (1) at (-2.37764,0.772542) {\Huge $1$};
\node[draw=gray,fill=yellow!50] (2) at (0.,2.5) {\Huge $2$};
\node[draw=gray,fill=yellow!50] (3) at (2.37764,0.772542) {\Huge $3$};
\node[draw=gray,fill=yellow!50] (4) at (1.46946,-2.02254) {\Huge $4$};
\node[draw=gray,fill=yellow!50] (5) at (-1.46946,-2.02254) {\Huge $5$};
\end{scope}
\begin{scope}[>={Stealth[black]},
              every edge/.style={draw=gray, thick}]
\path [-] (1) edge node {} (2);
\path [-] (1) edge node {} (3);
\path [-] (1) edge node {} (4);
\path [-] (1) edge node {} (5);
\path [-] (2) edge node {} (3);
\path [-] (2) edge node {} (4);
\path [-] (2) edge node {} (5);
\path [-] (3) edge node {} (4);
\path [-] (3) edge node {} (5);
\path [-] (4) edge node {} (5);
\end{scope}
\end{tikzpicture}
}
\caption{Bidirected networks with one agent (node 1) and four augmented states that can be used to address problem $\mathbf{P_1}$ and ensure time reversibility.\hfill$\circ$}
\label{fig:4copies_that_work}
\end{figure}
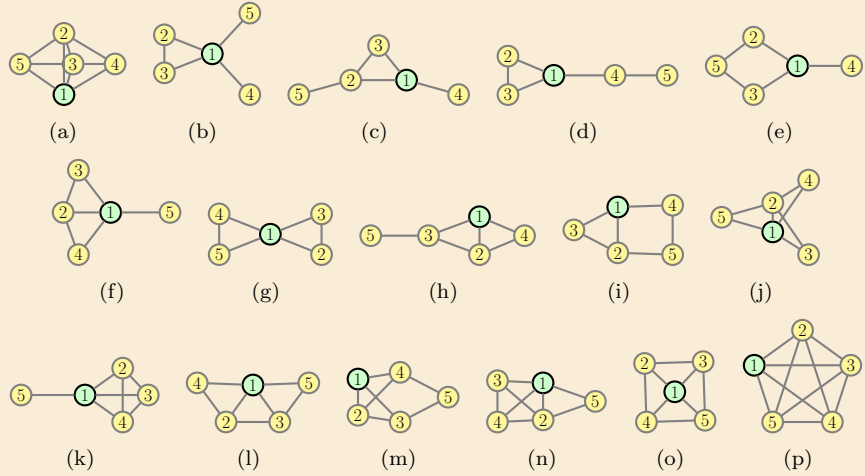
\end{lemma}
}
}
}
\vspace{1mm}

Next, we propose an algorithm to address problem $\mathbf{P_{1}^D}$ that produces a time-reversible matrix $A^\text{\textsf{P}}$ corresponding to a bidirected graph whenever the input matrix $A$ is row stochastic, time-reversible and corresponds to a bidirected graph. To this end we use the network of agent's augmented states depicted in Figure~\ref{fig:4copies_that_work}~(a). 
Notice that, in contrast with the previous solution, we require that the input is the matrix $A$, row stochastic, time-reversible and corresponding to a bidirected graph. 

\begin{algorithm}[!ht]
		{%\small
		\caption{Construction of the privacy dynamics matrix $A^\text{\textsf{P}}$ for consensus of bidirectional networks}
		\label{alg:buildAP_rev}
		\begin{algorithmic}[1]
			\STATE{\textbf{input}: row stochastic dynamics matrix $A\in\mathbb R^{N\times N}$, time-reversible and that corresponds to a bidirected graph}
			\STATE{\textbf{output}: dynamics matrix $A^\text{\textsf{P}}\in\mathbb R^{5N\times 5N}$}
		\STATE{\textbf{fill} the entries of $A^{\text{\textsf{P}}}$ with zeros}
		\STATE{\textbf{set}
		\hfill{\color{gray}$\rhd$
		 \underline{Copy the matrix $A$ to the first $n$ rows and $n$ columns}}
		$$A^{\text{\textsf{P}}}_{ij}=\frac{A_{ij}}{2}, \forall i,j\in\{1,\ldots,N\}\text{ and }i\neq j$$}
		\STATE{\textbf{for} $i=1,\ldots,N$\hfill{\color{gray}$\rhd$ \underline{Set additional entries values}}
		$$
		\begin{array}{ll}
		     A^{\text{\textsf{P}}}_{i, N + 4 i - 3}  = \frac{1}{12}, &  A^{\text{\textsf{P}}}_{i, N + 4 i - 2}  =  \frac{1}{8}, \\
	    	A^{\text{\textsf{P}}}_{i, N + 4 i - 1} = \frac{1}{4}, & A^{\text{\textsf{P}}}_{i, N + 4 i - 0} = \frac{1}{24}, 
		    \\
		    A^{\text{\textsf{P}}}_{N + 4 i - 2, N + 4 i - 3} = \frac{1}{2}, & A^{\text{\textsf{P}}}_{N + 4 i - 3, N + 4 i - 2} = \frac{3}{22}, 
            \\
		    A^{\text{\textsf{P}}}_{N + 4 i - 1, N + 4 i - 3} = \frac{1}{4}, & 
		    A^{\text{\textsf{P}}}_{N + 4 i - 3, N + 4 i - 1} = \frac{1}{11}, \\
		    A^{\text{\textsf{P}}}_{N + 4 i - 0, N + 4 i - 3} = \frac{15}{16}, &
		    A^{\text{\textsf{P}}}_{N + 4 i - 3, N + 4 i - 0} = \frac{15}{22}, 
		    \\
		    A^{\text{\textsf{P}}}_{N + 4 i - 3, i} = \frac{1}{11}, &
		    A^{\text{\textsf{P}}}_{ N + 4 i - 2,i} = \frac{1}{2}, 
		    \\
		    A^{\text{\textsf{P}}}_{ N + 4 i - 1,i} = \frac{3}{4}, &
		    A^{\text{\textsf{P}}}_{ N + 4 i - 0,i} = \frac{1}{16} 
		\end{array}
		$$
		}
% 		$$
% 		$$
% 		A^{\text{\textsf{P}}}_{i, N + 4 i - 3} = \frac{1}{12} 
% 		,\;
% % 		$$
% % 		$$
% 		A^{\text{\textsf{P}}}_{i, N + 4 i - 2} = \frac{1}{8} 
% 		$$
% 		$$
% 		A^{\text{\textsf{P}}}_{i, N + 4 i - 1} = \frac{1}{4} 
% 		$$
% 		$$
% 		A^{\text{\textsf{P}}}_{i, N + 4 i - 0} = \frac{1}{24} 
% 		$$
% 		$$
% 		A^{\text{\textsf{P}}}_{N + 4 i - 2, N + 4 i - 3} = \frac{1}{2} 
% 		$$	
% 		$$
% 		A^{\text{\textsf{P}}}_{N + 4 i - 3, N + 4 i - 2} = \frac{3}{22} 
% 		$$			

% 		$$
% 		A^{\text{\textsf{P}}}_{N + 4 i - 1, N + 4 i - 3} = \frac{1}{4} 
% 		$$	
% 		$$
% 		A^{\text{\textsf{P}}}_{N + 4 i - 3, N + 4 i - 1} = \frac{1}{11} 
% 		$$
		
% 		$$
% 		A^{\text{\textsf{P}}}_{N + 4 i - 0, N + 4 i - 3} = \frac{15}{16} 
% 		$$	
% 		$$
% 		A^{\text{\textsf{P}}}_{N + 4 i - 3, N + 4 i - 0} = \frac{15}{22} 
% 		$$

% 		$$
% 		A^{\text{\textsf{P}}}_{N + 4 i - 3, i} = \frac{1}{11} 
% 		$$		
% 		$$
% 		A^{\text{\textsf{P}}}_{ N + 4 i - 2,i} = \frac{1}{2} 
% 		$$
% 		$$
% 		A^{\text{\textsf{P}}}_{ N + 4 i - 1,i} = \frac{3}{4} 
% 		$$
% 		$$
% 		A^{\text{\textsf{P}}}_{ N + 4 i - 0,i} = \frac{1}{16} 
% 		$$
%		}
%		\STATE{\textbf{normalize} the rows of $A^{\text{\textsf{P}}}$ dividing by their sum}
		\end{algorithmic} 
		}
\end{algorithm}

\vspace{1mm} 
{\centering
\allowdisplaybreaks
\noindent \colorbox{c2!18}{\parbox{0.98\textwidth}{
\begin{lemma}\label{prop:rev_row_sto}
    Let $A\in\mathbb R^{N\times N}$ be a time-reversible matrix such that $\mathcal G(\bar A)$ is a bidirected graph. Then, the matrix $A^\text{\textsf{P}}\in\mathbb R^{5N\times 5N}$ that results from Algorithm~\ref{alg:buildAP_rev} is row stochastic.~\hfill$\circ$  
\end{lemma}
}
}
}
\vspace{1mm}

\begin{proof}
    To show that $A^\text{\textsf{P}}$ is row stochastic, we first analyse the first $N$ rows of $A^\text{\textsf{P}}$. 
    The first $N$ row columns of each of this rows sums up to $\frac{1}{2}$. This happens since the original $A$ is row stochastic and in step~4 of Algorithm~\ref{alg:buildAP_rev}, we divide each element of $A$ by $2$. 
    Subsequently, in row $i$, we further have that the entries different from $0$ are 
    $A^{\text{\textsf{P}}}_{i, N + 4 i - 3},A^{\text{\textsf{P}}}_{i, N + 4 i - 2},A^{\text{\textsf{P}}}_{i, N + 4 i - 1}$ and $A^{\text{\textsf{P}}}_{i, N + 4 i - 0}$. The sum of those entries is $\frac{1}{12}+\frac{1}{8}+\frac{1}{4}+\frac{1}{24}=\frac{1}{2}$. 
    Therefore, each of the $N$ first rows sum up to $1$. 
    Finally, we have to check that the remaining rows also sum up to $1$. 
    This yields by construction. 
    Notice that for each $i=1,\ldots,N$ the sum of the row:  
    \begin{itemize}
        \item $N+4i-3$ is $\frac{1}{11}+\frac{3}{22}+\frac{1}{11}+\frac{15}{22}=1$;
        \item $N+4i-2$ is $\frac{1}{2}+\frac{1}{2}=1$;
        \item $N+4i-1$ is $\frac{3}{4}+\frac{1}{4}=1$;
        \item $N+4i-0$ is $\frac{1}{16}+\frac{15}{16}=1$.\qedhere
    \end{itemize}
\end{proof}

\vspace{1mm} 
{\centering
\allowdisplaybreaks
\noindent \colorbox{c2!18}{\parbox{0.98\textwidth}{
\begin{lemma}\label{prop:rev_is_rev}
    Let $A\in\mathbb R^{N\times N}$ be a time-reversible matrix such that $\mathcal G(\bar A)$ is a bidirected graph. Then, the matrix $A^\text{\textsf{P}}\in\mathbb R^{5N\times 5N}$ that results from Algorithm~\ref{alg:buildAP_rev} is time-reversible.~\hfill$\circ$ 
\end{lemma}
}
}
}
\vspace{1mm}

\begin{proof}
    First, we can check That the submatrix of $A^\text{\textsf{P}}$ corresponding to the first $N$ rows and first $N$ columns is satisfies the equation $r_i A^\text{\textsf{P}}_{ij}=r_j A^\text{\textsf{P}}_{ji}$ for all $i=1,\ldots,N$ with $r_i\geq 0$. 
    Notice that, for the original $A$, we have that for all $i=1\ldots,N$ there is $s\in(\mathbb R^+_0)^N$ such that $s_i A_{ij}=s_j A_{ji}$ and $\displaystyle\sum_{i=1}^N s_i=1$. 
    Hence, it is also true that $s_i \frac{A_{ij}}{2}=s_j \frac{A_{ji}}{2}$ which is the same as $s_i A^\text{\textsf{P}}_{ij}=s_j A^\text{\textsf{P}}_{ji}$. 
    Next, we look at a block corresponding to the additional augmented states of node $i$. 
    Given $s_i$, we define:
    \begin{itemize}
        \item $s_{N+4i-3}=\frac{11}{12}s_i$;
        \item $s_{N+4i-2}=\frac{1}{4}s_i$;
        \item $s_{N+4i-1}=\frac{1}{3}s_i$;
        \item $s_{N+4i-0}=\frac{2}{3}s_i$. 
    \end{itemize}
    Now, replacing in the equation to check reversibility, we have that
    \begin{itemize}
        \item $s_{i}A^{\text{\textsf{P}}}_{i, N + 4 i - 3}=s_{N + 4 i - 3}A^{\text{\textsf{P}}}_{ N + 4 i - 3,i}$ becomes $s_{i}\frac{1}{12}=\frac{11}{12}s_i\frac{1}{11}$ which is equivalent to $s_i=s_i$;
        
        \item $s_{i}A^{\text{\textsf{P}}}_{i, N + 4 i - 2}=s_{N + 4 i - 2}A^{\text{\textsf{P}}}_{ N + 4 i - 2,i}$ becomes $s_{i}\frac{1}{8}=\frac{1}{4}s_i\frac{1}{2}$ which is equivalent to $s_i=s_i$;
        
        \item $s_{i}A^{\text{\textsf{P}}}_{i, N + 4 i - 1}=s_{N + 4 i - 1}A^{\text{\textsf{P}}}_{ N + 4 i - 1,i}$ becomes $s_{i}\frac{1}{4}=\frac{1}{3}s_i\frac{3}{4}$ which is equivalent to $s_i=s_i$;
        
        \item $s_{i}A^{\text{\textsf{P}}}_{i, N + 4 i - 0}=s_{N + 4 i - 0}A^{\text{\textsf{P}}}_{ N + 4 i - 0,i}$ becomes $s_{i}\frac{1}{24}=\frac{2}{3}s_i\frac{1}{16}$ which is equivalent to $s_i=s_i$.
    \end{itemize}  
    However, now we no longer have that $Z=\sum_{i=1}^{5N}s_i=1$, but we can simply normalize each $s_i$ to become $\tilde s_i=\frac{s_i}{Z}$ so that $\sum_{i=1}^{5N}\tilde s_i=1$ without changing the previous equations (as it corresponds to divide both sizes by the same quantity). 
    %and $\displaystyle\sum_{i=1}^N s_i=1$. 
\end{proof}

Next, we show, without loss of generality, that agent $1$, by observing all the agents in the original network and its own augmented states, cannot recover other agents' augmented states. 

\vspace{1mm} 
{\centering
\allowdisplaybreaks
\noindent \colorbox{c2!18}{\parbox{0.98\textwidth}{
\begin{lemma}\label{prop:rev_obs}
     Let $A\in\mathbb R^{N\times N}$ be a time-reversible matrix such that $\mathcal G(\bar A)$ is a bidirected graph. Then, the system $x[k+1]=A^\text{\textsf{P}} x[k]$ and $y[k]=Cx[k]$ , where the matrix $A^\text{\textsf{P}}\in\mathbb R^{5N\times 5N}$ results from Algorithm~\ref{alg:buildAP_rev} and 
\[
C =
\begin{bmatrix}
\bigI_{N+4} & \rvline  & \bigzero_{(N+4) \times (4N-4)} \\
\end{bmatrix}
\]
has the following properties: 
     \begin{enumerate}
         \item it is not observable;
         \item $\left( e_{5N}^j+\displaystyle\sum_{k=N+4j-3}^{N+4j} e_{5N}^{k}\right)\notin span\left( P_{O(A^\text{\textsf{P}},C)}^\lambda\right)$, for all $j\neq i$;
         \item $\displaystyle\left(\frac{1}{{(v_L)}_j} e_{5N}^j+\displaystyle\sum_{k=N+4j-3}^{N+4j} \frac{1}{{(v_L)}_k} e_{5N}^{k}\right)\notin span\left( P_{O(A^\text{\textsf{P}},C)}^\lambda\right)$, for all $j\neq i$,
     \end{enumerate}  
    where $v_L$ is the left-eigenvector of $A^\text{\textsf{P}}$ associated with the eigenvalue $1$.\hfill$\circ$
 \end{lemma}  
 }
 }
 }
 \vspace{1mm}
 
\begin{proof}
    We construct the observabillity matrix $P_{O(A^\text{\textsf{P}},C)}^0$ and compute its span. 
    \[
    \begin{array}{l}
        span\left(P_{O(A^\text{\textsf{P}},C)}^0\right)= \\
        span\left(\begin{bmatrix}
        \bigI_{N+4} & \rvline  & \bigzero_{(N+4) \times (4N-4)} \\ \hline
        \bigzero_{(4N-4)\times (N+4)} & \rvline & \bigR
        \end{bmatrix}
        \right),
        \end{array}
    \]
where $R\in\mathbb R^{(4N-4)\times (4N-4)}$ is $R=diag(D,\ldots,D)$ with
\[
D=
\begin{bmatrix}
 1 & 0 & 0 & 0 \\
 0 & 1 & 0 & \frac{22}{3} \\
 0 & 0 & 1 & -\frac{7}{2} \\
 0 & 0 & 0 & 0 \\
\end{bmatrix}
\]
Hence, $P_{O(A^\text{\textsf{P}},C)}^0$ is not full rank, $(1)$ holds, and $e_{5N}^{N+4j-3}+e_{5N}^{N+4j-2}+e_{5N}^{N+4j-1}+e_{5N}^{N+4j-0}\notin span\left( P_{O(A^\text{\textsf{P}},C)}^\lambda\right)$ for all $j\neq i$, $(2)$ holds. 
Finally, we have that the left-eigenvector $v_L$ associated with the eigenvalue $1$  corresponds to the $s$ obtained in Lemma~\ref{prop:rev_is_rev}. 
We have that 
\[
\left[
\begin{smallmatrix}
 {(v_L)}_{j} & {(v_L)}_{N+4j-3} & {(v_L)}_{N+4j-2} & {(v_L)}_{N+4j-1} & {(v_L)}_{N+4j-0} 
\end{smallmatrix}
\right]
=\]
\[
\left[
\begin{smallmatrix}
 1 & \frac{11}{12} & \frac{1}{4} & \frac{1}{3} & \frac{2}{3} 
\end{smallmatrix}\right].
\]
Therefore, 
$\frac{1}{{(v_L)}_j} e_{5N}^j+\displaystyle\sum_{k=N+4j-3}^{N+4j} \frac{1}{{(v_L)}_k} e_{5N}^{k}= e_{5N}^j+\frac{12}{11}e_{5N}^{N+4j-3}+4 e_{5N}^{N+4j-2}+3e_{5N}^{N+4j-1}+\frac{3}{2}e_{5N}^{N+4j-0}\notin span\left( P_{O(A^\text{\textsf{P}},C)}^\lambda\right)$.
\end{proof}

\begin{algorithm}[H]
		{%\small
		\caption{Solution to problem $\mathbf{P_1^D}$}
		\label{alg:sol_P2}
		\begin{algorithmic}[1]
			\STATE{\textbf{input}: row stochastic dynamics matrix $A\in\mathbb R^{N\times N}$ corresponding to a bidirectional connected graph, and initial agents state $x[0]\in\mathbb R^N$}
			\STATE{\textbf{output}: row stochastic dynamics matrix $A^\text{\textsf{P}}\in\mathbb R^{5N\times 5N}$, and initial agents state $\tilde{x}^{\text{\textsf{P}}}[0]\in\mathbb R^{5N}$}
		\STATE{\textbf{compute} $A^\text{\textsf{P}}$ with Algorithm~\ref{alg:buildAP_rev}}
		\STATE{\textbf{set} $s\in\mathbb R^{5N}$}
		\STATE{\textbf{set} $s_1=1$}
		\STATE{\textbf{set} $\mathcal N=\{2,\ldots,N\}$}
		\STATE{\textbf{set} $\mathcal Y=\{1\}$}
		\WHILE{$\mathcal N\neq \emptyset$}
		    \STATE{\textbf{find} $j\in\mathcal N$ such that $A^\text{\textsf{P}}_{ij}\neq 0$ for some $i\in\mathcal Y$}
		    \STATE{\textbf{set} $s_j=\displaystyle\frac{s_i A^\text{\textsf{P}}_{ij}}{A^\text{\textsf{P}}_{ji}}$}
		    \STATE{\textbf{set} $\mathcal N=\mathcal N\setminus\{j\}$}
		    \STATE{\textbf{set} $\mathcal Y=\mathcal Y\cup\{j\}$}
		\ENDWHILE
		\FOR{i=1,\ldots,N}
		    \STATE{\textbf{set} $s_{N+4i-3}=\frac{11}{12}s_i$}
		    \STATE{\textbf{set} $s_{N+4i-2}=\frac{1}{4}s_i$}
		    \STATE{\textbf{set} $s_{N+4i-1}=\frac{1}{3}s_i$}
		    \STATE{\textbf{set} $s_{N+4i-0}=\frac{2}{3}s_i$}
		\ENDFOR
        \STATE{\textbf{set} 
        \[
        v_L=
        \begin{bmatrix}
           \frac{s_1}{Z} & \ldots & \frac{s_{5N}}{Z}
        \end{bmatrix},\text{ with }Z=\sum_{i=1}^{5N}s_i
        \]
        }
        \FOR{i=1,\ldots,N}
		\STATE{\textbf{choose} $\alpha_1,\alpha_2,\alpha_3,\alpha_4>0$ such that 
		\[
		    \frac{1}{5}\sum_{j=1}^4 \alpha_j = x_i[0]
		\]
		}
		\STATE{\textbf{set} $\tilde{x}_i^{\text{\textsf{P}}}[0]\qquad\;\; = 0$}
		\STATE{\textbf{set} $\tilde{x}_{N+4i-3}^{\text{\textsf{P}}}[0] = \frac{Z}{s_{N+4i-3}}\alpha_1$}
		\STATE{\textbf{set} $\tilde{x}_{N+4i-2}^{\text{\textsf{P}}}[0] = \frac{Z}{s_{N+4i-2}}\alpha_2$}
		\STATE{\textbf{set} $\tilde{x}_{N+4i-1}^{\text{\textsf{P}}}[0] = \frac{Z}{s_{N+4i-1}}\alpha_3$}
		\STATE{\textbf{set} $\tilde{x}_{N+4i-0}^{\text{\textsf{P}}}[0] = \frac{Z}{s_{N+4i-0}}\alpha_4$}
		\ENDFOR
		\end{algorithmic} 
		}
\end{algorithm}

\vspace{1mm} 
{\centering
\allowdisplaybreaks
\noindent \colorbox{c3!10}{\parbox{0.98\textwidth}{
\begin{theorem}\label{th:main_P2}
    Algorithm~\ref{alg:sol_P2} produces a solution to problem $\mathbf{P_1^D}$.~\hfill$\circ$
\end{theorem}
}
}
}
\vspace{1mm}

\begin{proof}
    By Lemma~\ref{prop:rev_obs}, the matrix obtained with Algorithm~\ref{alg:buildAP_rev} ensures that an agent of the original system is not able to recover the initial state of other original agents nor the combination of other original agents with their augmented states. 
    By Lemma~\ref{prop:rev_row_sto}, we have that the matrix obtained with Algorithm~\ref{alg:buildAP_rev} is row stochastic. 
    Moreover, by Lemma~\ref{prop:rev_obs}, we ensure that an agent cannot recover the initial value of other agents. 
    Finally, by Proposition~\ref{th:consensus_reqs}, we have ensured that the output dynamics matrix $A^\text{\textsf{P}}$ and the output new initial state $\tilde{x}^{\text{\textsf{P}}}[0]$ ensure that the system reaches average consensus.
\end{proof}

\vspace{1mm} 
{\centering
\allowdisplaybreaks
\noindent \colorbox{c3!10}{\parbox{0.98\textwidth}{
\begin{theorem}\label{th:main_P2_complexity}
    The Algorithm~\ref{alg:sol_P2} has computational time-complexity of $\mathcal O(N^2)$, where $N\times N$ is the dimension of the input matrix $A$.\hfill$\circ$
\end{theorem} 
}
}
}
\vspace{1mm}

\begin{proof}
    The first main step of Algorithm~\ref{alg:sol_P2} is to compute $A^\text{\textsf{P}}$ with Algorithm~\ref{alg:buildAP_rev}, which has computation time-complexity of $\mathcal O(N^2)$ since we build a matrix with dimension $5N\times 5N$. 
    The second main computational step of the algorithm is the while loop, steps~8--13, which also has computation time-complexity of $\mathcal O(N^2)$ because in each iteration we remove an element of the set $\mathcal N$ (with initial size $|\mathcal N|=N-1$), and for each iteration we search for an element in $\mathcal Y$ (which maximum size is $|\mathcal Y|=N$), yielding $\mathcal O(N^2)$ operations. 
    Next, the for loop in steps~14--19 has computational time-complexity of $\mathcal O(N)$. 
    Similarly, step~20 has computational time-complexity of $\mathcal O(N)$. 
    Finally, the last for loop, steps~21--28, has also computational time-complexity of $\mathcal O(N)$. 
    The total computational time-complexity is given as the sum of the previous terms, which is $\mathcal O(N^2)+\mathcal O(N^2)+\mathcal O(N)+\mathcal O(N)+\mathcal O(N)=\mathcal O(N^2)$. 
\end{proof}

The next corollary follows immediately from the fact that the computation space-complexity cannot be larger than the time-complexity (see Section 4.1 of~\cite{arora2009computational}).

\vspace{1mm} 
{\centering
\allowdisplaybreaks
\noindent \colorbox{c4!15}{\parbox{0.98\textwidth}{
\begin{corollary}\label{cor:main_P2_space}
    The Algorithm~\ref{alg:sol_P2} has computational space-complexity of $\mathcal O(N^2)$, where $N\times N$ is the dimension of the input matrix $A$.\hfill$\circ$
\end{corollary}
}
}
}
\vspace{1mm}

In the next section, we present some illustrative examples of the proposed methods.

\section{Illustrative examples}\label{sec:ill_exp}

\subsection*{Example 1}
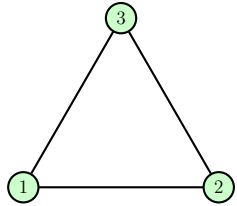
\begin{figure}[H]
\begin{center}
\begin{tikzpicture}[scale=.6, transform shape,node distance=1.5cm]
\begin{scope}[every node/.style={circle,thick,draw},square/.style={regular polygon,regular polygon sides=4}]
\node[fill=green!20] (1) at (-2.16506,-1.25) {\large $1$};
\node[fill=green!20] (2) at (2.16506,-1.25) {\large $2$};
\node[fill=green!20] (3) at (0.,2.5) {\large $3$};
\end{scope}
\begin{scope}[>={Stealth[black]},
              every edge/.style={draw=black, thick}]
\path [-] (1) edge node {} (2);
% \path [->] (2) edge[bend right=15] node {} (1);
\path [-] (2) edge node {} (3);
% \path [->] (3) edge[bend right=15] node {} (2);
\path [-] (1) edge node {} (3);
% \path [->] (3) edge[bend right=15] node {} (1);
\end{scope}
\end{tikzpicture}
\caption{Graph representation of a cycle network with 3 agents.}
\label{fig:graph_cycle_3_undirected}
\end{center}
\end{figure}

The first example we use is a cycle graph with 3 agents, corresponding to the matrix 
$
A=\left[\begin{smallmatrix}
  0 & \frac{1}{2} & \frac{1}{2} \\
 \frac{1}{2} & 0 & \frac{1}{2} \\
 \frac{1}{2} & \frac{1}{2} & 0 \\
\end{smallmatrix}\right],
$ 
with digraph representation depicted in Figure~\ref{fig:graph_cycle_3_undirected} and initial agents states $x[0]=\left[
\begin{smallmatrix}
 \frac{1}{2} & \frac{1}{3} & \frac{1}{5}  \\
\end{smallmatrix}
\right]^\intercal$. 
Notice that $A$ is trivially time-reversible using the vector $s=\left[\begin{smallmatrix}
   \frac{1}{3} & \frac{1}{3} & \frac{1}{3}
\end{smallmatrix}\right]^\intercal$ since $A$ is also column stochastic. 

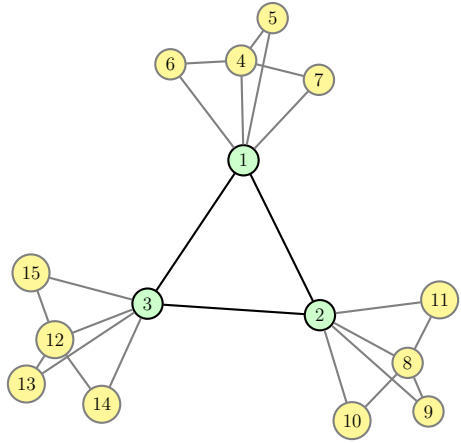
\begin{figure}[H]
\begin{center}
\begin{tikzpicture}[scale=.6, transform shape,node distance=1.5cm]
\begin{scope}[every node/.style={circle,thick,draw},square/.style={regular polygon,regular polygon sides=4}]
\node[fill=green!20] (1) at (4.81111,5.77276) {\large $1$};
\node[fill=green!20] (2) at (6.50401,2.3382) {\large $2$};
\node[fill=green!20] (3) at (2.68447,2.59388) {\large $3$};
\node[draw=gray,fill=yellow!50] (4) at (4.75906,7.98496) {\large $4$};
\node[draw=gray,fill=yellow!50] (5) at (5.45672,8.92496) {\large $5$};
\node[draw=gray,fill=yellow!50] (6) at (3.1916,7.90175) {\large $6$};
\node[draw=gray,fill=yellow!50] (7) at (6.48169,7.55724) {\large $7$};
\node[draw=gray,fill=yellow!50] (8) at (8.44719,1.28104) {\large $8$};
\node[draw=gray,fill=yellow!50] (9) at (8.90895,0.202946) {\large $9$};
\node[draw=gray,fill=yellow!50] (10) at (7.21577,0.) {\large $10$};
\node[draw=gray,fill=yellow!50] (11) at (9.15675,2.67711) {\large $11$};
\node[draw=gray,fill=yellow!50] (12) at (0.630753,1.79509) {\large $12$};
\node[draw=gray,fill=yellow!50] (13) at (0.,0.812305) {\large $13$};
\node[draw=gray,fill=yellow!50] (14) at (1.66965,0.362517) {\large $14$};
\node[draw=gray,fill=yellow!50] (15) at (0.102724,3.27607) {\large $15$};
\end{scope}
\begin{scope}[>={Stealth[black]},
              every edge/.style={draw=black, thick}]
\path [-] (1) edge node {} (2);
\path [-] (1) edge node {} (3);
\path [-] (1) edge[gray] node {} (4);
\path [-] (1) edge[gray] node {} (5);
\path [-] (1) edge[gray] node {} (6);
\path [-] (1) edge[gray] node {} (7);
\path [-] (2) edge node {} (3);
\path [-] (2) edge[gray] node {} (8);
\path [-] (2) edge[gray] node {} (9);
\path [-] (2) edge[gray] node {} (10);
\path [-] (2) edge[gray] node {} (11);
\path [-] (3) edge[gray] node {} (12);
\path [-] (3) edge[gray] node {} (13);
\path [-] (3) edge[gray] node {} (14);
\path [-] (3) edge[gray] node {} (15);
\path [-] (4) edge[gray] node {} (5);
\path [-] (4) edge[gray] node {} (6);
\path [-] (4) edge[gray] node {} (7);
\path [-] (8) edge[gray] node {} (9);
\path [-] (8) edge[gray] node {} (10);
\path [-] (8) edge[gray] node {} (11);
\path [-] (12) edge[gray] node {} (13);
\path [-] (12) edge[gray] node {} (14);
\path [-] (12) edge[gray] node {} (15);
\end{scope}
\end{tikzpicture}
\caption{Network obtained with Algorithm~\ref{alg:sol_P2} with input $A$, $\mathcal G(A^\text{\textsf{P}})$.}
\label{fig:cycle_3_rev_priv}
\end{center}
\end{figure}

Next, we apply Algorithm~\ref{alg:sol_P2} and obtain the matrix 
\[
\begin{array}{l}
A^\text{\textsf{P}}= \\
\left[
\begin{smallmatrix}
 0 & \frac{1}{4} & \frac{1}{4} & \frac{1}{12} & \frac{1}{8} & \frac{1}{4} & \frac{1}{24} & 0 & 0 & 0 & 0 & 0 & 0 & 0 & 0 \\
 \frac{1}{4} & 0 & \frac{1}{4} & 0 & 0 & 0 & 0 & \frac{1}{12} & \frac{1}{8} & \frac{1}{4} & \frac{1}{24} & 0 & 0 & 0 & 0 \\
 \frac{1}{4} & \frac{1}{4} & 0 & 0 & 0 & 0 & 0 & 0 & 0 & 0 & 0 & \frac{1}{12} & \frac{1}{8} & \frac{1}{4} & \frac{1}{24} \\
 \frac{1}{11} & 0 & 0 & 0 & \frac{3}{22} & \frac{1}{11} & \frac{15}{22} & 0 & 0 & 0 & 0 & 0 & 0 & 0 & 0 \\
 \frac{1}{2} & 0 & 0 & \frac{1}{2} & 0 & 0 & 0 & 0 & 0 & 0 & 0 & 0 & 0 & 0 & 0 \\
 \frac{3}{4} & 0 & 0 & \frac{1}{4} & 0 & 0 & 0 & 0 & 0 & 0 & 0 & 0 & 0 & 0 & 0 \\
 \frac{1}{16} & 0 & 0 & \frac{15}{16} & 0 & 0 & 0 & 0 & 0 & 0 & 0 & 0 & 0 & 0 & 0 \\
 0 & \frac{1}{11} & 0 & 0 & 0 & 0 & 0 & 0 & \frac{3}{22} & \frac{1}{11} & \frac{15}{22} & 0 & 0 & 0 & 0 \\
 0 & \frac{1}{2} & 0 & 0 & 0 & 0 & 0 & \frac{1}{2} & 0 & 0 & 0 & 0 & 0 & 0 & 0 \\
 0 & \frac{3}{4} & 0 & 0 & 0 & 0 & 0 & \frac{1}{4} & 0 & 0 & 0 & 0 & 0 & 0 & 0 \\
 0 & \frac{1}{16} & 0 & 0 & 0 & 0 & 0 & \frac{15}{16} & 0 & 0 & 0 & 0 & 0 & 0 & 0 \\
 0 & 0 & \frac{1}{11} & 0 & 0 & 0 & 0 & 0 & 0 & 0 & 0 & 0 & \frac{3}{22} & \frac{1}{11} & \frac{15}{22} \\
 0 & 0 & \frac{1}{2} & 0 & 0 & 0 & 0 & 0 & 0 & 0 & 0 & \frac{1}{2} & 0 & 0 & 0 \\
 0 & 0 & \frac{3}{4} & 0 & 0 & 0 & 0 & 0 & 0 & 0 & 0 & \frac{1}{4} & 0 & 0 & 0 \\
 0 & 0 & \frac{1}{16} & 0 & 0 & 0 & 0 & 0 & 0 & 0 & 0 & \frac{15}{16} & 0 & 0 & 0 \\
\end{smallmatrix}
\right],
\end{array}
\]
which corresponds to the graph representation depicted in Figure~\ref{fig:cycle_3_rev_priv}, the vector (left-eigenvector associated with the eigenvalue 1) $s=\left[
\begin{smallmatrix}
 \frac{2}{19} & \frac{2}{19} & \frac{2}{19} & \frac{11}{114} & \frac{1}{38} & \frac{2}{57} & \frac{4}{57} & \frac{11}{114} & \frac{1}{38} & \frac{2}{57} & \frac{4}{57} & \frac{11}{114} & \frac{1}{38} & \frac{2}{57} & \frac{4}{57} \\
\end{smallmatrix}
\right]^\intercal$, and the new vector of initial agents' states 
\[
\begin{split}
\tilde{x}^{\text{\textsf{P}}}[0]=\left[
\begin{smallmatrix}
 0 & 0 & 0 & 0.5894 & 0.8522 & 0.7909 & 0.8495 & 0.3415 & 0.7026 & 1.1778\end{smallmatrix}\right. \\
 \left.\begin{smallmatrix} 0.2615 & 0.0254 & 0.0305 & 1.1357 & 0.3357 \\
\end{smallmatrix}
\right]^\intercal.
\end{split}
\]

% Hence, another property that we have to ensure, is that agent 1 is not able to recover a combination of other agents copies that is given by the corresponding values in $s$. 
% Regarding agent 2, the agent 1 can recover $e_{15}^{8}+e_{15}^{9}+e_{15}^{10}+2 e_{15}^{11}$ with coefficients that are not multiple of the respective entries in $s$. 
% Analogously, regarding agent 3, the agent 1 can recover $e_{15}^{12}+e_{15}^{13}+e_{15}^{14}+2 e_{15}^{15}$ with coefficients that are not multiple of the respective entries in $s$. 

\begin{figure}[!ht]
\centering
\input{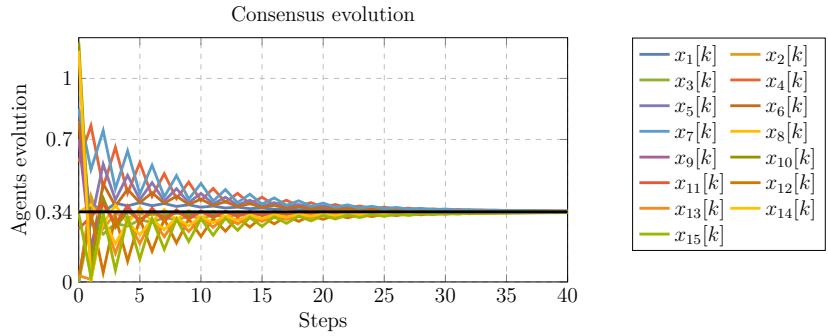}
\caption{Agents states evolution with $A^{\text{\textsf{P}}}$ and $\tilde{x}^{\text{\textsf{P}}}[0]$. The black line represents the average of the initial agents' states $x[0]$.}
\label{fig:cycle3_consensus_rev}
\end{figure}

In Figure~\ref{fig:cycle3_consensus_rev}, we depict the agents state evolution under network $\mathcal G(A^{\text{\textsf{P}}})$, weight matrix $A^{\text{\textsf{P}}}$, and vector of initial agents' states $\tilde{x}^{\text{\textsf{P}}}[0]$. 
Furthermore, we verified that the required conditions of Lemma~\ref{prop:rev_obs} hold, as envisioned. 

\subsection*{Example 2}

\begin{figure}[H]
\begin{center}
\begin{tikzpicture}[scale=.6, transform shape,node distance=1.5cm]
\begin{scope}[every node/.style={circle,thick,draw},square/.style={regular polygon,regular polygon sides=4}]
\node[fill=green!20] (1) at (1.89839,3.74605) {\large $1$};
\node[fill=green!20] (2) at (1.64044,5.52748) {\large $2$};
\node[fill=green!20] (3) at (2.35331,6.69686) {\large $3$};
\node[fill=green!20] (4) at (5.29166,7.62145) {\large $4$};
\node[fill=green!20] (5) at (4.55086,5.82772) {\large $5$};
\node[fill=green!20] (6) at (7.34204,4.64889) {\large $6$};
\node[fill=green!20] (7) at (9.9213,6.00252) {\large $7$};
\node[fill=green!20] (8) at (8.39954,6.68311) {\large $8$};
\node[fill=green!20] (9) at (9.33409,3.68274) {\large $9$};
\node[fill=green!20] (10) at (7.05564,2.13599) {\large $10$};
\node[fill=green!20] (11) at (4.01028,1.55303) {\large $11$};
\end{scope}
\begin{scope}[>={Stealth[black]},
              every edge/.style={draw=black, thick}]
\path [-] (1) edge node {} (2);
\path [-] (1) edge node {} (3);
\path [-] (1) edge node {} (11);
\path [-] (2) edge node {} (3);
\path [-] (2) edge node {} (5);
\path [-] (3) edge node {} (4);
\path [-] (4) edge node {} (5);
\path [-] (4) edge node {} (8);
\path [-] (5) edge node {} (6);
\path [-] (6) edge node {} (7);
\path [-] (6) edge node {} (10);
\path [-] (7) edge node {} (8);
\path [-] (8) edge node {} (9);
\path [-] (9) edge node {} (10);
\path [-] (10) edge node {} (11);
\end{scope}
\end{tikzpicture}
\caption{Graph representation of $A$, $\mathcal G(A)$.}
\label{fig:graph_rev_4nodes}
\end{center}
\end{figure}
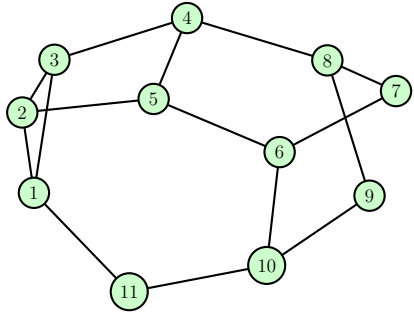

The second example is the network with 11 agents depicted in Figure~\ref{fig:graph_rev_4nodes}, with matrix 
\[
A=\left[
\begin{smallmatrix}
 0 & \frac{1}{3} & \frac{1}{3} & 0 & 0 & 0 & 0 & 0 & 0 & 0 & \frac{1}{3} \\
 \frac{1}{3} & 0 & \frac{1}{3} & 0 & \frac{1}{3} & 0 & 0 & 0 & 0 & 0 & 0 \\
 \frac{1}{3} & \frac{1}{3} & 0 & \frac{1}{3} & 0 & 0 & 0 & 0 & 0 & 0 & 0 \\
 0 & 0 & \frac{1}{3} & 0 & \frac{1}{3} & 0 & 0 & \frac{1}{3} & 0 & 0 & 0 \\
 0 & \frac{1}{3} & 0 & \frac{1}{3} & 0 & \frac{1}{3} & 0 & 0 & 0 & 0 & 0 \\
 0 & 0 & 0 & 0 & \frac{1}{3} & 0 & \frac{1}{3} & 0 & 0 & \frac{1}{3} & 0 \\
 0 & 0 & 0 & 0 & 0 & \frac{1}{2} & 0 & \frac{1}{2} & 0 & 0 & 0 \\
 0 & 0 & 0 & \frac{1}{3} & 0 & 0 & \frac{1}{3} & 0 & \frac{1}{3} & 0 & 0 \\
 0 & 0 & 0 & 0 & 0 & 0 & 0 & \frac{1}{2} & 0 & \frac{1}{2} & 0 \\
 0 & 0 & 0 & 0 & 0 & \frac{1}{3} & 0 & 0 & \frac{1}{3} & 0 & \frac{1}{3} \\
 \frac{1}{2} & 0 & 0 & 0 & 0 & 0 & 0 & 0 & 0 & \frac{1}{2} & 0 \\
\end{smallmatrix}
\right]
\]
and initial agents states \\$x[0]=\left[
\begin{smallmatrix}
 0.1 & 0.3 & 0.6 & 0.43 & 0.85 & 0.9 & 0.45 & 0.11 & 0.06 & 0.51 & 0.13 \\
\end{smallmatrix}
\right]^\intercal$. 
Notice that $A$ is time-reversible using the vector $s=\left[\begin{smallmatrix}
 \frac{1}{11} & \ldots & \frac{1}{11} \\
\end{smallmatrix}\right]^\intercal$ and $A$ is row stochastic. 

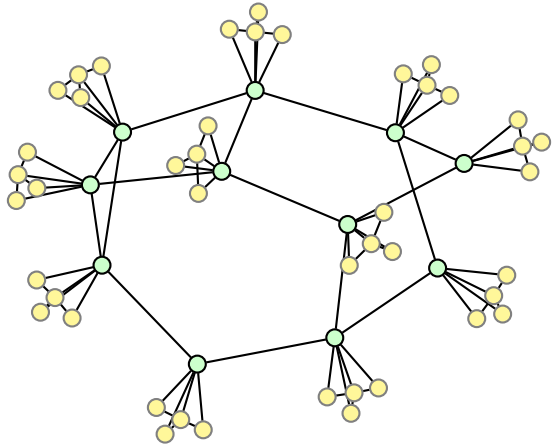
\begin{figure}[H]
\begin{center}
\begin{tikzpicture}[scale=.6, transform shape,node distance=1.5cm]
\begin{scope}[every node/.style={circle,thick,draw},square/.style={regular polygon,regular polygon sides=4}]
\node[fill=green!20] (1) at (1.89839,3.74605) {\,};
\node[fill=green!20] (2) at (1.64044,5.52748) {\,};
\node[fill=green!20] (3) at (2.35331,6.69686) {\,};
\node[fill=green!20] (4) at (5.29166,7.62145) {\,};
\node[fill=green!20] (5) at (4.55086,5.82772) {\,};
\node[fill=green!20] (6) at (7.34204,4.64889) {\,};
\node[fill=green!20] (7) at (9.9213,6.00252) {\,};
\node[fill=green!20] (8) at (8.39954,6.68311) {\,};
\node[fill=green!20] (9) at (9.33409,3.68274) {\,};
\node[fill=green!20] (10) at (7.05564,2.13599) {\,};
\node[fill=green!20] (11) at (4.01028,1.55303) {\,};
\node[draw=gray,fill=yellow!50] (12) at (0.857009,3.0279) {\,};
\node[draw=gray,fill=yellow!50] (13) at (0.451019,3.41756) {\,};
\node[draw=gray,fill=yellow!50] (14) at (0.534933,2.70356) {\,};
\node[draw=gray,fill=yellow!50] (15) at (1.23536,2.57485) {\,};
\node[draw=gray,fill=yellow!50] (16) at (0.0391827,5.75373) {\,};
\node[draw=gray,fill=yellow!50] (17) at (0.454116,5.45447) {\,};
\node[draw=gray,fill=yellow!50] (18) at (0.246642,6.25259) {\,};
\node[draw=gray,fill=yellow!50] (19) at (0.,5.18001) {\,};
\node[draw=gray,fill=yellow!50] (20) at (1.37864,7.96804) {\,};
\node[draw=gray,fill=yellow!50] (21) at (1.42743,7.46217) {\,};
\node[draw=gray,fill=yellow!50] (22) at (1.88445,8.16411) {\,};
\node[draw=gray,fill=yellow!50] (23) at (0.921747,7.62305) {\,};
\node[draw=gray,fill=yellow!50] (24) at (5.29356,8.91924) {\,};
\node[draw=gray,fill=yellow!50] (25) at (5.89227,8.85964) {\,};
\node[draw=gray,fill=yellow!50] (26) at (4.71871,8.97913) {\,};
\node[draw=gray,fill=yellow!50] (27) at (5.36475,9.36012) {\,};
\node[draw=gray,fill=yellow!50] (28) at (3.99828,6.20823) {\,};
\node[draw=gray,fill=yellow!50] (29) at (3.53043,5.95673) {\,};
\node[draw=gray,fill=yellow!50] (30) at (4.03446,5.33757) {\,};
\node[draw=gray,fill=yellow!50] (31) at (4.2482,6.83541) {\,};
\node[draw=gray,fill=yellow!50] (32) at (7.86389,4.22279) {\,};
\node[draw=gray,fill=yellow!50] (33) at (8.33371,4.04813) {\,};
\node[draw=gray,fill=yellow!50] (34) at (7.37631,3.74199) {\,};
\node[draw=gray,fill=yellow!50] (35) at (8.14353,4.91835) {\,};
\node[draw=gray,fill=yellow!50] (36) at (11.215,6.38649) {\,};
\node[draw=gray,fill=yellow!50] (37) at (11.379,5.81615) {\,};
\node[draw=gray,fill=yellow!50] (38) at (11.6391,6.47493) {\,};
\node[draw=gray,fill=yellow!50] (39) at (11.1172,6.96886) {\,};
\node[draw=gray,fill=yellow!50] (40) at (9.09967,7.73628) {\,};
\node[draw=gray,fill=yellow!50] (41) at (9.60296,7.50914) {\,};
\node[draw=gray,fill=yellow!50] (42) at (8.57339,7.98894) {\,};
\node[draw=gray,fill=yellow!50] (43) at (9.19665,8.19307) {\,};
\node[draw=gray,fill=yellow!50] (44) at (10.5765,3.06996) {\,};
\node[draw=gray,fill=yellow!50] (45) at (10.2011,2.5603) {\,};
\node[draw=gray,fill=yellow!50] (46) at (10.8663,3.52421) {\,};
\node[draw=gray,fill=yellow!50] (47) at (10.7752,2.66269) {\,};
\node[draw=gray,fill=yellow!50] (48) at (7.48235,0.915609) {\,};
\node[draw=gray,fill=yellow!50] (49) at (7.41497,0.461904) {\,};
\node[draw=gray,fill=yellow!50] (50) at (8.02404,1.02365) {\,};
\node[draw=gray,fill=yellow!50] (51) at (6.88694,0.825126) {\,};
\node[draw=gray,fill=yellow!50] (52) at (3.64222,0.321883) {\,};
\node[draw=gray,fill=yellow!50] (53) at (3.29574,0.) {\,};
\node[draw=gray,fill=yellow!50] (54) at (4.14109,0.0957859) {\,};
\node[draw=gray,fill=yellow!50] (55) at (3.1023,0.590875) {\,};
\end{scope}
\begin{scope}[on background layer,
              every edge/.style={draw=black, thick}]
\path [-] (1) edge node {} (2);
\path [-] (1) edge node {} (3);
\path [-] (1) edge node {} (11);
\path [-] (1) edge node {} (12);
\path [-] (1) edge node {} (13);
\path [-] (1) edge node {} (14);
\path [-] (1) edge node {} (15);
\path [-] (2) edge node {} (3);
\path [-] (2) edge node {} (5);
\path [-] (2) edge node {} (16);
\path [-] (2) edge node {} (17);
\path [-] (2) edge node {} (18);
\path [-] (2) edge node {} (19);
\path [-] (3) edge node {} (4);
\path [-] (3) edge node {} (20);
\path [-] (3) edge node {} (21);
\path [-] (3) edge node {} (22);
\path [-] (3) edge node {} (23);
\path [-] (4) edge node {} (5);
\path [-] (4) edge node {} (8);
\path [-] (4) edge node {} (24);
\path [-] (4) edge node {} (25);
\path [-] (4) edge node {} (26);
\path [-] (4) edge node {} (27);
\path [-] (5) edge node {} (6);
\path [-] (5) edge node {} (28);
\path [-] (5) edge node {} (29);
\path [-] (5) edge node {} (30);
\path [-] (5) edge node {} (31);
\path [-] (6) edge node {} (7);
\path [-] (6) edge node {} (10);
\path [-] (6) edge node {} (32);
\path [-] (6) edge node {} (33);
\path [-] (6) edge node {} (34);
\path [-] (6) edge node {} (35);
\path [-] (7) edge node {} (8);
\path [-] (7) edge node {} (36);
\path [-] (7) edge node {} (37);
\path [-] (7) edge node {} (38);
\path [-] (7) edge node {} (39);
\path [-] (8) edge node {} (9);
\path [-] (8) edge node {} (40);
\path [-] (8) edge node {} (41);
\path [-] (8) edge node {} (42);
\path [-] (8) edge node {} (43);
\path [-] (9) edge node {} (10);
\path [-] (9) edge node {} (44);
\path [-] (9) edge node {} (45);
\path [-] (9) edge node {} (46);
\path [-] (9) edge node {} (47);
\path [-] (10) edge node {} (11);
\path [-] (10) edge node {} (48);
\path [-] (10) edge node {} (49);
\path [-] (10) edge node {} (50);
\path [-] (10) edge node {} (51);
\path [-] (11) edge node {} (52);
\path [-] (11) edge node {} (53);
\path [-] (11) edge node {} (54);
\path [-] (11) edge node {} (55);
\path [-] (12) edge node {} (13);
\path [-] (12) edge node {} (14);
\path [-] (12) edge node {} (15);
\path [-] (16) edge node {} (17);
\path [-] (16) edge node {} (18);
\path [-] (16) edge node {} (19);
\path [-] (20) edge node {} (21);
\path [-] (20) edge node {} (22);
\path [-] (20) edge node {} (23);
\path [-] (24) edge node {} (25);
\path [-] (24) edge node {} (26);
\path [-] (24) edge node {} (27);
\path [-] (28) edge node {} (29);
\path [-] (28) edge node {} (30);
\path [-] (28) edge node {} (31);
\path [-] (32) edge node {} (33);
\path [-] (32) edge node {} (34);
\path [-] (32) edge node {} (35);
\path [-] (36) edge node {} (37);
\path [-] (36) edge node {} (38);
\path [-] (36) edge node {} (39);
\path [-] (40) edge node {} (41);
\path [-] (40) edge node {} (42);
\path [-] (40) edge node {} (43);
\path [-] (44) edge node {} (45);
\path [-] (44) edge node {} (46);
\path [-] (44) edge node {} (47);
\path [-] (48) edge node {} (49);
\path [-] (48) edge node {} (50);
\path [-] (48) edge node {} (51);
\path [-] (52) edge node {} (53);
\path [-] (52) edge node {} (54);
\path [-] (52) edge node {} (55);
\end{scope}
\end{tikzpicture}

\caption{Network obtained with Algorithm~\ref{alg:sol_P2} with input $A$, $\mathcal G(A^\text{\textsf{P}})$ (we omit the nodes' labels for visualization purposes).}
\label{fig:graph_rev_4nodes_priv}
\end{center}
\end{figure}

Next, we apply Algorithm~\ref{alg:sol_P2} and obtain a matrix $A^{\text{\textsf{P}}}$ that corresponds to the network in Figure~\ref{fig:graph_rev_4nodes_priv}, the vector (left-eigenvector associated with the eigenvalue 1) 
\[
\begin{split}
s=\left[
\begin{smallmatrix}
 \frac{3}{95} & \frac{3}{95} & \frac{3}{95} & \frac{3}{95} & \frac{3}{95} & \frac{3}{95} & \frac{2}{95} & \frac{3}{95} & \frac{2}{95} & \frac{3}{95} & \frac{2}{95} & \frac{11}{380} & \frac{3}{380} & \frac{1}{95} & \frac{2}{95} \end{smallmatrix}\right. \\
  \left.\begin{smallmatrix}
  \frac{11}{380} & \frac{3}{380} & \frac{1}{95} & \frac{2}{95} & \frac{11}{380} & \frac{3}{380} & \frac{1}{95} & \frac{2}{95} & \frac{11}{380} & \frac{3}{380} & \frac{1}{95} & \frac{2}{95} & \frac{11}{380} & \frac{3}{380} & \frac{1}{95} 
  \end{smallmatrix}\right. \\
  \left.\begin{smallmatrix}
  \frac{2}{95} & \frac{11}{380} & \frac{3}{380} & \frac{1}{95} & \frac{2}{95} & \frac{11}{570} & \frac{1}{190} & \frac{2}{285} & \frac{4}{285} & \frac{11}{380} & \frac{3}{380} & \frac{1}{95} & \frac{2}{95} & \frac{11}{570} 
  \end{smallmatrix}\right. \\
  \left.\begin{smallmatrix}
  \frac{1}{190} & \frac{2}{285} & \frac{4}{285} & \frac{11}{380} & \frac{3}{380} & \frac{1}{95} & \frac{2}{95} & \frac{11}{570} & \frac{1}{190} & \frac{2}{285} & \frac{4}{285} 
\end{smallmatrix}
\right]^\intercal,
\end{split}
\]
and the new vector of initial agents' states 
\[
\begin{split}
    \tilde{x}^{\text{\textsf{P}}}[0]=\left[
\begin{smallmatrix}
 0 & 0 & 0 & 0 & 0 & 0 & 0 & 0 & 0 & 0 & 0 & 0.0789501 & 0.358511 & 0.0528758 & 0.162382   \end{smallmatrix}\right. \\
  \left.\begin{smallmatrix} 0.143901 & 0.667626 & 0.9161 & 0.389181 & 0.772374 & 2.01193 & 1.53623 \end{smallmatrix}\right. \\
  \left.\begin{smallmatrix} 0.00630652 & 0.554657 & 0.965073 & 0.479013 & 0.492756 & 0.137887 & 2.89792 \end{smallmatrix}\right. \\
  \left.\begin{smallmatrix} 2.14223  & 1.32303 & 1.04491 & 0.596246 & 2.74642 & 0.852806 & 0.792219 \end{smallmatrix}\right. \\
  \left.\begin{smallmatrix} 0.802344 & 2.86735 & 0.0909197 & 0.106131 & 0.280592 & 0.173293 & 0.137202 \end{smallmatrix}\right. \\
  \left.\begin{smallmatrix} 0.0931979 & 0.560345   & 0.0542549 & 0.0232326 & 0.650005 & 2.73045 & 0.448699 \end{smallmatrix}\right. \\
  \left.\begin{smallmatrix} 0.0602469 & 0.0534166 & 0.548742 & 0.437764 & 0.343937
\end{smallmatrix}
\right]^\intercal.
\end{split}
\]

% NOVO ESTADO INICIAL
% \left(
% \begin{array}{ccccccccccccccccccccccccccccccccccccccccccccccccccccccc}
%  0 & 0 & 0 & 0 & 0 & 0 & 0 & 0 & 0 & 0 & 0 & 0.0789501 & 0.358511 & 0.0528758 & 0.162382 & 0.143901 & 0.667626 & 0.9161 & 0.389181 & 0.772374 & 2.01193 & 1.53623 & 0.00630652 & 0.554657 & 0.965073 & 0.479013 & 0.492756 & 0.137887 & 2.89792 & 2.14223 & 1.32303 & 1.04491 & 0.596246 & 2.74642 & 0.852806 & 0.792219 & 0.802344 & 2.86735 & 0.0909197 & 0.106131 & 0.280592 & 0.173293 & 0.137202 & 0.0931979 & 0.560345 & 0.0542549 & 0.0232326 & 0.650005 & 2.73045 & 0.448699 & 0.0602469 & 0.0534166 & 0.548742 & 0.437764 & 0.343937 \\
% \end{array}
% \right)

\begin{figure}[H]
\centering
\input{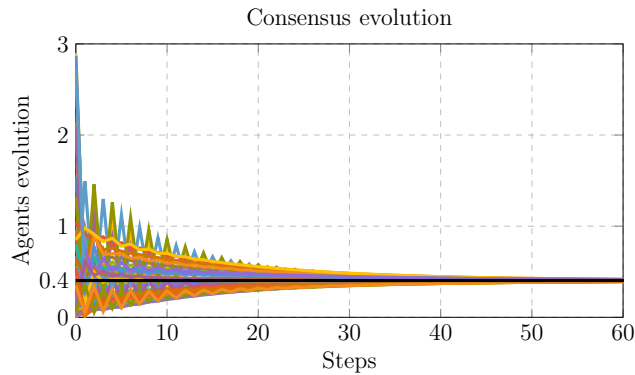}
\caption{Agents states evolution with $A^{\text{\textsf{P}}}$ and $\tilde{x}^{\text{\textsf{P}}}[0]$, where each color corresponds to the evolution of an agent state. The black line represents the average of the initial agents' states $x[0]$.}
\label{fig:graph_rev_4nodes_rev}
\end{figure}

 In Figure~\ref{fig:graph_rev_4nodes_rev}, we depict the agents state evolution under network $\mathcal G(A^{\text{\textsf{P}}})$, weight matrix $A^{\text{\textsf{P}}}$, and vector of initial agents' states $\tilde{x}^{\text{\textsf{P}}}[0]$. 
 Also, we verified that the required conditions of Lemma~\ref{prop:rev_obs} hold, as envisioned.

\section{Conclusions}\label{sec:conclusions}

In this paper, we tackled the problem of keeping the initial agents' states private in an average consensus network. 
We designed two algorithms to reach average consensus of a network of multi-agents that allows to keep their initial states private. 
The principal idea of the algorithms is to augment the state of each agent with a set of augmented states nodes that only interact locally. Then, each agent can choose how to distribute its initial state to the augmented states while ensuring that the other agents cannot recover the initial information. 

The first proposed algorithm achieves the privacy goal by designing three local augmented states of each agent and the weights of the adjacency matrix. 
A computational step required is to compute the left-eigenvector associated with the eigenvalue $1$ of the dynamics matrix, which in general may be only approximately computed with polynomial-time algorithms. 
Subsequently, to overcome this computational demanding task, we proposed a second algorithm that requires the additional assumption that the network of agents is bidirected. 

The second algorithm designs four local augmented states of each agent and the weights of the adjacency matrix, but it further ensures that the designed final network is time-reversible. 
This additional property allows computing the left-eigenvector associated with the eigenvalue 1 of the dynamics matrix explicitly with linear-time complexity. 
We observe that we can increase the proposed methods robustness to attack using the general method proposed in~\cite{ramos2020ijoc}. 

Finally, we illustrated the proposed algorithms with examples. 
Future research directions include the design of the network parameters that ensures a faster convergence rate of the proposed methods.

% \begin{ack}% Place acknowledgements here.
% \end{ack}

% \appendices
% \section{Proof of the First Zonklar Equation}
% Appendix one text goes here.

% % you can choose not to have a title for an appendix
% % if you want by leaving the argument blank
% \section{}
% Appendix two text goes here.

% use section* for acknowledgment
% \section*{Acknowledgment}

% The authors would like to thank...

% Can use something like this to put references on a page
% by themselves when using endfloat and the captionsoff option.
% \ifCLASSOPTIONcaptionsoff
%   \newpage
% \fi

\bibliographystyle{IEEEtran}
\bibliography{autosam}   % and a bib file to produce the 

\end{document}